\documentclass[11pt]{article}
\usepackage[margin=1.15in]{geometry}
\usepackage{graphicx} 
\usepackage{soul,
amsmath,
amssymb,
amsthm,
color,
ytableau,
bm,
mathrsfs,
colortbl,
thmtools,
thm-restate,
enumerate,
mathtools,
bbm,
shuffle,
comment,
float,
tikz,
tabularx,
booktabs,
graphicx
}
\makeatletter
\@ifundefined{newcounteralias}{}{%
  \renewcommand\thmt@autorefsetup{%
    \@xa\def\csname\thmt@envname autorefname\@xa\endcsname\@xa{\thmt@thmname}%
  }%
}
\makeatother

\usetikzlibrary{
positioning, 
arrows.meta, 
matrix,
backgrounds
}

\definecolor{faceblue}{RGB}{198,219,239}
\definecolor{light-blue}{HTML}{e9f1f8}
\definecolor{grayish}{RGB}{230,230,230}
\definecolor{light-grayish}{RGB}{250,250,250}

\usepackage{hyperref}
\hypersetup{
    colorlinks=true,
    linkcolor=violet,
    filecolor=blue,      
    urlcolor=teal,
    citecolor=teal
}
\usepackage[nameinlink,capitalize]{cleveref}
\makeatletter
\@ifundefined{newcounteralias}{}{%
  \renewcommand\thmt@autorefsetup{%
    \@xa\def\csname\thmt@envname autorefname\@xa\endcsname\@xa{\thmt@thmname}%
  }%
}
\makeatother

\crefname{thmintro}{theorem}{theorems}
\Crefname{thmintro}{Theorem}{Theorems}
\Crefname{condition}{Condition}{Conditions}

\newtheorem{theorem}{Theorem}[section]
\newtheorem{lemma}[theorem]{Lemma}
\newtheorem{proposition}[theorem]{Proposition}
\newtheorem{corollary}[theorem]{Corollary}

\theoremstyle{definition}
\newtheorem{definition}[theorem]{Definition}

\newtheorem{example}[theorem]{Example}
\newtheorem{remark}[theorem]{Remark}

\theoremstyle{plain}

\renewcommand{\thethmintro}{\Alph{thmintro}}

\newcommand{\bC}{\mathbb{C}}
\newcommand{\bN}{\mathbb{N}}

\newcommand{\bK}{\mathbb{K}}

\newcommand{\sn}{S_n}

\newcommand{\fl}{\mathrm{Fl}}

\newcommand{\del}{\mathrm{del}}
\newcommand{\lk}{\mathrm{lk}}
\newcommand{\depth}{\mathrm{depth}}
\renewcommand{\star}{\mathrm{star}}
\newcommand{\hilb}{\mathrm{Hilb}}
\newcommand{\supp}{\mathrm{supp}}
\newcommand{\id}{\mathrm{id}}
\newcommand{\lcm}{\mathrm{lcm}}

\newcommand{\boldx}{\mathbf{x}}

\newcommand{\calW}{\mathcal{W}}
\newcommand{\calP}{\mathcal{P}}

\newcommand{\fraks}{\mathfrak{S}}

\newcommand{\Tor}{\operatorname{Tor}}
\newcommand{\Ann}{\operatorname{Ann}}

\DeclareMathOperator{\Sd}{Sd}
\DeclareMathOperator{\wt}{wt}

\newcommand{\tim}[1]{{\color{red}Tim: #1}}

\title{Face rings and the cohomology of unions of Schubert varieties}
\author{Ashleigh Adams-Jessmore, Sean T. Griffin, and Tim Ryan}
\date{\today}

\begin{document}
\maketitle
\begin{abstract}
The cohomology ring of a union of Schubert varieties $X_{\calW}=\bigcup_{w\in\calW}X_w,$ for $\calW\subseteq\sn$, in the complete flag variety is well known to be a quotient of the coinvariant ring $R_n$. 
We determine exactly when the ideal defining this quotient is generated by squarefree monomials, i.e.~when it can be identified with the image of a Stanley-Reisner ideal $I_\Delta$ in $R_n$. 
We prove that a single Schubert variety $X_w$ has this property exactly when $w\in\sn$ avoids a set of seven patterns, and in this case we describe the corresponding \emph{doubly-pruned} simplicial complex $\Delta$. We then show this is equivalent to $X_w$ being cut out of the complete flag variety by particular inclusion conditions on the pieces of the flag. 
We prove these complexes are vertex decomposable, and are Cohen--Macaulay precisely when $X_w$ is smooth. We then give a closed formula for the Poincar{\'e} polynomial of $X_w$.
For unions $X_\calW$, pattern avoidance is necessary but no longer sufficient, and we prove a necessary and sufficient condition we call \emph{prunable} that describes how the individual complexes associated to each $w\in \calW$ may be glued in the union. 
Quotients of $R_n$ by squarefree monomial ideals appear naturally in the Adams--Reiner conjecture on resolutions of Stanley-Reisner rings over certain parameter rings. We prove their conjecture for a family of simplicial complexes that includes the doubly-pruned complexes.
\end{abstract}
\section{Introduction}
\label{sec: introduction}
Write $S = \bK[x_1,\ldots,x_n]$ and $I_n = (e_1,\ldots,e_n)$ for the ideal of $S$ generated by the elementary symmetric polynomials. 
A classical result of Borel \cite{Borel1953} states that for any field $\bK,$ the cohomology ring of the complete flag variety $\fl_n$ has a presentation
\[
H^*(\fl_n;\bK)\cong R_n,
\]
where $R_n = S/I_n$ is the \emph{coinvariant algebra}. 
For any Schubert variety $X_w\subseteq \fl_n$, the restriction map $H^*(\fl_n;\bK)\twoheadrightarrow H^*(X_w;\bK)$ is surjective and its kernel is the ideal $J_w\subseteq R_n$ linearly spanned by the Schubert polynomials $\fraks_v$ with $v\not\leq w.$ 
Hence, we have a presentation of the cohomology ring of a Schubert variety as $H^*(X_w;\bK)\cong R_n/J_w.$

While the Schubert polynomials $\fraks_v$ for $v\not\leq w$ give a convenient $\bK$-linear spanning set of $J_w$ (modulo $I_n$), it is still an open problem to give a minimal generating set of the ideal $J_w$. Toward this goal, Reiner, Woo, and Yong \cite{ReinerWooYong} found a relatively small generating set of $J_w$ consisting of Schubert polynomials indexed by bigrassmannian permutations.
By work of St. Dizier-Yong \cite{stdizier-yong}, $J_w$ can be written as a sum of \emph{basic ideals}, and each basic ideal is minimally generated by a subset of Reiner--Woo--Yong generators.

On the other hand, Gasharov and Reiner gave a particularly small generating set of order $n^2$ for $J_w$ in the case that $w$ is \emph{defined by inclusions}, meaning that its associated Schubert variety $X_w$ is cut out of $\fl_n$ by inclusion conditions on the parts of its flag. 
Gasharov and Reiner characterize $w\in\sn$ defined by inclusions as those avoiding the patterns $4231,35142,42513,$ and $351624$ and show that $J_w$ is generated by the elementary symmetric polynomials in initial and final blocks of the variables. 
Ding \cite{ding01} has also given a family of Schubert varieties in partial flag varieties, each of whose cohomology ideal has a generating set with up to $n$ generators. 
Note, while there is some overlap, Ding's class does not contain the class of Schubert varieties discussed in this paper and vice versa. 

We fix notation before stating our results. 
Given a set of permutations $A$, denote by $S_n(A)$ the set of permutations in $S_n$ that avoid all patterns in $A$. 
Let
\begin{align*}
\calP &\coloneqq \{ 2143,\ 2341,\ 2413,\ 2431,\ 4123,\ 4213,\ 4231 \}.
\end{align*}
For $w\in S_n$, we sometimes say \emph{$w$ is $\calP$-avoiding} to mean $w\in S_n(\calP)$. 
For a simplicial complex $\Delta$ on a vertex set $[n]:=\{1,\ldots,n\},$ let $I_\Delta\subseteq S$ be it's \emph{Stanley-Reisner ideal}, which is generated by the squarefree monomials $x_G$ determined by the nonfaces $G\subseteq [n]$ of $\Delta.$ The corresponding quotient ring $R_{\Delta} = S/I_\Delta$ is called the \emph{Stanley-Reisner ring} of $\Delta$.

In this article, we introduce a family of simplicial complexes called \emph{doubly pruned} characterized by their minimal nonfaces being initial or final subsets of their support, see Definition~\ref{def: doubly pruned simplicial complexes} for the precise description. To each $w\in S_n(\calP)$, we attach a doubly pruned simplicial complex $\Delta_w$ and prove that the ideal $J_w\subseteq R_n$ is generated by square-free monomials. In fact, we prove these are exactly the Schubert varieties whose cohomology ring has this property. In our first main theorem, we prove that $J_w$ being generated by squarefree monomials in $R_n$ is equivalent to the pattern-avoidance condition above and to a set of inclusion conditions on the Schubert variety.

\begin{restatable}{thmintro}{MainA}
\label{thm:main-A}
    For $w\in S_n$, the following are equivalent:
    \begin{enumerate}
    \renewcommand{\theenumi}{\thethmintro\arabic{enumi}}
    \renewcommand{\labelenumi}{(\arabic{enumi})}
    \item\label[condition]{cond:A1} (Algebraic) The ideal $J_w$ is generated by squarefree monomials in $R_n$. 
    \item\label[condition]{cond:A2} (Combinatorial) $w\in\sn(\calP)$.
    \item\label[condition]{cond:A3}  (Geometric) The Schubert variety $X_w$ is defined by the following set \footnote{For $w=\id$, $X_\id$ is a single point defined by $V_k=\bC^k$ for all $k$.} of inclusions for some $i\leq j$
        \begin{itemize}
            \item $V_1 = \bC^1$,\dots, $V_{i-1} = \bC^{i-1}$,
            \item $\bC^i\subseteq V_{w^{-1}(i)}$,
            \item $V_{w^{-1}(j)-1}\subseteq \bC^{j-1}$,
            \item $V_{j} =\bC^{j}$,\dots, $V_{n-1} = \bC^{n-1}$
        \end{itemize}
    \end{enumerate}
\end{restatable}

Thus, if any of the equivalent conditions in \Cref{thm:main-A} holds, then $H^*(X_w;\bK)\cong S/J_w$ where $J_w$ is generated by squarefree monomials modulo $I_n$. Equivalently, $J_w = I_n +I_\Delta$ for some simplicial complex $\Delta$ on $[n]$. In this case, the set of generators of $J_w$ is particularly simple to describe and has size at most $n$.

In Section~\ref{sec: intro to doubly pruned simplicial complexes}, we give equivalent combinatorial conditions for a simplicial complex to be doubly pruned. Furthermore, the collection of doubly pruned complexes admits nice poset-theoretic properties. For example, in the assignment $w\in S_n(\calP)$ mapping to $\Delta_w$ gives an injection of the subposet of Bruhat order on $S_n(\calP)$ into the Boolean lattice, as we prove in \Cref{prop:Inclusion}. The image of this map of posets is properly contained in the collection of doubly pruned complexes. The subposet of doubly pruned complexes forms a sublattice of the Boolean lattice, as we prove in \Cref{cor: lattice}.

We extend \Cref{thm:main-A} to the setting of unions of Schubert varieties $X_\calW = \bigcup_{w\in \calW}X_w$. Since containment of Schubert varieties is governed by strong Bruhat order on their corresponding permutations, unions of Schubert varieties are indexed by antichain $\calW\subseteq S_n$ under Bruhat order. Set $J_\calW$ to be kernel of the surjective map $R_n = H^*(\fl_n;\bK)\to H^*(X_\calW;\bK)$. By standard arguments using dual bases, $J_\calW = \bigcap_{w\in \calW}J_w$. Further, set
\[
\Delta_\calW := \bigcup_{w\in\calW}\Delta_w.
\]
A basic property of Stanley--Reisner ideals is that the intersection of a collection of Stanley--Reisner ideals is equal to the Stanley--Reisner ideal of the union of the complexes, see~\cite{MillerSturmfels2005}. Thus, 
\[
I_{\Delta_\calW} = \bigcap_{w\in\calW}I_{\Delta_w}\subseteq S.
\]

In view of \Cref{thm:main-A}, one might expect that for $\calW\subseteq S_n(\calP)$,
\begin{equation}\label{eq: hopes and dreams}
J_\calW =\bigcap_{w\in \calW} J_w= \bigcap_{w\in \calW} (I_n + I_{\Delta_w}) =  I_n + \bigcap_{w\in \calW} I_{\Delta_w}= I_n + I_{\Delta_\calW},
\end{equation}
and thus $J_\calW$ is generated by squarefree monomials (modulo $I_n$). However, the third equality in \eqref{eq: hopes and dreams} does not always hold, even in the case of $|\calW|=2$, as the next example shows.

\begin{example} 
Let $n = 4$ and the commuting simple reflections $u = s_1 = 2134$ and $v = s_3 = 1243,$ so that $X_u$ and $X_v$ are projective curves meeting at $X_\id.$ Both lie in $\mathrm{S}_4(\calP),$ so $J_u$ and $J_v$ are generated by squarefree monomials, while $J_\calW$ is not. Indeed, since $s_2\not\leq s_1$ and $s_2\not\leq s_3,$ the class $\fraks_{s_2} = x_1 + x_2$ lies in $J_\calW,$ while no single variable is in $J_\calW$, as the reader can check. The failure of \eqref{eq: hopes and dreams} can then be seen as follows. The complex $\Delta_\calW$ is four isolated points and $S/(I_{\Delta_\calW}+I_4)$ has Hilbert series $1+3t,$ whereas for $X_\calW,$ we have two curves meeting at a point, which has cells indexed by $\id, s_1,$ and $s_3,$ and thus has Poincar\'e polynomial $1+2t.$  
\end{example}

The failure of \eqref{eq: hopes and dreams} is due to the fact that, in general, 
\[
I_n + (I_{\Delta_u}\cap I_{\Delta_v}) \subsetneq (I_n + I_{\Delta_u}) \cap (I_n + I_{\Delta_v}).
\]
We translate this failure into homological terms, as follows, and prove that it is equivalent to a topological condition on the associated complexes. Since $I_{\Delta_u}\cap I_{\Delta_v} = I _{\Delta_u\cup\Delta_v}$ and $I_{\Delta_u}+I_{\Delta_v} = I_{\Delta_u\cap\Delta_v},$ by lifting to the face rings, one may construct the Mayer-Vietoris sequence
\[
0\longrightarrow
S/I_{\Delta_u\cup\Delta_v}
\longrightarrow
S/I_{\Delta_u}\oplus S/I_{\Delta_v}
\longrightarrow
S/I_{\Delta_u\cap\Delta_v}
\longrightarrow
0.
\]
The cohomology rings of $X_u,X_v,$ and their union and intersection fit into an analogous short exact sequence, since Schubert varieties have no odd degree cohomology.
However, passing from face rings to cohomology rings means tensoring with $R_n = S/I_n,$ which is only right exact. The failure above is exactly the cokernel of the induced map on $\Tor_1^S(R_n,-)$, which we explain in \Cref{subsec:Prune_Equality}.

Our second main theorem (\Cref{thm:main-B}) shows that this failure is the only obstruction. More importantly, it shows that the only time that \Cref{eq: hopes and dreams} holds is for a Bruhat antichain $\calW\subseteq\sn$ such that $J_\calW$ is generated by squarefree monomials in $R_n.$ We also give a combinatorial-geometric interpretation.  
In order to prove this result (\Cref{sec: proofs of the main theorems}), we use Koszul resolutions to show that the map induced by $\Tor_1^S(R_n,-)$ is surjective if and only if the complexes $\Delta_w$ are glued together in a way reminiscent of a shelling on simplicial complexes (\Cref{def:prunable}).

\begin{restatable}[]{thmintro}{MainB}
\label{thm:main-B}
Let \(\mathcal W\subseteq S_n\) be a strong-Bruhat antichain. Then the following are equivalent.

\begin{enumerate}[(1)]
    \renewcommand{\theenumi}{\thethmintro\arabic{enumi}}
    \renewcommand{\labelenumi}{(\arabic{enumi})}
\item\label[condition]{cond:B1} (Algebraic) The ideal $J_{\mathcal W}$ is generated by squarefree monomials in \(R_n\) where
\[ J_{\mathcal W}:=\bigcap_{w\in\mathcal W} J_w\subseteq R_n. \]

\item\label[condition]{cond:B2} (Cohomological) We have \(\mathcal W\subseteq S_n(\calP)\) and 
\[ \bigcap_{w\in\mathcal W}(I_n+I_{\Delta_w}) = I_n+\bigcap_{w\in\mathcal W}I_{\Delta_w}\]
as ideals in $S=\bK[x_1,\dots,x_n]$.

\item\label[condition]{cond:B3} (Combinatorial) We have $\mathcal{W}\subseteq S_n(\mathcal{P})$, and the family of simplicial complexes $\{\Delta_w:w\in \mathcal{W}\}$ is \textbf{prunable}, which we define below in Definition~\ref{def:prunable}.\label[(3)]{item:3}
\end{enumerate}
\end{restatable}


The topological property we need in \eqref{item:3} above is as follows. Note that it is \textbf{not} a purely topological property of the union $\bigcup_{w\in \mathcal W}\Delta_w$ but rather depends on how the $\Delta_w$ are glued together.

\begin{definition}\label{def:prunable}
    Given a Bruhat antichain $\mathcal{W}\subseteq S_n(\mathcal{P})$, we say $\{\Delta_w:w\in \mathcal{W}\}$ is \textbf{prunable} if \(\mathcal W\) admits an ordering $w_1,\ldots,w_k$ such that the following two conditions hold for every \(2\le r\le k\), where
    \[\Omega_r\coloneqq\Delta_{w_r}\cap(\Delta_{w_1}\cup\cdots\cup\Delta_{w_{r-1}}).\]
    \begin{enumerate}
    \item[(a)]\(\Omega_r\) is pure and $\dim\Omega_r = \dim\Delta_{w_r}$.
    \item[(b)] For every nonempty subset \(T\subseteq\{1,\ldots,r-1\}\), there is some \(t\in T\) such that $$\Delta_{w_r}\cap \left(\bigcup_{i\in T}\Delta_{w_i}\right) = \Delta_{w_r}\cap \Delta_{w_t}.$$ 
    \end{enumerate}
\end{definition}

\begin{example}
In the case of $|\calW| = 2$, $\calW = \{v,w\}$ is prunable if and only if the intersection $\Delta_v\cap \Delta_w$ is pure of dimension $\min(\dim\Delta_v,\dim\Delta_w)$. As an example of how \emph{prunable} is not solely a topological condition on the union $\Delta_v\cup \Delta_w$, observe that for $v = 3412$, $w=3241$, $\Delta_v\cap\Delta_w$ is not pure dimensional, so $\{v,w\}$ is not prunable. However, $\Delta_v\cup \Delta_w = \Delta_{3421}$, which is a doubly pruned complex, hence trivially prunable.

\[
\begin{tikzpicture}[scale=0.56, x=1cm, y=1cm,
  vtx/.style={circle, fill=black, inner sep=1.15pt},
  ghost/.style={circle, draw=black!30, fill=white, inner sep=1.0pt},
  ed/.style={line width=0.7pt, black!80}]
  \coordinate (v2) at (0.000,1.300);
  \coordinate (v3) at (-1.220,-0.780);
  \coordinate (v4) at (1.220,-0.780);
  \coordinate (v1) at (0.000,-0.060);
  \fill[faceblue] (v1) -- (v2) -- (v4) -- cycle;
  \fill[faceblue] (v1) -- (v3) -- (v4) -- cycle;
  \draw[ed] (v1) -- (v2);
  \draw[ed] (v1) -- (v3);
  \draw[ed] (v1) -- (v4);
  \draw[ed] (v2) -- (v3);
  \draw[ed] (v2) -- (v4);
  \draw[ed] (v3) -- (v4);
  \node[vtx] at (v1) {};
  \node[above right, inner sep=1.6pt, font=\scriptsize] at (v1) {$1$};
  \node[vtx] at (v2) {};
  \node[above, inner sep=1.6pt, font=\scriptsize] at (v2) {$2$};
  \node[vtx] at (v3) {};
  \node[below left, inner sep=1.6pt, font=\scriptsize] at (v3) {$3$};
  \node[vtx] at (v4) {};
  \node[below right, inner sep=1.6pt, font=\scriptsize] at (v4) {$4$};
  \node at (0,-2) {$\Delta_{3412}$};
\end{tikzpicture}
\qquad
\begin{tikzpicture}[scale=0.56, x=1cm, y=1cm,
  vtx/.style={circle, fill=black, inner sep=1.15pt},
  ghost/.style={circle, draw=black!30, fill=white, inner sep=1.0pt},
  ed/.style={line width=0.7pt, black!80}]
  \coordinate (v2) at (0.000,1.300);
  \coordinate (v3) at (-1.220,-0.780);
  \coordinate (v4) at (1.220,-0.780);
  \coordinate (v1) at (0.000,-0.060);
  \fill[faceblue] (v1) -- (v2) -- (v3) -- cycle;
  \fill[faceblue] (v1) -- (v2) -- (v4) -- cycle;
  \draw[ed] (v1) -- (v2);
  \draw[ed] (v1) -- (v3);
  \draw[ed] (v1) -- (v4);
  \draw[ed] (v2) -- (v3);
  \draw[ed] (v2) -- (v4);
  \node[vtx] at (v1) {};
  \node[above right, inner sep=1.6pt, font=\scriptsize] at (v1) {$1$};
  \node[vtx] at (v2) {};
  \node[above, inner sep=1.6pt, font=\scriptsize] at (v2) {$2$};
  \node[vtx] at (v3) {};
  \node[below left, inner sep=1.6pt, font=\scriptsize] at (v3) {$3$};
  \node[vtx] at (v4) {};
  \node[below right, inner sep=1.6pt, font=\scriptsize] at (v4) {$4$};
    \node at (0,-2) {$\Delta_{3241}$};
\end{tikzpicture}
\quad
\begin{tikzpicture}[scale=0.56, x=1cm, y=1cm,
  vtx/.style={circle, fill=black, inner sep=1.15pt},
  ghost/.style={circle, draw=black!30, fill=white, inner sep=1.0pt},
  ed/.style={line width=0.7pt, black!80}]
  \coordinate (v2) at (0.000,1.300);
  \coordinate (v3) at (-1.220,-0.780);
  \coordinate (v4) at (1.220,-0.780);
  \coordinate (v1) at (0.000,-0.060);
  \fill[faceblue] (v1) -- (v2) -- (v4) -- cycle;
  \draw[ed] (v1) -- (v2);
  \draw[ed] (v1) -- (v3);
  \draw[ed] (v1) -- (v4);
  \draw[ed] (v2) -- (v3);
  \draw[ed] (v2) -- (v4);
  \node[vtx] at (v1) {};
  \node[above right, inner sep=1.6pt, font=\scriptsize] at (v1) {$1$};
  \node[vtx] at (v2) {};
  \node[above, inner sep=1.6pt, font=\scriptsize] at (v2) {$2$};
  \node[vtx] at (v3) {};
  \node[below left, inner sep=1.6pt, font=\scriptsize] at (v3) {$3$};
  \node[vtx] at (v4) {};
  \node[below right, inner sep=1.6pt, font=\scriptsize] at (v4) {$4$};
    \node at (0,-2) {$\Delta_{3412}\cap\Delta_{3241}$};
\end{tikzpicture}
\qquad
\begin{tikzpicture}[scale=0.56, x=1cm, y=1cm,
  vtx/.style={circle, fill=black, inner sep=1.15pt},
  ghost/.style={circle, draw=black!30, fill=white, inner sep=1.0pt},
  ed/.style={line width=0.7pt, black!80}]
  \coordinate (v2) at (0.000,1.300);
  \coordinate (v3) at (-1.220,-0.780);
  \coordinate (v4) at (1.220,-0.780);
  \coordinate (v1) at (0.000,-0.060);
  \fill[faceblue] (v1) -- (v2) -- (v3) -- cycle;
  \fill[faceblue] (v1) -- (v2) -- (v4) -- cycle;
  \fill[faceblue] (v1) -- (v3) -- (v4) -- cycle;
  \draw[ed] (v1) -- (v2);
  \draw[ed] (v1) -- (v3);
  \draw[ed] (v1) -- (v4);
  \draw[ed] (v2) -- (v3);
  \draw[ed] (v2) -- (v4);
  \draw[ed] (v3) -- (v4);
  \node[vtx] at (v1) {};
  \node[above right, inner sep=1.6pt, font=\scriptsize] at (v1) {$1$};
  \node[vtx] at (v2) {};
  \node[above, inner sep=1.6pt, font=\scriptsize] at (v2) {$2$};
  \node[vtx] at (v3) {};
  \node[below left, inner sep=1.6pt, font=\scriptsize] at (v3) {$3$};
  \node[vtx] at (v4) {};
  \node[below right, inner sep=1.6pt, font=\scriptsize] at (v4) {$4$};
    \node at (0,-2) {$\Delta_{3421}$};
\end{tikzpicture}
\]
\end{example}

Along the way, we prove that each simplicial complex $\Delta_w$ for $w\in S_n(\mathcal{P})$ has the following topological properties.

\begin{restatable}[]{thmintro}{MainC}
\label{thm:main-C}
    For each $w\in\sn(\calP),$ the simplicial complex $\Delta_w$ is vertex decomposable, hence shellable and sequentially Cohen-Macaulay. When $\Delta_w$ is pure, it is Cohen-Macaulay. 
\end{restatable}

Moreover, we give simple conditions to determine the purity of a simplicial complex $\Delta_w$ for $w\in S_n(\mathcal P)$.

In \Cref{sec: the combinatorics of doubly pruned simplicial complexes}, we analyze the homological properties of the cohomology ring $R_n/J_w \cong H^*(X_w;\bK)$ and the Stanley--Reisner ring $R_{\Delta_w} \coloneqq S/I_{\Delta_w}$ when $w\in S_n(\calP)$. In particular, we compute the minimal free resolution, the depth and dimension, and the Hilbert series of these rings. The bulk of the work occurs in the case when the squarefree monomial relations ``cross.'' Furthermore, since the elementary symmetric polynomials do not form a regular sequence on $R_{\Delta_w}$, the Hilbert series cannot readily be read off of a Koszul complex computation.

The same failure of freeness is the focus of our last result, which is a more general statement than simply doubly pruned complexes, but does include them as a specialization. After passing to $R_\Delta,$ the relations $I_n$ become canonical. The image of $e_k$ is the sum of $x_F$ over the $k$-element faces of $\Delta,$ the $k$th universal parameter of $\Delta$ denoted $\theta_k$ in the sense of Herzog-Moradi \cite{herzog-moradi} and Adams-Reiner \cite{adams-reiner}, so that $R_{\Delta,n} = R_{\Delta}/(\theta),$ where $\theta = (\theta_1,\ldots,\theta_d).$ Adams and Reiner conjecture \cite[Conj.~6.1]{adams-reiner} that the minimal free resolution of $R_\Delta$ over $\bK[\theta]$ is governed by the barycentric subdivision of $\Delta,$
\[
\dim_\bK
\Tor_r^{\bK[\theta]}(R_\Delta,\bK)_N
=
\sum_{
\substack{
T\subseteq[d] \\
\wt T = N
}
}
\dim_\bK\tilde{H}_{\#T-r-1}(\Sd_T(\Delta);\bK),
\]
where $d = \dim\Delta +1, \wt T = \sum_{k\in T}k,$ and $\Sd_T(\Delta)$ is the rank-selected subcomplex. The conjecture is open in general.

\begin{restatable}[]{thmintro}{MainD}
\label{thm:main-D}
        The Adams-Reiner conjecture (restated in \Cref{eq: Adams-Reiner conj}) holds over any field for every simplicial complex with at most two minimal nonfaces that are of cardinality at least two. In particular, it holds for all doubly pruned complexes, and so for every $\Delta_w$ with $w\in\sn(\calP).$ Moreover, $\Tor_r^{\bK[\theta]}(R_\Delta,\bK) = 0$ for $r\geq 2,$ and $\Tor_1$ vanishes unless the two nonfaces meet. In the latter case,
    \[
    \sum_N
    \dim_\bK \Tor_1^{\bK[\theta]}(R_\Delta,\bK)_N t^N 
    = t^{m-1+\gamma}~[\alpha]_t[\beta]_t[m-2]_t!.
    \]
\end{restatable}

 For $w\in\sn(\calP),$ $\Tor_1$ measures the failure of $H^*(X_w;\bK)$ to be free over $\bK[\theta],$ is nonzero precisely in the crossing case, and has total dimension $(q-i)(j-p)(j-i-1)!.$ Crossing is therefore the failure of Cohen-Macaulayness, the obstruction to the Koszul computation, and the only source of higher $\Tor$ over the parameter ring.

\textbf{Organization.}
This paper is organized as follows. In \Cref{sec: background}, we collect relevant background. In \Cref{sec: ideals in coinvariant rings}, we prove several lemmas on squarefree monomials in $R_n.$ In \Cref{sec: intro to doubly pruned simplicial complexes}, we introduce doubly pruned simplicial complexes and study their properties, including characterizing them by pattern avoidance. In \Cref{sec: proofs of the main theorems}, we prove \Cref{thm:main-A} and \Cref{thm:main-B}, the latter using Koszul cycle arguments. In \Cref{sec: the combinatorics of doubly pruned simplicial complexes}, we prove \Cref{thm:main-C} and \Cref{thm:main-D} and compute the Hilbert series.  In \Cref{sec: example}, we work out the example case where $n=4.$ 

\textbf{Future work.} Squarefree monomials have been studied previously in Schubert geometry in the context of degenerations of matrix Schubert varieties and their coordinate rings, see \cite{bertiger,brion-lakshmibai, hamaker-pechenik-weigandt,klein-diagonal,knutson-lam-speyer-richardson, knutson-miller-grobner}. In contrast, this article characterizes when Schubert cohomology ideals are generated by squarefree monomials. It would be interesting to find parallels between the simplicial complexes in that setting and the doubly pruned complexes introduced in our work.

Product formulas for Poincar\'e polynomials have precedents in Gasharov \cite{gasharov-factoring} and in the fibre bundle structures of Schubert varieties \cite{alland-richmond,richmond-slofstra}, whereas ours come from a free resolution. Closest to Section~6 are the complexes built from lower Bruhat intervals by Bolognini-Sentinelli \cite{bolognini-sentinelli}.

Finally, one question these results raise is whether \Cref{thm:main-B} holds $T$-equivariantly, with squarefree products of the linear forms $x_r-y_s$ in place of squarefree monomials. One implication is immediate upon specializing $y=0,$ and Stanley-Reisner presentations of equivariant cohomology are known in the GKM setting \cite{goertsches-solomadin}, of which $\fl_n$ is an instance.

\section*{Acknowledgments}
Adams was supported by Japan Science and Technology Agency (JST), Excellence Program for Engaging Research Talent – Japan (EXPERT-J), Grant Number JPMJEX2501, NSF grant DMS-2247089, and NSF Graduate Fellowship Grant No. 2034612.
\section*{AI disclosure}
The mathematical content of this work, including all statements and proofs, is due to the authors. Two generative AI tools (Claude: Opus 5.5, OpenAI: GPT 5.5) were used for exploration, proofreading, for secondary computer verification of many of the results for $n\leq 6$ (the second author, Griffin, also wrote SageMath code for primary verification and pattern finding), for cleaning up \Cref{fig: s4 simplicial complex poset}, and for literature searches.
\section{Background}
\label{sec: background}
\subsection{Simplicial complexes and Stanley-Reisner rings}
\label{subsec: simplicial complexes and Stanley-Reisner rings}
This section will outline the theory of simplicial complexes and Stanley-Reisner theory.
\begin{definition}
    A (finite, regular) \emph{simplicial complex} $\Delta$ on a vertex set $V = [n] :=\{1,\ldots,n\}$ is a collection of subsets $F\subseteq V$ called \emph{faces} such that if $F$ is a face of $\Delta$ and $G\subseteq F,$ then $G\in\Delta.$ The \emph{dimension} of a face $F$ is $\dim F = |F| - 1.$
\end{definition}
\begin{remark}
    We note that traditionally the definition of simplicial complexes also includes the stipulation that $\Delta$ is supported on $V$, meaning that every singleton $\{v\}$ for $v\in V$ is a face of $\Delta.$ 
    If $v\in V,$ but $\{v\}$ is not a face of $\Delta,$ then there exists no face of $\Delta$ containing $\{v\}.$ Thus, $\Delta$ is simply a simplicial complex on the vertex set $V\setminus\{v\}.$ Algebraically, this corresponds to including the variable $x_v$ in the Stanley-Reisner ideal of $\Delta$ (Definition~\ref{def: SR ideal}). 
    In this article, we allow singletons to be nonfaces in our definition of a simplicial complex. 
\end{remark}
We now give a series of short definitions required for working with simplicial complexes. For further background on simplicial complexes and Stanley-Reisner theory, we refer the reader to \cite{BrunsHerzog1993, MillerSturmfels2005,Stanley1996}. A \emph{nonface} is a subset $G\subseteq V$ such that $G$ is not a face in $\Delta$ and it is a \emph{minimal nonface} if there does not exist a nonface $H\subsetneq G$.  
A \emph{facet}, or \emph{maximal face}, is a face in $\Delta$ that is not properly contained in any other face of $\Delta$. The \emph{dimension} of a face $F$ is $\dim F = |F|-1$, and
\[
\dim\Delta = \max\{\dim F : F\text{ is a facet of }\Delta\}.
\]
If all facets of $\Delta$ have the same dimension, then $\Delta$ is called \emph{pure}.

If every subset of $V$ is a face of $\Delta,$ equivalently $\Delta$ has no nonfaces (other than vertex nonfaces), then $\Delta$ is called a \emph{simplex}. In this case, $V$ is the unique maximal facet of $\Delta$ and $\dim\Delta = |V|-1.$
A collection of faces $\Delta'\subseteq\Delta$ that is closed under subsets is also a simplicial complex and is called a \emph{subcomplex} of $\Delta.$ Given any face $F$ of a simplicial complex $\Delta,$ the collection $\overline{F} = \{G:G\subseteq F\}$ of all subsets of $F$ is closed under subsets, so $\overline{F}$ is also a subcomplex of $\Delta.$ In particular, it is the simplex on the vertices of $F.$ 
Two other important subcomplexes are the deletion of a vertex and link of a face.
\begin{definition}
\label{def: deletion and link of a simplicial complex}
Given a simplicial complex $\Delta$ on a vertex set $V$ and for any $v\in V$ the \emph{deletion of $v$} is 
\[
\del_v(\Delta) = \{F\in\Delta: v\not\in F\}
\]
and the \emph{link of a face $F\in\Delta$} is
\[
\lk_F(\Delta) = \{G\in\Delta : G\cap F = \varnothing \text{ and }G\cup F\in\Delta\}.
\]    
When $F = \{v\}$ is a singleton we write $\lk_v(\Delta)$ for $\lk_{\{v\}}(\Delta).$
\end{definition}
\begin{definition}[\cite{bjorner1996shellable,provan1980decompositions}]
\label{def: vertex-decomposable}
    A (finite, regular) simplicial complex $\Delta$ on $n$ vertices is \emph{vertex decomposable} if
    \begin{itemize}
        \item $\Delta$ is a simplex; or
        \item there exists a vertex $v\in V$ such that
        \begin{enumerate}
            \item $\del_v(\Delta)$ is vertex decomposable,
            \item $\lk_v(\Delta)$ is vertex decomposable, and 
            \item $v$ is a \emph{shedding vertex}, meaning no facet of $\lk_v(\Delta)$ is a facet of $\del_v(\Delta).$
        \end{enumerate}
    \end{itemize}
\end{definition}
If $\Delta$ is vertex decomposable, then it is also shellable \cite{bjorner1996shellable}. 
\begin{definition}[\cite{bjorner1996shellable}]
\label{def: shellable}
    A finite simplicial complex $\Delta$ is \emph{shellable} if its facets can be arranged in a linear order $F_1,F_2,\ldots,F_t$ such that for all $2\leq j\leq t$, the subcomplex
    \[
    \left(
    \bigcup_{i=1}^{j-1}\overline{F_i}
    \right)
    \cap
    \overline{F_j}
    \]
    is pure of dimension $\dim F_j - 1$. Such an order is called a \emph{shelling} of $\Delta$.
\end{definition}
We move from the combinatorics of simplicial complexes into the algebra of these topological objects. Let $\bK$ be a field and let $S = \bK[x_1,\ldots,x_n]$ be a polynomial ring. 
\begin{definition}
\label{def: SR ideal}
If $\Delta$ is a simplicial complex on $n$ vertices then define the \emph{Stanley-Reisner ideal} of $\Delta$ to be
\[
I_\Delta = (\boldx_G~:~G\text{ is a nonface of }\Delta),
\]
where $\boldx_G = \prod_{v\in G}x_v,$ and the \emph{Stanley-Reisner ring} to be the quotient $R_\Delta := S/I_\Delta.$
\end{definition}
We observe that the Stanley-Reisner ideal is generated by the minimal nonfaces of $\Delta,$ since if $G$ is a nonface of $\Delta$ that is not minimal, then $G\supsetneq H$ for some minimal nonface $H\neq\varnothing,$ and $\boldx_H$ divides $\boldx_G.$
The Stanley-Reisner ring is a graded ring of Krull dimension $\dim\Delta + 1.$ 
\begin{example}
\label{ex: Stanley-Reisner ring - simplex}
    The Stanley-Reisner ring of the simplex $\Delta$ on $n$ vertices (often written as $\Delta_n$) is the polynomial ring $\bK[x_1,\ldots,x_n].$
\end{example}
 
There is a natural way to translate unions and intersections of simplicial complexes into their Stanley-Reisner rings. Given simplicial complexes $\Delta_1$ and $\Delta_2$ on the same vertex set $V = [n],$ with Stanley-Reisner rings $R_{\Delta_1}$ and $R_{\Delta_2},$ then 
\begin{itemize}
    \item $I_{\Delta_1} + I_{\Delta_2} = I_{\Delta_1\cap\Delta_2}$
    \item $I_{\Delta_1}\cap I_{\Delta_2} = I_{\Delta_1\cup\Delta_2}.$
\end{itemize}
Let $\Delta$ be a simplicial complex of dimension $d$ on a vertex set $V = [n]$
with Stanley-Reisner ring $R_\Delta = S/I_\Delta$, where $S = \bK[x_1,\ldots,x_n]$ and $I_\Delta$ is the Stanley-Reisner ideal of $\Delta$. Then $R_\Delta$ is a finitely generated graded $S$-module, and each graded piece $(R_\Delta)_a$, $a\in\bN$, is a finite-dimensional $\bK$-vector space. Recording these dimensions in a generating function gives the \emph{Hilbert series}
\[
\hilb_S(R_\Delta;t) = \sum_{a\geq 0}\dim_\bK\!\left((R_\Delta)_a\right)t^a.
\]
We may rewrite the Hilbert series in the closed form
\[
\hilb_S(R_\Delta;t) = \frac{h(t)}{(1-t)^{d+1}},
\]
where $h(t)$ is called the \emph{h-polynomial}.
\begin{definition}[\cite{BrunsHerzog1993}]   
     A simplicial complex $\Delta$ is \emph{Cohen-Macaulay over $\bK$} if its corresponding Stanley-Reisner ring is Cohen-Macaulay.
\end{definition}
Passing to the ring to verify Cohen-Macaulayness is not always the most effective check. Instead, we may use Reisner's criterion, which allows one to look only at the simplicial complex itself to determine whether it is Cohen-Macaulay or not.
\begin{theorem}[\cite{Reisner1976}]
    Let $\Delta$ be a simplicial complex and let $\bK$ be a field. Then $\Delta$ is Cohen-Macaulay over $\bK$ if and only if for every face $F\in\Delta,$ $\Tilde{H}_i(\lk_F(\Delta);\bK) = 0$ for all $i<\dim\lk_F(\Delta).$ This is called \emph{Reisner's criterion}.
\end{theorem}
\begin{definition}[{\cite[Ch.~III]{Stanley1996}}]
\label{def: sequentially CM}
    A finitely generated graded $R_\Delta$-module $M$ is \emph{sequentially Cohen--Macaulay}
    if there is a filtration of graded submodules
    $0 = M_0 \subsetneq M_1 \subsetneq \cdots \subsetneq M_r = M$
    such that each quotient $M_i/M_{i-1}$ is Cohen--Macaulay and
    $\dim(M_1/M_0) < \dim(M_2/M_1) < \cdots < \dim(M_r/M_{r-1})$,
    where $\dim$ denotes Krull dimension.
    A simplicial complex $\Delta$ is sequentially Cohen--Macaulay over $\bK$ if $R_\Delta$ is sequentially Cohen--Macaulay as a module over itself.
\end{definition}
If $\Delta$ is shellable, then it is also sequentially Cohen--Macaulay
\cite[Thm.~11.3]{BjornerWachsII}. If $\Delta$ is both pure and shellable, then it is
Cohen--Macaulay \cite[Ch.~III]{Stanley1996}.
Summarizing the previous implications:
\[
\text{vertex decomposable}~\implies~\text{shellable}~\implies~\text{sequentially Cohen-Macaulay}
\]
and, if the simplicial complex is pure, shellable $\implies$ Cohen-Macaulay \cite{bjorner1996shellable,Stanley1996}.

\subsection{Flag varieties and the coinvariant algebra}
\label{subsec: Flag varieties and the coinvariant algebra}
We now review the combinatorial and geometric background required following \cite{BjornerBrenti} and \cite{FultonYoungTableaux}. Throughout this section, $\bK$ denotes an arbitrary field.

Let $\sn$ denote the symmetric group on $[n]:=\{1,\ldots,n\}$ written in one-line notation and generated by the simple transpositions we denote $s_i = (i~i+1)$ for $1\leq i\leq n-1.$
The \emph{length} $\ell(w)$ of a permutation $w\in\sn$ is the number of inversions of $w,$ that is the number of pairs $(a~b)$ with $a < b$ and $w(a)>w(b).$ The element $w_0 = n(n-1)\cdots 1$ is the unique longest element in $\sn.$ In general, we write $\leq$ for the strong Bruhat order on $\sn,$ which is characterized by $v\leq w$ if and only if some reduced word for $v$ is a subword of some reduced word for $w$, see for example \cite{BjornerBrenti}. Alternatively, we may consider the Bruhat order in an equivalent numerical form.
Throughout, we shall encode $w\in \sn$ as its permutation matrix with a 1 in matrix entry $(w(j),j)$ (row $w(j)$ and column $j$) for each $j$ and 0s elsewhere.
\begin{definition}
    Let $w\in\sn.$ The \emph{rank matrix} of $w$ is the $n\times n $ matrix $r_w$ whose entry in row $a$ and column $b$ is
    \[
    r_w(a,b) = \#\{k\leq b : w(k)\leq a\}.
    \]
    In other words, $r_w(a,b)$ is the number of 1s in the upper-left submatrix of $w$ formed by its top $a$ rows and first $b$ columns.
\end{definition}
Using Ehresmann's criterion \cite{ehresmann}, we are able to put an ordering on the rank matrices that matches the strong Bruhat order on their associated permutations.
\begin{proposition}\cite{BjornerBrenti}
\label{prop: Ehresmann's criterion}
    Let $v,w\in\sn.$ Then $v\leq w$ if and only if $r_v(a,b)\geq r_w(a,b)$ for all $1\leq a,b\leq n.$
\end{proposition}

Let $\fl_n$ denote the variety of complete flags
\[
0 = V_0\subset V_1\subset \cdots\subset V_n = \bC^n,\quad \dim_\bC V_i = i.
\]
Fixing the standard basis $\epsilon_1,\ldots,\epsilon_n$ of $\bC^n,$ we write $E_\bullet$ for the standard flag and $F_\bullet$ for the opposite flag with $F_i = \langle \epsilon_n,\ldots,\epsilon_{n-i+1}\rangle.$ Identifying $\fl_n$ with $GL_n(\bC)/B,$ where $B$ is the Borel subgroup (i.e., the subgroup of invertible upper triangular matrices), the Bruhat decomposition expresses $\fl_n$ as a disjoint union of $B$-orbits
\[
\fl_n = \bigsqcup_{w\in\sn}BwB/B.
\]
We write $C_w = BwB/B$ for the \emph{Schubert cell} of $w$ and $X_w = \overline{C_w}$ for the corresponding \emph{Schubert variety.} It is known that $\dim_\bC X_w = \ell(w),$ where $\ell(w)$ is the length of $w.$
The Bruhat order records containment of Schubert varieties and is given by the rule: $X_v\subseteq X_w$ if and only if $v\leq w$ \cite{FultonYoungTableaux}. Combined with \Cref{prop: Ehresmann's criterion}, we can translate containment of Schubert varieties into a finite list of inequalities among rank matrix entries. 

Given $w$, the Schubert variety $X_w$ is cut out of $\fl_n$ by a list of rank conditions on its parts. First, the \emph{Rothe diagram} of $w$ is
\[
D(w)\coloneqq \{(i,j)\in [n]\times [n] : j< w^{-1}(i)\text{ and }w(j) < i\}
\]
We note that this convention for $D(w)$ is not the standard one but is adapted well for Schubert varieties as we have defined them above. 
Pictorially, $D(w)$ is the set of cells in the $n\times n$ grid remaining after eliminating all cells immediately North and immediately East of each $1$ in $w$ (including the cell containing each $1$). The \emph{essential set} \cite{FultonYoungTableaux} of $w$ is
\[
E(w) \coloneqq \{(a,b) \in D(w): (a-1,b)\notin D(w)\text{ and }(a,b+1)\notin D(w)\}.
\]
In other words, $E(w)$ is the set of cells of $D(w)$ that are maximally northwest in its connected component of $D(w)$. 
The rank conditions defining $X_w$ are then
\[
X_w = \{V_\bullet\in \fl_n : \dim_\bC(V_b\cap E_a) \geq r_w(a,b)\text{ for }(a+1,b)\in E(w)\}
\]

\begin{example}\label{ex:Essential}
For $w = 1657423$, the permutation matrix, rank matrix, Rothe diagram (shaded purple), and essential set (highlighted with a black border) are
\[
\begin{bmatrix}
    1 & 0 & 0 & 0 & 0 & 0 & 0\\
    0 & 0 & 0 & 0 & 0 & 1 & 0\\
    0 & 0 & 0 & 0 & 0 & 0 & 1\\
    0 & 0 & 0 & 0 & 1 & 0 & 0\\
    0 & 0 & 1 & 0 & 0 & 0 & 0\\
    0 & 1 & 0 & 0 & 0 & 0 & 0\\
    0 & 0 & 0 & 1 & 0 & 0 & 0\\
\end{bmatrix}\quad
\begin{bmatrix}
    1 & 1 & 1 & 1 & 1 & 1 & 1\\
    1 & 1 & 1 & 1 & 1 & 2 & 2\\
    1 & 1 & 1 & 1 & 1 & 2 & 3\\
    1 & 1 & 1 & 1 & 2 & 3 & 4\\
    1 & 1 & 2 & 2 & 3 & 4 & 5\\
    1 & 2 & 3 & 3 & 4 & 5 & 6\\
    1 & 2 & 3 & 4 & 5 & 6 & 7\\
\end{bmatrix}\quad
\begin{tikzpicture}[scale=0.46,baseline={([yshift=-2.9mm]current bounding box.center)}]
\foreach \c/\r in {
    1/2,
    1/3, 6/3,
    1/4,
    1/5,
    1/6,
    1/7, 2/7, 3/7
}{
    \fill[blue!30] (\c-1,7-\r) rectangle (\c,8-\r);
}
\draw[step=1,black,thin] (0,0) grid (7,7);
\foreach \c/\r in {
    1/1,
    2/6,
    3/5,
    4/7,
    5/4,
    6/2,
    7/3
}{
    \fill[black] (\c-0.5,7.5-\r) circle (0.13);
}
\foreach \i in {1,...,7}{
    \node[above] at (\i-0.5,7) {\small \i};
}
\foreach \i in {1,...,7}{
    \node[left] at (0,7.5-\i) {\small \i};
}
\foreach \c/\r in {1/2, 3/7, 6/3} {\draw[black, line width=1.5pt] (\c-1,7-\r) rectangle (\c,8-\r);}
\end{tikzpicture}
\]
Thus, $X_w$ is defined by the rank conditions $\dim_\bC(V_1\cap E_1)\geq 1$, $\dim_\bC(V_3\cap E_6)\geq 3$, and $\dim_\bC(V_6\cap E_2)\geq 2$. Equivalently, $V_1 = E_1$, $V_3 \subset E_6$, and $E_2\subset V_6$.
\end{example}

We now move to defining the cohomology ring of $\fl_n.$ Let $S = \bK[x_1,\ldots,x_n],$ let $e_1,\ldots,e_n$ denote the elementary symmetric polynomials in $x_1,\ldots,x_n$ variables, and let $I_n = (e_1,\ldots,e_n)$ denote the ideal they generate. The quotient $R_n:=S/I_n$ is the \emph{coinvariant algebra} of $\sn.$

\begin{theorem}[\cite{Borel1953}]
\label{thm: Borel - flag variety isomorphism}
    There is an isomorphism of graded rings $H^*(\fl_n,\bK)\cong R_n,$ under which $x_i$ has cohomological degree $2.$ In particular, the odd cohomology of $\fl_n$ vanishes.
\end{theorem}

Concretely, let $0 = \mathcal{S}_0\subset\mathcal{S}_1\subset\cdots\subset\sn = \mathcal{O}^{\oplus n}$ denote the tautological flag of subbundles of $\fl_n,$ and let $L_i = \mathcal{S}_i/\mathcal{S}_{i-1}$ be the associated line bundles. Under the isomorphism of \Cref{thm: Borel - flag variety isomorphism}, the variable $x_i = -c_1(L_i),$ and the relations $e_i = \cdots = e_n = 0$ correspond to the triviality of the total Chern class of $\bigoplus_i L_i\cong\mathcal{O}^{\oplus n}.$

Since $e_1,\dots, e_n$ are algebraically independent, they form a regular sequence. Taking the Koszul resolution of $R_n$ over $S$ and computing the Hilbert series, we obtain
\[
\hilb(R_n;t) = \prod_{i=1}^n(1-t^i)\cdot(1-t)^{-n} = \prod_{i=1}^n(1+t+\cdots+t^{i-1}) = [n]_t!,
\]
where $[k]_t := 1+t+\cdots + t^{k-1}$ and $[k]_t! := [1]_t\cdots [k]_t.$ 
In particular, $\dim_\bK R_n = n!,$ the top nonvanishing graded piece of $R_n$ sits in degree $\binom{n}{2}$ and is one dimensional. $R_n$ also satisfies Poincar{\'e} duality between graded pieces in degrees $k$ and $\binom{n}{2}-k.$ Comparing the above with the standard identity $\sum_{w\in\sn}t^{\ell(w)}=[n]_t!$ verifies that the graded dimensions of $R_n$ are the number of permutations of each length. Furthermore, $H^*(\fl_n;\bK)$ has a basis indexed by $\sn$ and given by the classes of the Schubert varieties $X_w.$ (We make this explicit in the next section.) 

\subsection{Schubert polynomials and the cohomology ideals \texorpdfstring{$J_w$}{Jw} of Schubert varieties}
\label{subsec: Schubert stuff}
The coinvariant algebra $R_n = S/(e_1,\ldots,e_n),$ where $e_i$ is the $i$th elementary symmetric function, has a basis given by the Schubert polynomials $\fraks_w$ (\Cref{def: schubert polynomials}), which are in bijection with the permutations of $\sn.$ 
This basis has a direct correspondence to the geometry of the flag variety \cite{Borel1953}, as we recall below. 

We begin this section by introducing the divided difference operator, which are used to define the Schubert polynomials. We then define the ideal $J_w$ and record some cases in which a convenient generating set of $J_w$ is known. In particular, we highlight work of Gasharov-Reiner \cite{GasharovReiner} which studies the case of $X_w$ being defined by inclusion conditions. In that case, $J_w$ is generated by elementary symmetric polynomials in initial and final blocks of the variables.

\begin{definition}
\label{def:div-diff}
For \(1\leq i<n\), the divided-difference operator is
\[
\partial_i(f(x))=\frac{f(x)-s_i(f(x))}{x_i-x_{i+1}},
\]
where $s_i(f(x))$ is the result of swapping variables $x_i$ and $x_{i+1}$ in $f(x)$. 
If $v = s_{i_1}\cdots s_{i_p}$ is a reduced expression for a permutation $v$, define $\partial_v=\partial_{i_1}\cdots\partial_{i_p}$.
This definition is independent of the choice of reduced word \cite[Chapter II, Equation (2.5)]{Macdonald}.  
\end{definition}

Since $I_n$ is generated by symmetric polynomials and each $\partial_i$ is linear over the ring of symmetric polynomials in $x_i,\ldots,x_n$, then $\partial_i(I_n)\subseteq I_n.$ 
Each $\partial_i$ therefore descends to a well-defined operator on $R_n,$ and hence so does every $\partial_v.$
The divided difference operators act on squarefree monomials in a predictable and simple way, which we record here for later calculations.

\begin{remark}
\label{rem:sqfr-mon-div-diff}    
On a squarefree monomial,
\begin{equation}
\label{eq:squarefree-dd}
 \partial_i(x_F) = 
 \begin{cases}
 x_{F\setminus\{i\}},&i\in F,\ i+1\notin F,\\
 -x_{F\setminus\{i+1\}},&i\notin F,\ i+1\in F,\\
 0,&\text{otherwise}.
 \end{cases}
\end{equation}
\end{remark}

The Schubert polynomials are defined by applying divided difference operators to the \emph{top degree staircase monomial}.

\begin{definition}[\cite{lascoux-schutzenberger}]
\label{def: schubert polynomials}
    Let $w_0 = n(n-1)\cdots 1$ be the longest element of $\sn$ and define
    \[
    \fraks_{w_0} := x_1^{n-1}x_2^{n-2}\cdots x_{n-1}.
    \]
    The \emph{Schubert polynomial} of a permutation $w\in\sn$ is
    \[
    \quad \fraks_w:=\partial_{w^{-1}w_0}(\fraks_{w_0}).
    \]
    Since $\ell(w^{-1}w_0) = \binom{n}{2} - \ell(w)$ and $\partial_v$ lowers degree by $\ell(v),$ the polynomial $\fraks_w$ is homogeneous of degree $\ell(w).$
\end{definition}

The following records the two properties of this construction that we use, namely that it produces a basis and that it is compatible with grading.

\begin{proposition}[\cite{Macdonald,lascoux-schutzenberger}]
\label{prop: schuberts form a basis for Rn}
    The set $\{\fraks_w:w\in\fraks_n\}$ is a $\bK$-basis for $R_n.$ In particular, $\dim_\bK(R_n)_k$ is the number of permutations in $\fraks_n$ of length $k.$
\end{proposition}

We will also need the following relative divided-difference operators in the proof of our main theorem.
\begin{definition}\cite{Macdonald}
\label{def:relative-diff-oper} 
The relative/skew divided-difference operators \(\partial_{u/v}:S\to S\) are defined as the unique linear operators which obey the generalized Leibniz rule 
\begin{equation}
\label{eq:generalized-leibniz}
 \partial_u(fg) = \sum_{v\leq u}
 v\!\left(\partial_{u/v}f\right)\partial_vg.
\end{equation}
\end{definition}

\begin{remark}
    Note that the operators $\partial_u$ and $\partial_{u/v}$ lower the degree by $\ell(u)$ and $\ell(u)-\ell(v),$ respectively \cite[Ch.~II,~Eq.~2.17]{Macdonald}.
\end{remark}

\begin{definition}
\label{def:linear-functionals}
As the Schubert polynomials form a basis for $R_n$ as a $\mathbb{K}$-vector space, there are functionals $\lambda_v:R_n \to \mathbb{K}$ such that  $\lambda_v(\mathfrak{S}_u) = \delta_{u,v}$, where $\delta_{u,v}$ is the Kronecker delta function.
\end{definition}

\begin{remark}
\label{rem:lin-func-length}
If \(f\in(R_n)_{\ell(v)}\) where $\ell(v)$ is the minimum length of $v$ as a reduced word as a product of simple transpositions, then by \cite[Ch.~IV,~Eq.~4.2]{Macdonald} we have
\begin{equation}
\label{eq:coefficient-dd}
 \lambda_v(f)=\partial_v(f).
\end{equation}
\end{remark}

A \textbf{Schubert variety} $X_w$ is a subvariety of $\fl_n$ cut out by certain rank conditions of the form $\dim(V_a\cap E_b)\geq c$ on the moving flag $V_\bullet$ \cite[Ch.~1]{Handbook}. 
Each Schubert polynomial $\fraks_w\in R_n$ is a representative of the cohomology class $[X_w]$ of the Schubert variety corresponding to $w$.

The Schubert varieties $X_w$ are the cell closures of a cell decomposition of $\fl_n$ by Schubert cells. Since each Schubert variety is the union of the Schubert cells that are contained in it, the map on cohomology $\iota^*:R_n\cong H^*(\fl_n;\bK)\to H^*(X_w;\bK)$ induced by the inclusion $\iota:X_w\subseteq \fl_n$ is surjective. Define $J_w\coloneqq \ker(\iota^*)$. By general facts about cohomology, see for example~\cite{ReinerWooYong}, 
\begin{equation}\label{eq: Jw}
J_w = \mathrm{span}_\bK\{\fraks_v:v\in\sn,v\not\leq w\}\subseteq R_n.
\end{equation}

We can also extend the concept to intersections of ideals.
\begin{lemma}
\label{lem: surjective map on homology from flag var to union}
    Let $\calW\subseteq\sn$ be a Bruhat antichain, let $X_\calW = \bigcup_{w\in\calW}X_w,$ and let $J_\calW$ be the kernel of the restriction map $\iota^*:R_n\cong H^*(\fl_n;\bK) \to H^*(X_\calW;\bK).$ Then $\iota^*$ is surjective and 
    \[
    J_\calW = \bigcap_{w\in\calW} J_w = \mathsf{span}_{\bK}\{\fraks_v:v\not\leq w \text{ for every }w\in\calW\}+I_n.
    \]
\end{lemma}
\begin{proof}
    The union of $X_\calW$ is the union of Schubert cells $C_v,$ for $v\leq w$ and $w\in \calW$, and every cell has even real dimension. Therefore, $H_*(X_\calW;\bK)$ has a basis $\{[X_v] : v\leq w,\, w\in \calW\}$ and $H^*(X_\calW;\bK)$ is its dual under the cap product pairing. Since $X_\calW$ is a closed union of cells, the map
    \[
    H^*(\fl_n;\bK)\to H^*(X_\calW;\bK)
    \]
    sends each class $[X_v]$ in $H^*(\fl_n;\bK)$ to either $[X_v]$ in $H^*(X_\calW;\bK)$ if $X_v\subseteq X_{\calW}$ or to $0$ otherwise. Thus, the map is surjective with kernel spanned by the $[X_v]=\fraks_v$ for $v\not\leq w$ for all $w\in \calW$.
\end{proof}

Recall that the essential set of $w$ identifies a minimal collection of incidence conditions defining $X_w.$ We are concerned with certain permutations for which these conditions take on a rather simple form.

\begin{definition}[\cite{GasharovReiner}]
    A Schubert variety $X_w$ is \emph{defined by inclusions} if it is cut out by conditions of the form $E_a\subseteq V_b$ and $V_b\subseteq E_a$ alone.
\end{definition}

Remarkably, this class of Schubert varieties is detected by pattern avoidance. 

\begin{theorem}[\cite{GasharovReiner}]
    The Schubert variety $X_w$ is defined by inclusions if and only if $w$ avoids the patterns $4231,35142,42513,$ and $351624.$
\end{theorem}

\begin{proposition}[\cite{ReinerWooYong}]\label{prop:e-relations}
    Suppose $X_w$ is defined by inclusions. Then $J_w$ is generated by
    \[
    \{e_m(x_1,\ldots,x_r): m = r-s+1,\ldots,r\}\quad\text{for each condition }E_s\subseteq F_r,
    \]
    \[
    \{e_m(x_{r+1},\ldots,x_n): m = s-r+1,\ldots,n-r\}\quad\text{for each condition }F_r\subseteq E_s.
    \]
\end{proposition}

Finally, we isolate the part of an ideal generated by its squarefree monomials, which is of particular importance for us. We note that the following definitions are non-standard.
\begin{definition}
\label{def:sqfr-mon-subideal}
For a homogeneous ideal \(I\subseteq R_n\), define the following ideal in $R_n$,
\[ \mathrm{SqfMon}(I)=\left( x_F:F\subseteq[n],\ x_F\in I\right). \]
\end{definition}

\begin{remark}
\label{rem:sqfr-mon-subideal}
In other words, \(\mathrm{SqfMon}(I)\) is the ideal generated by all squarefree monomials  in \(I\), so
\[ I\text{ is generated by squarefree monomials}\quad\Longleftrightarrow\quad \mathrm{SqfMon}(I)=I.
\]
\end{remark}

\begin{definition} \label{def:forbidden-ascents} 
Fix \(w\in S_n\) with \(w\neq\mathrm{id}\), where \(w(k)=k\) for \(k\notin[i,j]\), $w(i)\neq i$, and $w(j)\neq j$. 
Set 
\[ A(w):=\mathrm{Des}_L(w)^c\cap[i+1,j-2]\]
where $\mathrm{Des}_L(w)$ is the set of left descents of $w$, i.e, 
\[\mathrm{Des}_L(w) := \left\{q\in[n-1]\;:\;w^{-1}(q+1)<w^{-1}(q)\right\}.\]
\end{definition}

\begin{remark}
\label{rem:forbidden-ascents}
For $q\in[i+1,j-2]$,  \(q\in A(w)\) exactly when\footnote{as this is exactly when $\ell(s_q(w)) = \ell(w)+1$} \(w< s_qw\), and $w\in S_n(\mathcal{P})$ exactly when $A(w) = \varnothing$. (See \Cref{lem:AvoidEquivDecrease}.)
\end{remark}

\section{Ideals in coinvariant rings}
\label{sec: ideals in coinvariant rings}
This section describes the behaviour of squarefree monomials in the coinvariant ring, which is what the monomial computations of \Cref{sec: proofs of the main theorems} rest on. The essential feature is that the highest elementary symmetric function, $x_1 \cdots x_m = e_m$ vanishes in $R_m$: it is the only squarefree monomial that does, the annihilator of $x_A$ is generated by the complementary monomial $x_{A^c}$, and divisibility still detects membership in an ideal generated by one or two squarefree monomials.

\begin{lemma}
\label{lem:ann-squarefree-proofonly}
Suppose $x_A$ is a squarefree monomial in $\bK[x_1,\dots,x_m]$, where $A\subseteq [m]$. Then 
$\Ann_{R_m}(x_A)=(x_{A^c})$.
\end{lemma}

\begin{proof}
    Since $x_Ax_{A^c} = x_1\cdots x_m\in I_m$, then $(x_{A^c})\subseteq \Ann_{R_m}(x_A)$. 
    Conversely, given $f\in \Ann_{R_m}(x_A)$ it has a unique decomposition as a polynomial as $f = f_1 + f_2x_{A^c}$, where $f_1$ has no terms divisible by $x_{A^c}$. Then $x_Af_1\in I_m$ as well, and thus by hypothesis on $f_1$, we have $x_Af_1\in (e_1,\dots, e_{m-1})$. Indeed, otherwise if $x_Af_1$ has a nonzero term of the form $e_mg$ for some $g$, then $f_1$ has a nonzero term $x_{A^c}g$, contradicting our assumption on $f_1$.
    
    Since $e_1,\dots, e_{m-1},e_m$ is a regular sequence, $e_m$ is a nonzerodivisor on $S/(e_1,\dots, e_{m-1})$. Since $x_A|e_m$ as polynomials, then $x_A$ is also a nonzerodivisor on $S/(e_1,\dots, e_{m-1})$, thus `canceling' $x_A$ from $x_Af_1 = 0\mod (e_1,\dots, e_{m-1})$ gives $f_1\in (e_1,\dots, e_{m-1})$. Therefore, the image of $f_1$ in $R_m$ is zero. Hence, $f = f_2x_{A^c}\in (x_{A^c})$ as elements in $R_m$. 
\end{proof}

\begin{lemma}
\label{lem:quotient-and-annihilators}
\(R_m/(x_{i_1},\dots,x_{i_s})\) is the ring \(R_{m-s}\) up to relabeling the variables.
\end{lemma}

\begin{proof}
This follows by setting the chosen variables equal to zero in each
\(e_r\).  
\end{proof}

\begin{lemma}
\label{lem:contraction}
Let \(x_F,x_G\in \mathbb{C}[x_1,\ldots,x_m]\) be distinct squarefree monomials of degree at least two, with \(|F\cap G|\neq1\).  
For every squarefree \(x_H\in \mathbb{C}[x_1,\ldots,x_m]\),
\begin{enumerate}
    \item $ \overline{x_H}\in(\overline{x_F})\subseteq R_m
 \quad\Longleftrightarrow\quad
 x_F\mid x_H$
 \item $\overline{x_H}\in(\overline{x_F},\overline{x_G})  \quad\Longleftrightarrow\quad  x_F\mid x_H$ or $x_G\mid x_H.$
\end{enumerate}

\end{lemma}

\begin{proof}
For both statements, the reverse implication is trivially true.

To prove (1) by contrapositive, suppose \(x_F\nmid x_H\), i.e., inside $[m]$, $H^c \cap F\neq \varnothing$.

If \(|H^c \cap F|\geq 2\), choose \(a\in H^c \cap F\) and consider the quotient \(\pi_a: R_m \to R_m/(x_a)\).  
Then \(\pi_a(\overline{x_F})\) becomes zero, while \(\pi_a(\overline{x_H})\) remains a proper squarefree monomial because it still omits
the variables in \((H^c \cap F)\setminus\{a\}\). 
Thus \(\overline{x_H}\not\in(\overline{x_F})\subseteq R_m\).

If \(|H^c \cap F|=1\), i.e. $H^c\cap F = \{a\}$, write
\[ x_F=x_ax_I,\qquad x_H=x_Ix_J.\]
Here \(|I|\geq 1\), since \(\deg x_F\geq2\).
Without loss of generality, assume all $x_i$ appear in $H$ or $F$ so the complementary monomial to \(x_I\) is \(x_ax_J\). 
If \(\overline{x_H}\in (\overline{x_F})\), then there exists an $f$ such that \(\overline{x_H} = \overline{f}\overline{x_F}\), then
\[  \overline{x_I}\overline{(x_J-fx_a)}=\overline{0}. \]
\Cref{lem:ann-squarefree-proofonly} gives \(\overline{x_J-fx_a}\in(\overline{x_ax_J})\), and hence \(\overline{x_J}\in(\overline{x_a})\). 
Modding out by \(\overline{x_a}\) would make \(\overline{x_J}\) zero.  
This is impossible by \Cref{lem:ann-squarefree-proofonly} as \(\overline{x_J}\) is nonzero and proper there because it omits the nonempty factor \(\overline{x_I}\) and $x_a\not\vert x_J$ (as $x_F\not\vert x_H$).  
This proves Part (1).

Now suppose that neither \(x_F\) nor \(x_G\) divides \(x_H\).  
Choose an inclusion-minimal subset \(I\subseteq H^c\subset [m]\) meeting both \(F\) and \(G\).

If \(I\subsetneq H^c\), mod out $R_n$ by the variables indexed by \(I\).  
Both \(x_F\) and \(x_G\) vanish, while the image of \(x_H\) remains a proper squarefree monomial because it still omits every variable in \(H^c\setminus I\).  
Hence, \(\overline{x_H}\notin(\overline{x_F},\overline{x_G})\).

If \(I= H^c\), then either $I = \{a\} \subseteq F\cap G$ or $I = \{b,c\}$ where  $b\in F\setminus (F\cap G)$ and $c\in G\setminus (F\cap G)$.  
 
If $I= \{a\}$, $x_{F\cap G}=\gcd(x_F,x_G)$ has degree at least $2$ as $1\leq |F\cap G| \neq 1$.  
As $a\in F\cap G \cap I \subseteq H^c$, we have \(x_{F\cap G}\nmid x_H\).
Then Part 1 gives \(\overline{x_H}\notin (\overline{x_{F\cap G}})\supseteq(\overline{x_F},\overline{x_G})\).

If $I = \{b,c\}$, then mod out by $x_b$; label the quotient map $\pi_b$.
Then \(\pi_b(\overline{x_F})\) vanishes, while  \(\pi_b(\overline{x_G})\) remains squarefree of degree at least two and does not divide \(\pi_b(\overline{x_H})\), since \(x_c\) is missing.  
Part 1 shows that \(\pi_b(\overline{x_H})\not\in(\pi_b(\overline{x_G}))\).
Thus, \(\overline{x_H}\notin(\overline{x_F},\overline{x_G})\).

We finally see that in all cases,  \(\overline{x_H}\notin(\overline{x_F},\overline{x_G})\) so the contrapositive of the forward implication of Part (2) is proved.
\end{proof}

\begin{lemma}
\label{lem:endpoint-intersection-proofonly}
Let $A, B$ be disjoint subsets of $[m]$ such that $A\neq \emptyset$, and let $x_{A}$ and $x_{B}$ be the corresponding squarefree monomials in $\bK[x_1,\dots,x_m]$.
Then 
\begin{equation}\label{eq:ideal_int_xa}
    \bigcap_{a\in A}\left(x_{a},\;x_{A\setminus a}\,x_B\right) \subseteq (x_A,x_B)\qquad\text{inside }R_m.
\end{equation}
\end{lemma}

\begin{proof}
Induct on $|A|$.  
The case $|A|=1$ is immediate. 
For $|A|>1$, let $a_0$ be an arbitrary element of $A$ and write $x_A=x_{A\setminus a_0}\, x_{a_0}$. 

Given $f\in R_m$ that's in the intersection on the left-hand side of \eqref{eq:ideal_int_xa}, then in particular $f\in (x_{a_0},x_{A\setminus a_0}\,x_B)$, so
\[f=x_{a_0}g+x_{A\setminus a_0}\,x_Bh\]
for some $g,h$. 
For $a\in A\setminus a_0$, since $x_{a}$ divides the second term of this expansion, the condition $f\in(x_{a},x_{A\setminus a}\,x_B)$ gives $x_{a_0}g\in(x_{a},x_{A\setminus a}\,x_B)$. Reducing modulo $x_{a}$ gives
\[ 
\overline x_{a_0}\,\overline g\in (\overline{x_{A\cup B\setminus a}}) \subseteq R_{m-1}, 
\]
where $R_{m-1}$ is the coinvariant algebra with variables $x_i$ for $i\in [m]\setminus a$. Thus, there is some $\overline{\ell}\in R_{m-1}$ such that $\overline{x_{a_0}}\,\overline{g} = \overline{\ell}\,\overline{x_{A\cup B\setminus a}}$. Thus, $\overline{0} = \overline{x_{a_0}}(\overline{g}-\overline{\ell}\,\overline{x_{A\cup B\setminus\{a,a_0\}}})$ so $\overline{g}-\overline{\ell}\,\overline{x_{A\cup B\setminus\{a,a_0\}}}\in \Ann_{R_{m-1}}(\overline{x_{a_0}})$ and hence $\overline{g}-\overline{\ell}\,\overline{x_{A\cup B\setminus\{a,a_0\}}}\in (\overline{x_{[m]\setminus\{a,a_0\}}})$ by \Cref{lem:ann-squarefree-proofonly}. Therefore, $\overline{g}\in (\overline{x_{A\cup B\setminus \{a,a_0\}}})$. Lifting this last statement back to $R_m$, we have $g \in (x_a,x_{A\cup B\setminus \{a,a_0\}})$. Since this is for all $a\in A\setminus a_0$, we have $g\in \bigcap_{a\in A\setminus a_0} (x_a,x_{A\cup B\setminus \{a,a_0\}})$. Since $A\setminus a_0\neq \emptyset$, by our inductive hypothesis $g\in (x_{A\setminus a_0},x_B)$, which implies $f\in (x_A,x_B)$. Our induction is complete.
\end{proof}

\begin{lemma}
\label{lem:local-endpoint-colon-proofonly}
For $1\leq p<q\leq m+1$ with $p<m$,
\[(x_{[q,m]}:x_{[1,p]})\subseteq(x_{[p+1,q-1]},x_{[q,m]})\]
in $R_m$, where by convention $x_{\emptyset}=0$.
\end{lemma}

\begin{proof}
If $q=m+1$, 
\Cref{lem:ann-squarefree-proofonly} gives $(0:x_{[1,p]})=(x_{[p+1,m]})$ and we are done.

Now suppose $q\leq m$.
Let $f \in (x_{[q,m]}:x_{[1,p]})$, meaning that $fx_{[1,p]}\in (x_{[q,m]})$.
Quotienting by $x_a$ for some $a\geq q$ gives $\overline{f}\,\overline{x_{[1,p]}}\in (\overline{x_{[q,m]}}) =(\overline{0}) \subseteq R_{m-1}$, where $R_{m-1}$ is the coinvariant ring on variables $x_{a'}$ for $a'\in [m]\setminus a$. Thus, $\overline{f}\in \mathrm{Ann}_{R_{m-1}}(\overline{x_{[1,p]}})$
so by \Cref{lem:ann-squarefree-proofonly},
$\overline{f}\in(\overline{x_{[p+1,m]\setminus a}})$. Therefore, we have $f\in(x_a,x_{[p+1,m]\setminus a})$
for every $a\geq q$.  Applying \Cref{lem:endpoint-intersection-proofonly} with $A = [q,m]$ and $B = [p+1,q-1]$ then gives
$f\in(x_{[p+1,q-1]},x_{[q,m]})$.
\end{proof}

\section{Doubly pruned simplicial complexes}
\label{sec: intro to doubly pruned simplicial complexes}
In this section we introduce a family of simplicial complexes that we call \emph{doubly pruned}\footnote{We call these complexes doubly pruned because each complex is a simplex on $[i, j]$ with two pieces cut away: the faces containing an initial interval $[i, p],$ and the faces containing a final interval $[q, j].$}. These are simplicial complexes $\Delta_w$ that are determined by two interval nonfaces and singletons outside of it (\Cref{def: doubly pruned simplicial complexes}). For now, this section will focus on the properties of these complexes that are required in the proofs of the main theorem. In particular, we establish a dictionary between $\sn(\calP)$ and these complexes that match Bruhat order with containment and the closure properties needed to pass from a single permutation to an antichain. 
The combinatorial and homological consequences of these complexes are deferred to \Cref{sec: the combinatorics of doubly pruned simplicial complexes}.

\begin{definition}
\label{def: doubly pruned simplicial complexes}
A simplicial complex $\Delta$ on vertex set $[n]$ is \emph{doubly pruned} if $\Delta = \emptyset$ or for some $1\leq i\leq j\leq n$ its nonfaces are generated by the singletons $\{k\}$ for $k\in[n]\setminus[i,j]$ together with at most one interval $[i,p]$ and at most one interval $[q,j],$ where $i\leq p\leq j$ and $i\leq q\leq j.$
\end{definition}
Recall from the introduction that $S_n(\calP)$ is the set of permutations avoiding the patterns found in $\calP,$ where
\begin{align*}
\calP &\coloneqq \{ 2143,\ 2341,\ 2413,\ 2431,\ 4123,\ 4213,\ 4231 \}.
\end{align*}
By \Cref{lem: doubly pruned complexes are Schubert complexes}, a simplicial complex $\Delta_w$ for $w\in\sn(\calP)$ is doubly pruned, since structurally, its minimal nonfaces are exactly the two intervals of \Cref{def:Deltaw}, together with the singletons outside the interval $[i(w),j(w)].$ By working with this more combinatorial definition, rather than with $\Delta_w$ directly, allows us to discard the singleton nonfaces (a trivial operation) and prove the relevant properties. In particular, we show in \Cref{lem:double pruned intersections} show that intersections of doubly pruned complexes are also doubly pruned. This includes the intersections $\Delta_v\cap\Delta_w$ arising from antichains in $\sn(\calP).$ 
In \Cref{lem: doubly pruned complexes are Schubert complexes} we show that $w\mapsto(i(w),j(w),p(w),q(w))$ is an injection from $\sn(\calP)$ to the doubly pruned complexes on $[n].$ 
Moreover, \Cref{prop:Inclusion} identifies Bruhat order on $\sn(\calP)$ with containment of the corresponding complexes via the rank matrix $R_w.$ 
This gives a dictionary between the poset $(\sn(\calP),\leq)$ and a subset of the doubly pruned complexes ordered by inclusion. 
This allows us to show, in \Cref{thm: doubly pruned complexes imply vertex decomposable}, that every doubly pruned complex, and in particular every that $\Delta_w,$ is vertex decomposable. Thus, they are also shellable and sequentially Cohen--Macaulay. \Cref{prop: facet description of doubly pruned complexes} and \Cref{cor: purity of doubly pruned complexes} describe the facets explicitly and show that $\Delta_w$ fails to be pure exactly when
$i(w)<q(w)\leq p(w)<j(w)$, in which case its facets have two different dimensions (\Cref{ex: 3412} exhibits this directly). 

In \Cref{subsec:hilb} we use the same description, the ideal $I_{\Delta_w}$ generated by two interval monomials together with a run of variables) in order to compute
the Hilbert series of $R_{\Delta_w}=S/I_{\Delta_w}$
(\Cref{thm: hilbert series - single Schubert variety}). Here we split the proof into the case where
the two interval nonfaces are disjoint and the case where they overlap.
\subsection{Pattern avoidance and the complexes \texorpdfstring{$\Delta_w$}{Delta w}}
We begin by analyzing the structure of the permutations in $\sn(\calP).$
\begin{lemma}\label{lem:AvoidEquivDecrease}
    The set $S_n(\calP)$ is equal to the set of all permutations $w$ such that if $i$ is the smallest non-fixed point of $w$ and $j$ is the largest non-fixed point of $w$, then $i+1,i+2,\dots, j-1$ appear in decreasing order in $w$. 
\end{lemma}
\begin{proof}
    We prove the first containment via the contrapositive. Suppose $i+1,i+2\dots, j-1$ do not appear in decreasing order. Without loss of generality, we may assume $i=1$ and $j=n$ since none of the pattern start with $1$ or end with $4$ so initial and final fixed points will not be part of any pattern in $\calP$. Thus, our assumptions are that $2,3,\dots, n-1$ do not appear in decreasing order, $w(1)\neq 1$, and $w(n)\neq n$. Then there is some $1<k<\ell<n$ such that $w^{-1}(k)<w^{-1}(\ell)$. If $w^{-1}(n)<w^{-1}(1)$, then the letters $1,k,\ell,n$ will form a $2341,2413,2431,4123,4213$, or $4231$ pattern depending on the six positions of $k,\ell$ relative to $n,1$. 
    
    Otherwise, $w^{-1}(1)<w^{-1}(n)$. If there exists a letter $a$ to the left of $1$ and letter $b$ to the right of $n$ such that $a<b$, then $[a,1,n,b]$ forms a $2143$ pattern. Otherwise, we have that every letter left of $1$ is greater than every letter right of $n$ (and there must exist at least one of each). Since $2,3,\dots, n-1$ are not in decreasing order, then we again have some $1<k<\ell<n$ such that $k$ is left of $\ell$. If we assume further that $k$ is left of $1$ and $\ell$ is between $1$ and $n$, then $[k,\ell,n,b]$ forms a $2341$ pattern, where $b$ is any letter right of $n$, which must have $b<k$ by our assumption above. If instead $k$ is between $1$ and $n$ and $\ell$ is right of $n$, then $[a,1,k,\ell]$ forms a $4123$ pattern, where $a$ is any letter left of $1$, which has the property that $a>\ell$ by our assumption above. Finally, the last case is when $k$ and $\ell$ are both between $1$ and $n$. Then $[a,k,\ell,b]$ forms a $4231$ pattern, where $a$ is any letter left of $1$ and $b$ is any letter right of $n$. Thus, in every case $w$ has a pattern in the set $\calP$.
    To prove the other containment, assume that $i+1,i+2,\dots, j-1$ appear in $w$ in decreasing order. By way of contradiction, assume that $w$ contains one of the patterns in $\calP$. We may again assume without loss of generality that $i=1$ and $j=n$. If $w$ contains one of the patterns in $\calP$, then the $2$ and $3$ in the pattern cannot both be from $2,3,\dots, n-1$ since these letters decrease in $w$ but increase in the pattern. This results in a contradiction, since it forces the $2$ or $3$ to represent either $1$ or $n$ in $w$. Therefore, $w$ must avoid $\calP$.
\end{proof}

This leads to our definition of the complexes that are the main focus of this paper, which we will call $\Delta_w,$ with $w\in\sn(\calP).$ Immediately afterwards, in \Cref{lem: doubly pruned complexes are Schubert complexes}, we show that any $\Delta_w$ is a doubly pruned simplicial complex. Thereby, allowing us to prove results for doubly pruned simplicial complexes and then specialize to these $\Delta_w.$

\begin{definition}
\label{def:Deltaw}
    Given $w\in S_n(\calP)$, define
    \[
    i(w) = \min\{k:w(k)\neq k\}\quad\text{ and }\quad j(w) =\max\{k:w(k)\neq k\},
    \]
    the smallest and largest non-fixed points of $w$. (Note that if $w$ is the identity set $i = n+1, j = 0.$) Set $p(w) = w^{-1}(i(w))$ and $q(w) = w^{-1}(j(w))$ so that $i(w)\leq p(w)\leq j(w)$ and $i(w)\leq q(w)\leq j(w)$.
    
    Define $\Delta_w$ to be the simplicial complex on ground set $[n]$ whose minimal nonfaces are
\begin{itemize}
    \item singletons: $\{1\},\ldots\{i(w)-1\},\{j(w)+1\},\ldots,\{n\};$
    \item intervals: $[i(w),p(w)]$ and $[q(w),j(w)].$
\end{itemize}
\end{definition}
In the next lemma, for an interval $A\subseteq[n],$ we use the notation $\mathscr{P}(A)$ to mean the \emph{power set of} $A.$
\begin{lemma}
\label{lem: doubly pruned complexes are Schubert complexes}
    A  simplicial complex $\Delta$ on $[n]$ is doubly pruned if and only if one of the following (mutually exclusive) conditions holds:
    \begin{enumerate}[(1)]
        \item $\Delta$ is the simplex $\Delta = \mathscr{P}([i,j])$ for some $1\leq i\leq j\leq n$,
        \item $\Delta$ has minimal nonfaces $\{1\}$,$\{2\}$,\dots, $\{i-1\}$, $\{j+1\},\dots, \{n\}$, $[i,p]$, and $[p,j]$ for some $1\leq i<p<j\leq n$, or
        \item $\Delta = \Delta_w$ for some $w\in S_n(\calP)$.
    \end{enumerate}
\end{lemma}
\begin{proof}
    If $\Delta = \varnothing$, then $\Delta = \Delta_{id}$ by convention, so we assume $\Delta\neq \varnothing$.
    
    In all three cases, they clearly satisfy the definition of doubly pruned complex; in Case 3, given $\Delta_w$ for $w\in S_n(\calP)$, take $i = i(w), j = j(w), p = p(w),$ and $q = q(w)$ in the definition of doubly-pruned complexes.
    Conversely, given any doubly-pruned complex $\Delta$ on $[n]$, if it has no minimal nonfaces that are not vertices, it clearly has the form $\Delta=\mathscr{P}([i,j])$.
    If it has two minimal nonfaces of the form $[i,p]$ and $[p,j]$, then we are in Case 2.
    If it has a single minimal nonface that is not a vertex, either $[i,p]$ (in which case define $q=i$) or it is $[q,j]$ (in which case define $p=j$) , or two minimal nonfaces with $p\neq q$, then define $w$ to be the unique permutation such that $i(w)=i$, $j(w)=j$, $p(w)=p$, and $q(w)=q$.
    These cases are mutually exclusive: in Case (1) there is no minimal nonface of cardinality at least two, in Case (2) there are two and they meet in the single vertex $p$, and in Case (3) this cannot occur, since $p(w) = q(w)$ would force $w^{-1}(i(w)) = w^{-1}(j(w))$.
\end{proof}

\begin{example}
\label{ex: 142536}
    Let $w = 142536 \in S_6(\calP)$, so that $w^{-1} = 135246$. The non-fixed points
    of $w$ are $2,3,4,5$, so $i(w) = 2$ and $j(w) = 5$. The values strictly between
    $i(w)$ and $j(w)$ are $3$ and $4$, and they appear in $w$ in the order $4,3$,
    which is decreasing, so $w$ is $\calP$-avoiding by
    \Cref{lem:AvoidEquivDecrease}. Reading off the remaining data of
    \Cref{def:Deltaw},
    \[
    p(w) = w^{-1}(i(w)) = w^{-1}(2) = 3,
    \qquad
    q(w) = w^{-1}(j(w)) = w^{-1}(5) = 4,
    \]
    so the minimal nonfaces of $\Delta_w$ are the singletons $\{1\}$ and $\{6\}$,
    coming from $[6]\setminus[2,5]$, together with the two intervals
    \[
    [i(w),p(w)] = \{2,3\}
    \qquad\text{and}\qquad
    [q(w),j(w)] = \{4,5\}.
    \]
    In particular $\Delta_w$ is doubly pruned, with $i = 2$, $j = 5$, $p = 3$ and
    $q = 4$ in \Cref{def: doubly pruned simplicial complexes}. A subset
    $F\subseteq[6]$ is a face of $\Delta_w$ exactly when $F\subseteq\{2,3,4,5\}$ and
    $F$ contains at most one of $2,3$ and at most one of $4,5$, so the facets are
    \[
    \{2,4\},\qquad \{2,5\},\qquad \{3,4\},\qquad \{3,5\},
    \]
    all of dimension $1$. Pictorially, $\Delta_w$ is a $4$-cycle,
    \begin{center}
    \begin{tikzpicture}[scale=0.4, vtx/.style={circle,fill=black,inner sep=1.7pt}]
      \coordinate (v2) at (-2,0);
      \coordinate (v4) at (0,2);
      \coordinate (v3) at (2,0);
      \coordinate (v5) at (0,-2);
      \draw[thick] (v2) -- (v4) -- (v3) -- (v5) -- cycle;
      \node[vtx] at (v2) {}; \node[vtx] at (v4) {};
      \node[vtx] at (v3) {}; \node[vtx] at (v5) {};
      \node[left]  at (v2) {\small $2$};
      \node[above] at (v4) {\small $4$};
      \node[right] at (v3) {\small $3$};
      \node[below] at (v5) {\small $5$};
      \node[anchor=east] at (-3,0) {\large $\Delta_w \;=$};
    \end{tikzpicture}
    \end{center}
    the join of the two-point complex on $\{2,3\}$ with the two-point complex on
    $\{4,5\}$.
\end{example}
As one would hope, the Bruhat order carries over into containment of simplicial complexes.

\begin{proposition}
\label{prop:Inclusion}
    Given $v,w\in S_n(\calP)$, the following are equivalent:
    \begin{itemize}
        \item $v\leq w$,
        \item $i(w) \leq i(v)$, $j(v)\leq j(w)$, if $i(v)=i(w)$ then $p(v)\leq p(w)$, and if $j(v)=j(w)$ then $q(w) \leq q(v)$,
        \item $\Delta_v\subseteq \Delta_w$.
    \end{itemize}
\end{proposition}
\begin{proof}
    We first show the top two bullets points are equivalent. Suppose $v,w\in S_n(\calP)$ such that $v\leq w$. Then $r_v(p,q) \geq r_w(p,q)$ for all $p,q$. Since for all $\ell<i(w)$, $r_v(\ell,\ell)\geq r_w(\ell,\ell) = \ell$, then $v$ must have fixed points up to index $i(w)-1$, meaning $i(v)\geq i(w)$. A similar argument shows $j(w)+1,\dots,n$ must be fixed points of $v$ so $j(v)\leq j(w)$.
    Suppose $i(v)=i(w)$, and set $i$ to be this common value. Then $i=r_v(p(v),i) \geq r_w(p(v),i)$ implies that the position of $i$ in $w$ is at or to the right of $p(v)$, meaning $p(w) \geq p(v)$. Supposing $j\coloneqq j(v)=j(w)$, since $r_v(q(w)-1,j-1)\geq r_w(q(w)-1,j-1) = q(w)-1$, this forces $v(\ell)\leq j-1$ for all $\ell<q(w)-1$, which implies $j$ must be weakly to the right of position $q(w)-1$, so $q(v)\geq q(w)$.
    Conversely, suppose $v,w\in S_n(\calP)$ and the three conditions hold. Then we have $r_v(a,b)= r_w(a,b)$ automatically for all $(a,b)$ except for those in $[i(w),j(w)-1]\times [i(w),j(w)-1]$. Thus, we may assume without loss of generality that $i(w)=1$ and $j(w)=n$, and we only need to check the inequality for $(a,b)$ in the square $[i(v),j(v)-1]\times [i(v),j(v)-1]$.
    It can be checked that for any $a,b<n$, 
    \[
    r_w(a,b) = \delta_{p(w)\leq a} + \max(0,a+b-n+1-\delta_{p(w)\leq a} - \delta_{q(w)\leq a})
    \]
    and shifting this formula by $i(v)-1$, we get
    \begin{align*}
    r_v(a,b) 
    & = i(v)-1 + \delta_{p(v)\leq a} + \max(0,(a-i(v)+1) \\ 
    & + (b-i(v)+1) - (j(v)-i(v)+1)+1 - \delta_{p(v)\leq a} - \delta_{q(v)\leq a}).\\
    &= i(v)-1 + \delta_{p(v)\leq a} + \max(0,a+b-j(v)-i(v)+2 - \delta_{p(v)\leq a} - \delta_{q(v)\leq a}).
    \end{align*}
    When both nonzero maxima are achieved, we have $r_w(a,b) = a+b-n+1-\delta_{q(w)\leq a}$ and $r_v(a,b) = a+b-j(v)+1-\delta_{q(v)\leq a}$. Then $r_w(a,b)\leq r_v(a,b)$ follows easily from the fact that if $j(w)=n=j(v)$ then $q(w)\leq q(v)$. If the max in $r_w(a,b)$ is $0$, then we have $r_w(a,b) = \delta_{p(w)\leq a}$, which is less than or equal to $i(v)-1+\delta_{p(v)\leq a}$ since $i(w)=1=i(v)$ implies $p(v)\leq p(w)$. Finally, if the max in $r_v(a,b)$ is $0$, then $a+b+2-j(v)-\delta_{q(v)\leq a}<i(v) + \delta_{p(v)\leq a}$. Using this inequality, we see $r_v(a,b) = i(v)-1+\delta_{p(v)\leq a} \geq a+b-n+1-\delta_{q(w)\leq a} \geq r_w(a,b)$. Thus, $r_w(a,b)\leq r_v(a,b)$ for all $a,b$, and thus $v\leq w$.
    
    The fact that the second and third bullet points are equivalent follows easily from the definition of $\Delta_v$ in terms of its nonfaces.
\end{proof}

\subsection{Intersections of doubly pruned complexes}
Doubly pruned complexes are closed under intersection. We prove this, and then
describe explicitly the intersections that arise from an antichain in $S_n(\calP)$.
\begin{proposition}
\label{prop:intersection is doubly pruned}
    The intersection of two doubly pruned simplicial complexes on $[n]$ is doubly pruned.
\end{proposition}
\begin{proof}
    Let $\Delta_1$ and $\Delta_2$ be doubly pruned, and for $r\in\{1,2\}$ write the nonfaces of $\Delta_r$ as generated by the singletons $\{k\}$ with $k\in[n]\setminus[i_r,j_r]$ together with at most one interval $[i_r,p_r]$ and at most one interval $[q_r,j_r]$. 
    A subset of $[n]$ is a nonface of $\Delta_1\cap\Delta_2$ if and only if it is a nonface of $\Delta_1$ or a nonface of $\Delta_2.$ Thus, the nonfaces of the intersection are generated by the union of the two generating sets. 
    Set $I = \max\{i_1,i_2\}$ and $J = \min\{j_1,j_2\}$, so that $[i_1,j_1]\cap[i_2,j_2] = [I,J]$.

    If $I>J$, then every $k\in[n]$ lies outside $[i_1,j_1]$ or outside $[i_2,j_2]$, so every singleton is a nonface and $\Delta_1\cap\Delta_2 = \emptyset$, which is doubly pruned by \Cref{def: doubly pruned simplicial complexes}. Assume then that $I\leq J$. We know that the singleton generators of the intersection are the $\{k\}$ with $k\not\in[i_1,j_1]$ together with the $\{k\}$ with $k\notin[i_2,j_2],$ which means they are precicely the $\{k\}$ with $k\notin[I,J]$.

    Now consider the at most four interval generators.
    If $i_r<I$, then $i_r\not\in[I,J]$, so $\{i_r\}$ is a singleton generator of the intersection contained in $[i_r,p_r]$. Thus, $[i_r,p_r]$ is not a minimal nonface and may be discarded. It is analogous for $[q_r,j_r]$ when $j_r>J$. Each surviving interval generator, therefore, has left endpoint $I$ or right endpoint $J$.

    There are at most two surviving generators with left endpoint $I,$ and any two intervals sharing a left endpoint are nested, so at most one of them is minimal.
    The same argument applies to the at most two surviving generators with right endpoint $J$. Hence the minimal nonfaces of $\Delta_1\cap\Delta_2$ are the singletons $\{k\}$ with $k\notin[I,J]$, together with at most one interval $[I,p]$ and at most one interval $[q,J]$, so $\Delta_1\cap\Delta_2$ is doubly
    pruned.
\end{proof}

\begin{corollary}\label{cor: lattice}
    The collection of doubly pruned simplicial complexes on $[n]$ forms a lattice under containment of simplicial complexes, in which the meet of $\Delta_1$ and $\Delta_2$ is $\Delta_1\cap\Delta_2$ and the join is the smallest doubly pruned complex containing $\Delta_1\cup\Delta_2.$
\end{corollary}

\begin{proof}
    By \Cref{prop:intersection is doubly pruned}, the intersection of two doubly pruned simplicial complexes on $[n]$ is a doubly pruned simplicial complex, hence the poset of doubly pruned simplicial complexes is a meet-semilattice, in which the meet is given by intersection. Since there is a unique top element of the finite poset, namely the simplex supported on $[1,n]$, it is therefore a lattice \cite[Prop.~3.3.1]{StanleyEC1}. 
    The doubly pruned complexes containing both $\Delta_1$ and $\Delta_2$ form a nonempty finite collection, as it contains the simplexes supported on $[1,n],$ so their intersection is the smallest doubly pruned complex containing $\Delta_1\cup\Delta_2.$
\end{proof}

\begin{remark}
    The subposet $\{\Delta_w : w\in S_n(\calP)\}$ is not a lattice. Indeed, observe that meets and joins do not exist when $n=4$; see Figure~\ref{fig: s4 simplicial complex poset}.
\end{remark}

\begin{corollary}
\label{lem:double pruned intersections}
    Let $\{v,w\}\subseteq S_n(\calP)$ be an antichain such that $i(v)\geq i(w)$, and if they are equal $p(v)< p(w)$. Then the nonfaces of $\Delta_v\cap\Delta_w$  are generated by 
    \[1,\dots,i(v)-1,\qquad j(w)+1,\dots,n,\qquad [i(v),p(v)],  [q(w),j(w)].\]
\end{corollary}
\begin{proof}
    As the intersection, the nonfaces of $\Delta_v\cap \Delta_w$ are generated by the minimal nonfaces of $\Delta_v$ union the nonfaces of $\Delta_w$.
    Thus, the intersection's minimal nonfaces are generated by 
    \begin{itemize}
        \item $1,\dots,\max(i(v)-1,i(w)-1),$
        \item $\min(j(v)+1,j(w)+1),\dots,n,$ and
        \item $[i(v),p(v)], [q(v),j(v)],[i(w),p(w)], [q(w),j(w)].$
    \end{itemize}
    By assumption, $\max(i(v)-1,i(w)-1) = i(v)-1$.
    
    If $j(w)>j(v)$, then either ($i(v)>i(w)$ and  $j(w)>j(v)$) or ($i(v)=i(w)$, $p(v)< p(w)$, and  $j(w)>j(v)$).
    In either case, by \Cref{prop:Inclusion}, $\Delta_v\subseteq  \Delta_w$, so $v\leq w$ which violates the assumption that they are an antichain.
    Thus, $\min(j(v)+1,j(w)+1)=j(w)+1$.
    If $i(w)<i(v)$, then $\{i(w)\}$ is a minimal nonface of $\Delta_v$ strictly contained in $[i(w),p(w)]$.
    If $i(w)=i(v)$, then $p(v)<p(w)$ so $[i(v),p(v)]$ is a minimal nonface of $\Delta_v$ strictly contained in $[i(w),p(w)]$.
    In either case, $[i(w),p(w)]$ is not a minimal nonface of the intersection.
    If $j(w)<j(v)$, then $j(v)$ is a minimal nonface of $\Delta_w$ strictly contained in $[q(v),j(v)]$.
    If $j(w)=j(v)$ and $q(w)\leq q(v)$, then $v\leq w$ by \Cref{prop:Inclusion} which violates the assumption that the pair is an antichain.
    If $j(w)=j(v)$ and $q(w)> q(v)$,
    then $[q(w),j(w)]$ is a minimal nonface of $\Delta_w$ strictly contained in $[q(v),j(v)]$.
    In both possible cases, $[q(v),j(v)]$ is not a minimal nonface of the intersection.
    Thus, the (minimal) nonfaces of $\Delta_v\cap \Delta_w$ are generated by
    \[1,\dots,i(v)-1,\qquad j(w)+1,\dots,n,\qquad [i(v),p(v)],  [q(w),j(w)].\qedhere\]
\end{proof}

\begin{corollary}
    \label{lem:antichain ordering}
    Let $\calW = \{w_1,\dots,w_m\} \subseteq S_n(\calP)$ be an antichain ordered lexicographically by the tuples \([i(w_k),j(w_k),-p(w_k),-q(w_k)]\).
    Then for all $k\in [1,m-1]$
    \[\Delta_{w_k}\cap \Delta_{w_m}\subseteq \Delta_{w_{k+1}}\cap \Delta_{w_m}\quad\text{ and }\quad \Delta_{w_k}\cap \Delta_{w_1}\supseteq \Delta_{w_{k+1}}\cap \Delta_{w_1}.\]
\end{corollary}
\begin{proof}
    As the inclusion $\Delta_{w_{m-1}}\cap \Delta_{w_m} \subseteq \Delta_{w_m}\cap \Delta_{w_m}$ is immediate, let $k\in[1,m-2]$.
    In order to conclude the inclusion $\Delta_{w_k}\cap \Delta_{w_m}\subseteq \Delta_{w_{k+1}}\cap \Delta_{w_m}$, it suffices to show that a set which generates the nonfaces of $\Delta_{w_k}\cap \Delta_{w_m}$ generates a set that generates the nonfaces of $\Delta_{w_{k+1}}\cap \Delta_{w_m}$.
    
    By the ordering either $i(w_m)> i(w_k)$ or, if they are equal, then $j(w_m)\geq j(w_k)$.
    If $i(w_m)= i(w_k)$ and $j(w_m)= j(w_k)$, then $-p(w_m)\geq -p(w_k)$; if those were equal, $\calW$ would not be an antichain so $p(w_k)> p(w_m)$.
    If $i(w_m)= i(w_k)$ and $j(w_m)> j(w_k)$, then in order for $\calW$ to be an antichain, $p(w_k)> p(w_m)$.
    Thus, in all cases, we may apply \Cref{lem:double pruned intersections}, and so the nonfaces of $\Delta_{w_k}\cap \Delta_{w_m}$
    \[1,\dots,i(w_m)-1,\qquad j(w_k)+1,\dots,n,\qquad [i(w_m),p(w_m)],  [q(w_k),j(w_k)].\]
    The exact same argument works replacing $k$ with $k+1$, and so the  nonfaces of $\Delta_{w_{k+1}}\cap \Delta_{w_m}$
    \[1,\dots,i(w_m)-1,\qquad j(w_{k+1})+1,\dots,n,\qquad [i(w_m),p(w_m)],  [q(w_{k+1}),j(w_{k+1})].\]
    In order to conclude the inclusion from those lists, it suffices to show that the first list generates the second list.
    In other words, it suffices to show that $j(w_k)\leq j(w_{k+1})$ and if they are equal, then $q(w_k)\geq q(w_{k+1})$ for $k\leq m-2$.
    By the ordering $i(w_k)\leq i(w_{k+1})$.
    If $i(w_k)= i(w_{k+1})$, then by the ordering $j(w_k)\leq j(w_{k+1})$.
    If $i(w_k)< i(w_{k+1})$ and $j(w_k)>j(w_{k+1})$, then $w_{k+1}\leq w_k$ by \Cref{prop:Inclusion}, so if $i(w_k)< i(w_{k+1})$, then $j(w_k)\leq j(w_{k+1})$.
    In either case, we have $j(w_k)\leq j(w_{k+1})$.
    If $j(w_k)< j(w_{k+1})$, then the conclusion holds, so we assume $j(w_k)= j(w_{k+1})$.
    It now suffices to show $q(w_k)\geq q(w_{k+1})$.
    As $i(w_k)\leq i(w_{k+1})$ by the ordering, we deal with it in two cases depending on if that is an equality or not.
    
    We deal with the equality case first.
    If $i(w_k)= i(w_{k+1})$, then as $j(w_k)= j(w_{k+1})$, by the ordering $p(w_k)\geq p(w_{k+1})$.
    If $p(w_k)= p(w_{k+1})$, then they are comparable by \Cref{prop:Inclusion} as the $i$ and $j$ coordinates were already equal, which is impossible so $p(w_k)> p(w_{k+1})$.
    If $p(w_k)> p(w_{k+1})$ and $q(w_k)\leq q(w_{k+1})$, then they are again comparable by \Cref{prop:Inclusion}, so $q(w_k)>  q(w_{k+1})$.
    Now we deal with the case of a strict inequality, i.e., $i(w_k)< i(w_{k+1})$.
    Recall we still have $j(w_k)= j(w_{k+1})$.
    If $q(w_k)\leq q(w_{k+1})$, then they are again comparable by \Cref{prop:Inclusion}, so $q(w_k)>  q(w_{k+1})$.
    
    An entirely symmetric argument proves the other inclusion.
\end{proof}
\section{Proofs of Theorems A and B}
\label{sec: proofs of the main theorems}
In this section, we prove our main theorems, \Cref{thm:main-A} and \Cref{thm:main-B}. We do this by first proving the most immediate directions of the equivalences in \Cref{thm:main-A}. We then prove \Cref{thm:main-B} in several lemmas and then deduce the remaining directions of \Cref{thm:main-A} as a corollary.

\subsection{Proofs of some parts of \texorpdfstring{\Cref{thm:main-A}}{Theorem A}}

In this subsection, we prove that implications in \Cref{thm:main-A} that follow easily from work of Gasharov--Reiner~\cite{GasharovReiner}, namely \Cref{cond:A2} $\implies$ \Cref{cond:A3} $\implies$ \Cref{cond:A1}. We complete the final implication \Cref{cond:A1} $\implies$ \Cref{cond:A2} by directly proving the more general statement that \Cref{cond:B1} $\implies$ $\calW\subseteq S_n(\calP)$ in \Cref{prop:1-implies-avoidance}.

\begin{lemma}[\Cref{cond:A3} implies \Cref{cond:A1}]
\label{lem: A3 implies A1}
    If $X_w$ is defined by the inclusions stated in \Cref{cond:A3}, then either $w=\id$ and $J_\id=(x_1,\dots, x_n)$ or $w\neq \id$ and $J_w$ is generated by the following squarefree monomials, where $i = i(w)$ and $j=j(w)$,
    \begin{itemize}
        \item The variables $x_1,\dots,x_{i-2}$, and $x_{i-1}$,
        \item The products $x_i\cdots x_{w^{-1}(i)}$ and $x_{w^{-1}(j)}\cdots x_j$, and
        \item The variables $x_{j+1},\dots,x_{n-1}$, and $x_n$.
    \end{itemize}
\end{lemma}

\begin{proof}
    If $w=\id$, then as ideals in $R_n$,
    \[
    J_w = (\fraks_{s_i}:1\leq i<n) = (x_1+x_2+\cdots + x_i:1\leq i<n)= (x_1,\dots, x_n).
    \]
    Otherwise, suppose $w\neq \id$ and $X_w$ is defined by the inclusion stated in \Cref{cond:A3}. Applying \Cref{prop:e-relations} to our situation, $J_w$ is generated by the following elementary symmetric polynomials\footnote{Each of these relations can alternatively be obtained from a simple Chern class analysis on the tautological quotient line bundles on $\fl_n$.}:
    \begin{itemize}
        \item $e_m(x_1,\dots, x_r)$ for $r<i$, $m = 1,\dots, r$,
        \item $e_m(x_1,\dots, x_{w^{-1}(i)})$ for $m=i-w^{-1}(i)+1,\dots, w^{-1}(i)$,
        \item $e_m(x_{w^{-1}(j)},\dots, x_n)$ for $m = j-w^{-1}(j)+1,\dots, n-w^{-1}(j)-1$,
        \item $e_m(x_{r+1},\dots, x_n)$ for $r>j$, $m=1,\dots, n-r$.
    \end{itemize}
    The first and last set of generators can be simplified to generators $x_1,\dots, x_{i-1}$ and $x_{j+1},\dots, x_n$. Modding out by these variables, we see that the second family of generators then simplifies to $e_{i-w^{-1}(i)+1}(x_i,\dots, x_{w^{-1}(i)}) = x_i\cdots x_{w^{-1}(i)}$. Furthermore, the third set of generators reduces to $e_{j-w^{-1}(j)+1}(x_{w^{-1}(j)},\dots, x_j) = x_{w^{-1}(j)}\cdots x_j$. 
    Thus, $J_w$ is generated by the squarefree monomials stated above.
\end{proof}

\begin{lemma}[\Cref{cond:A2} implies \Cref{cond:A3}]\label{lem:2implies3}
\label{lem: A2 implies A3}
    Suppose $w\in S_n(\calP)$, then $X_w$ is defined by the inclusions in \Cref{cond:A3}.
\end{lemma}

\begin{proof}
    If $w\in S_n(\calP)$, then $w$ also avoids \(\{4231, 35142, 42513,  351624\}\). Thus, $X_w$ is defined by inclusions in the sense of Gasharov--Reiner \cite{GasharovReiner}. In fact, the inclusion conditions for $w\in S_n(\calP)$ take the following precise form.

    If $w=\id$, then $X_\id$ is a single point consisting of the unique flag defined by $V_b=E_a$ for all $i$, where $E_\bullet$ is the standard flag. Otherwise, if $w\neq \id$ then by \Cref{lem:AvoidEquivDecrease}, the letters $i+1,i+2,\dots, j-1$ are in decreasing order in $w$, where $i$ and $j$ are the smallest and largest non-fixed-points of $w$, respectively. It is not hard to check that if $w^{-1}(i)\neq j$ and $w^{-1}(j)\neq i$, then the essential set $\mathrm{Ess}(w)$ of $w$ is
    \begin{gather*}
    \{(2,1),(3,2),\dots, (i,i-1)\}\,\cup\, \{(i+1,w^{-1}(i))\}\\ 
    \cup\,  \{(j,w^{-1}(j)-1)\}\, \cup \,\{(j+1,j),(j+2,j+1),\dots, (n-1,n)\}.
    \end{gather*}
    See \Cref{ex:Essential} for the case of $w=1657423$ where $i=2$ and $j=7$. 
    
    Each of these sets immediately translates into the inclusion conditions in \Cref{cond:A3}. For example, for $k<i$, $(k+1,k)\in \mathrm{Ess}(w)$ and $r_w(k,k) = k$. Thus, a defining condition of $X_w$ is $\dim_{\bC}(V_k\cap E_k)\geq k$, forcing $V_k = E_k$. For $(i+1,w^{-1}(i))$, $r_w(i,w^{-1}(i)) = i$, so $\dim_\bC(V_{w^{-1}(i)}\cap E_i)\geq i$, forcing the containment $E_i\subseteq V_{w^{-1}(i)}$. The other pairs in the essential set similarly translate to the other two sets of containment conditions in \Cref{cond:A3}. 
    
    In the case $w^{-1}(j)=i$, the pair $(j,w^{-1}(j)-1)$ is not in the essential set. This is simply because the corresponding condition $\dim_{\bC}(V_{i-1}\cap E_{j-1})\geq i-1$ translates to $V_{i-1}\subseteq E_{j-1}$, which is redundant from the condition $V_{i-1} = E_{i-1}$. Similarly, if $w^{-1}(i) = j$ then the containment condition becomes $E_i\subseteq V_j$, which is redundant from $V_j = E_j$.
\end{proof}

We record here a consequence of combining the previous two lemmas. Recall that $\pi:S\to R_n$ is the usual projection map.

\begin{lemma}
    \label{lem:defined by inclusions}
    If $w\in S_n(\calP)$, then $\pi^{-1}(J_w) = I_{\Delta_w}+I_n$.
\end{lemma}
\begin{proof}
    Suppose $w\in S_n(\calP)$. Combining \Cref{lem: A2 implies A3} and \Cref{lem: A3 implies A1}, we have that $J_w$ is generated by the squarefree monomials listed in \Cref{lem: A3 implies A1}. 
    %
    By inspection, these are exactly the monomials corresponding to the minimal nonfaces of $\Delta_w$, and so the result follows.
\end{proof}

Before continuing to complete the proofs of \Cref{thm:main-A} and \Cref{thm:main-B}, we count the exact number of minimal generators of $J_v$ for $v\in \sn(\calP)$.
\begin{proposition}
\label{prop: min gens of Iv}
    Let $v\in \sn(\calP)$ with $v\neq \mathrm{id}$, and let $i=i(w)$, $j=j(w)$, $p=p(w)$, and $q=q(w)$. 
    If $p \neq j$ or  $q\neq i$, then the minimal generating set of $I_{\Delta_v}$ given by the minimal nonfaces of $\Delta_v$ descends to a minimal generating set of $J_v$ in $R_n$. 
    If $p = j$ and $q= i$, then the images of the variables $x_1,\dots, x_{i-1}, x_{j+1},\dots, x_n$ in $R_n$  
    are a minimal generating set of $J_v$.
\end{proposition}
\begin{proof}
    By \Cref{lem:defined by inclusions}, the generating set of $I_{\Delta_v}$ descends to a generating set of $J_v$ so it suffices to show that they form a minimal generating set. 
    As there must be at least two variables that are not generators, $e_1$ cannot make any of the variable generators redundant.
    Thus, without loss of generality, we may assume that $i=1$ and $j=n$.
    
    If $p=j$ and $q=i$, then the only generator of $I_{\Delta_v}$ is $e_n$ so $J_v = (0)$, and the result holds in this case.

    If $p\neq j$ or $q\neq i$, it suffices to show that the images of the generators of $I_{\Delta_v}$ are nonredundant.
    Recall that any $w\in \calP$ is determined by the tuple $[i_w, j_w, p_w,q_w]$.
    Then $x_{[i,p]}$, $x_{[q,j]}$, or both is a minimal generating set.
    
    If $p\neq j$, then $x_{[i,p]}$ is a member of a minimal generating set of $I_{\Delta_v}$.
    If it was not a member of a minimal generating set of $J_v$, then $J_v \subseteq J_w$ where $w\in \calP$ is  determined by the tuple $[i_w, j_w, p_w,q_w] = [i_v, j_v, j_v,q_v]$.
    As $v\leq w$, $J_w\subset J_v$ as well by \Cref{prop:Inclusion} and \Cref{lem:defined by inclusions}, so $J_v=J_w$, which is impossible as $v\neq w$.
    Thus, $x_{[i,p]}$ is a member of a minimal generating set of $J_v$ (if it was a member of a minimal generating set of $I_{\Delta_v}$).

    A symmetric argument shows that $x_{[q,j]}$ is a member of a minimal generating set of $J_v$ (if it was a member of a minimal generating set of $I_{\Delta_v}$).
\end{proof}

\begin{remark}
\label{rem: number of min gens of Iv}
    It follows from the previous result that $J_w$ has a set of $n-j+i+1-\delta_{jp}-\delta_{iq}$ minimal generators in $R_n$. 
    For example, if $w = 3412,$ then $I_{\Delta_w}$ has two generators: $x_1x_2x_3$ and $x_2x_3x_4,$ since $i = 1, p = 3, q = 2,$ and $j = 4$ as in \Cref{def:Deltaw}. Thus, $n - j + i + 1 - \delta_{jp} - \delta_{iq} = n - m + 2 = 2,$ where $m := j - i + 1 = \#\supp(\Delta_w).$ One can also see that $J_w$ has a set of $n+\delta_{q<p}$ minimal generators in $R_n$ as follows.

    Note that any variable that is a generator of $J_w$ makes exactly one of the generators of $I_n$ redundant, so without loss of generality, we may assume that $J_v$ has no variable generators.
    If $J_w$ has a single generator, it is either $e_n$ or makes $e_n$ redundant (and no other generators of $I_n$ as each contains a monomial divisible by the generator) so the result holds in this case.  
    
    If $J_w$ has two minimal generators with $p<q$, then they make $e_{n-1}$ and $e_n$ redundant as any monomial missing at most one variable must be divisible by one or the other.
    No other generators of $I_n$ are made redundant as again they all contain a monomial not divisible by either generator (any monomial excluding one variable from each generator).
    Thus, the result holds in this case as well.
    
    If $J_w$ has two minimal generators with $p>q$, then they $e_n$ redundant, but do not make any other generator of $I_n$ redundant as again they all contain a monomial not divisible by either generator (any monomial excluding a variable from the overlapping variables). 
\end{remark}

We can now prove the first implication of \Cref{thm:main-B}.

\begin{proof}[Proof that \Cref{cond:B2} implies \Cref{cond:B1}]
   Let $\mathcal{W}\subseteq S_n(\mathcal{P})$ such that 
   \[\bigcap_{w\in\mathcal{W}} (I_{\Delta_w}+I_n) = \left(\bigcap_{w\in \mathcal{W}} I_{\Delta_w}\right)+I_n.\]
   Letting $\Delta_\mathcal{W} = \bigcup_{w\in\mathcal{W}}\Delta_w$, we then have 
   \begin{align*}
   \pi^{-1} \left( J_{\mathcal{W}} \right) 
   &= \pi^{-1}\left(\bigcap_{w\in\mathcal{W}} J_w\right)\\
   &= \bigcap_{w\in\mathcal{W}} \pi^{-1}(J_{w})\\
   &= \bigcap_{w\in\mathcal{W}} \left(I_{\Delta_w}+I_n\right) \text{ by \Cref{lem:defined by inclusions}}\\
   &= \left(\bigcap_{w\in\mathcal{W}} I_{\Delta_w}\right)+I_n \\
   & = I_{\Delta_\mathcal{W}}+I_n.
   \end{align*}
Applying the surjective map $\pi$ then gives 
\[J_\mathcal{W} = \pi\left(\pi^{-1}(J_{\mathcal{W}})\right) = \pi\left( I_{\Delta_\mathcal{W}}+I_n \right) = \overline{I_{\Delta_\mathcal{W}}},\]
so $J_\mathcal{W}$ is generated by squarefree monomials as desired.
\end{proof}

\subsection{Squarefree generation implies ideal equality}
We prove this implication of \Cref{thm:main-B} in two steps.
We first show the following.

\begin{proposition}
\label{prop:1-implies-avoidance}
If \Cref{cond:B1} holds, then $\mathcal W\subseteq S_n(\mathcal{P})$.
\end{proposition}

Following the proof of this proposition, we give that proof that \Cref{cond:B1} implies all of \Cref{cond:B2}.

\begin{proof}
    Let $\mathcal{W}\subseteq S_n$ be a set of permutations for which $J_{\mathcal{W}} = \bigcap_{w\in\mathcal W} J_w$ is generated by squarefree monomials as an ideal of $R_n$.
    Let $w\in\mathcal{W}$.
    Suppose for contradiction that $w\not\in S_n(\mathcal{P})$.

    Set \(\ell=\ell(w)\). 
    Recall to each permuation $w\in S_n$, we associate the set $A(w)=\mathrm{Des}_L(w)^c\cap[i+1,j-2]$.
    Using this, define a linear functional on \((R_n)_{\ell+1}\) by 
    \[  \Lambda_w:=\sum_{q\in A(w)}\lambda_{s_qw}. \]
    By \Cref{lem:AvoidEquivDecrease}, there exists some $q\in A(w)$ so this sum is non-empty.
    Fix a $q\in A(w)$.
    By definition of the $\lambda_*$, we have
    \[ \Lambda_w(\mathfrak{S}_{s_qw})=\sum_{q'\in A(w)}\lambda_{s_q'w}(\mathfrak{S}_{s_qw}) = \lambda_{s_qw}(\mathfrak{S}_{s_qw})=1. \]
    Thus, in order to find a contradiction, it suffices to show that we can also prove $\Lambda_w(\mathfrak{S}_{s_qw}) = 0$.

    In order to show that, we first need to show that $\mathfrak{S}_{s_qw}\in \mathrm{SqfMon}(J_w)$, which we do in this paragraph, and then show that $\Lambda_\mathcal{W}$ evaluates to $0$ on that entire ideal.
    If $s_qw\leq v\in\mathcal{W}$, then $w< s_qw\leq v$ which contradicts that $\mathcal{W}$ is an antichain (or the reflexivity of the strong Bruhat order).
    Thus, $s_qw\not\leq v$ for all $v\in \mathcal{W}$, so $\mathfrak{S}_{s_qw}\in J_{\mathcal{W}}$.
    As $J_{\mathcal{W}}\subseteq J_w$ and $J_{\mathcal{W}}$ is generated by squarefree monomials, $J_{\mathcal{W}}\subseteq \mathrm{SqfMon}(J_w)$.
    Thus, $\mathfrak{S}_{s_qw}\in \mathrm{SqfMon}(J_w)$.
    Thus, in order to find a contradiction by showing $\Lambda_w(\mathfrak{S}_{s_qw}) = 0$, it suffices to show that for all $f\in \mathrm{SqfMon}(J_w)$,  $\Lambda_{w}(f) = 0$.
    
    As $\mathrm{SqfMon}(J_w)$ is generated by squarefree monomials, $\mathrm{SqfMon}(J_w)_{\ell+1}$ is generated as a vector space by products of the form $hx_F$ where $x_F\in \left(\mathrm{SqfMon}(J_w)\right)_d \subseteq (J_w)_d$  and $h\in (R_n)_{\ell+1-d}$.
    Thus, to show that for all $f\in \mathrm{SqfMon}(J_w)$,  $\Lambda_{w}(f) = 0$, it suffices to show that $\Lambda_{w}(hx_F) = 0$ for all such pairs.

    By definition of $\Lambda_w$,  Remark \ref{rem:lin-func-length}, and \cite[Formula (2.7)]{Macdonald}, and Definition \ref{def:relative-diff-oper}, respectively, we have

    \begin{align}
    \label{eq: functional string of identities}
    \Lambda_w(hx_F)
    & = \sum_{q\in A(w)}\!\!\lambda_{s_qw}(hx_F) \\
    & = \sum_{q\in A(w)}\!\!\partial_{s_qw}(hx_F) \\
    & = \sum_{q\in A(w)}\!\!\partial_q\partial_w(hx_F) \\ 
    & = \sum_{q\in A(w)}\!\!\partial_q\sum_{r\leq w}r\!\left(\partial_{w/r}h\right)\partial_r(x_F)    
    \end{align}
    
    Now, $\deg\!\left(\partial_r(x_F)\right) = d-\ell(r)$ and \[ \deg\!\left(r(\partial_{w/r}h)\right) = \deg\!\left(\partial_{w/r}h\right) = \deg(h)-\bigl(\ell-\ell(r)\bigr) = \ell(r)+1-d.\] 
    In order for $r$ to possibly contribute a nonzero term, those must both be nonnegative so $d-1\leq \ell(r) \leq d$.

    When $\ell(r) =d$, $\partial_r(x_F) = \lambda_r(x_F)$ by Remark \ref{rem:lin-func-length}. 
    Since $x_F\in J_w$ and $J_w$ is spanned as a vector space by Schubert polynomials $\mathfrak{S}_{u}$ for $u\not\leq w$, then $\lambda_r(x_F) = 0$ for all $r\leq w$. 
    Thus, the terms in \eqref{eq: functional string of identities} with $\ell(r) =d$ contribute $0$, so \[\Lambda_w(hx_F) = \sum_{q\in A(w)}\!\!\partial_q\!\sum_{\substack{r\leq w\\\ell(r)=d-1}}\!\!\!r\left(\partial_{w/r}h\right)\partial_r(x_F).\]
    
    Then $\partial_{w/r}h$ has degree $\ell(r)-\ell+\ell+1-d = 0$, so it is a scalar, and hence is fixed by $r$, so the sums can be reordered as 
    \[
    \Lambda_w(hx_F) = \sum_{\substack{r\leq w\\\ell(r)=d-1}}\!\!\!r\left(\partial_{w/r}h\right)\sum_{q\in A(w)}\!\!\partial_q\partial_r(x_F)=0.
    \]
    Thus, it in order to show $\Lambda_{w}(hx_F) = 0$ for all relevant pairs $h$ and $x_F$, it  suffices to show $\sum_{q\in A(w)}\partial_q\partial_r(x_F)=0$ for all $\ell(r) = d-1$ and $r\leq w$. 
    If \(q'\in A(w)^C \subseteq [i+1,j-2]\), then \(s_{q'}w<w\), then \(\partial_{q'}\partial_w=0\) by \cite[Equation 2.7]{Macdonald} so it further suffices to show \[\displaystyle\sum_{q\in [i+1,j-2]}\partial_q\partial_r(x_F)=0\] when $\ell(r) = d-1$ and $r\leq w$ (and $x_F\in J_w$ is a squarefree monomial of degree $d$).

    Fix $r\leq w$ with $\ell(r) = d-1$ and $x_F \in (J_w)_d$ a squarefree monomial.
    Writing $r$ as a reduced word and repeatedly applying Remark \ref{rem:sqfr-mon-div-diff} gives that $\partial_r(x_F) = c_F x_a$ for some $c_F\in \{0,\pm 1\}$ and some $a\in [n]$.
    If $c_F = 0$, we have reached our conclusion so assume without loss of generality that $c_F \neq 0$.

    If $s_{a-1}$ or $s_{a}$ are in any reduced word expression for $r$, then  $c_F = 0$ as either of those would have eliminated that variable $x_a$ from $x_{F}$ either by dividing it out or by zeroing it out.
    Thus, we have 
    \begin{equation}
    \label{eq:two-terms}
     \partial_q\partial_r(x_F)  =
     \begin{cases}
     -c_F,&q=a-1,\\
     c_F,&q=a,\\
     0,&\text{otherwise}.
     \end{cases}
     \end{equation}
In other words, the sum simplifies to 
\[\displaystyle\sum_{q\in [i+1,j-2]}\partial_q\partial_r(x_F)=\displaystyle\sum_{q\in [a-1,a]\cap [i+1,j-2]}\partial_q\partial_r(x_F).\]
If $[a-1,a]\cap [i+1,j-2] = \varnothing$, then the sum is trivially zero as desired.
If $[a-1,a]\cap [i+1,j-2] = [a-1,a]$, then $1$ and $-1$ each appear once in the sum so it is zero as desired.

In order to show that the sum is $0$, it now suffices to show that if $|[a-1,a]\cap [i+1,j-2]| = 1$, then there is a contradiction.
Suppose for contradiction that $|[a-1,a]\cap [i+1,j-2]| = 1$.
Then either $a=i+1$ or $a-1 = j-2$; in the second case, $a=j-1$.
Set $b = i = a-1$ or $b =j-1=a$, respectively.
Let $ \mathrm{supp}(v)$ denote the set of $i$ such that $s_i$ appears in some/every reduced word giving $v$.
Since \(r\leq w\),
\[\mathrm{supp}(r)\subseteq \mathrm{supp}(w)\subseteq\{i,\ldots,j-1\}. \]
As $c_F\neq 0$, neither $a$ nor $a-1$ is in $\mathrm{supp}(r)$.
So in our two cases, $i,i+1\not\in \mathrm{supp}(r)$ and $j-2,j-1\not\in \mathrm{supp}(r)$, respectively.
In either case, $s_b$ commutes with every reduced word for $r$.

By commutativity and the fact that $b\not\in\mathrm{supp}(r)$, appending $s_b$ to a reduced word for $r$ gives a reduced word for $s_br$.
As $b\in\{i,j-1\}$ and $w$ fixes every $k\in[i,j]$ but neither $i$ nor $j,$ the permutation $w$ does not fix $[1,b],$ so $s_b$ occurs in every reduced word of $w.$ Fix a reduced word for $w.$ Since $r\leq w,$ by the subword property \cite[Corollary 2.2.3]{BjornerBrenti} some subword of it is a reduced word for $r,$ and as $b\not\in\supp(r),$ this subword uses no occurrence of $s_b.$ If we include one occurrence of $s_b$ into the subword, then we obtain a subword of the reduced word for $w,$ whose product is $s_br,$ since $s_b$ commutes with every letter of $r,$ and which is reduced by the previous sentence. Hence, $s_b\leq w.$

Recall that the Schubert polynomials form a basis for $R_n$ and the Schubert polynomials for $v\not\leq w$ form a basis of $J_w$.
Because \(x_F\in J_w\), it can be written in terms of that subbasis (so its coefficients on the other Schubert polynomials are $0$).
Thus,  \(\partial_{s_br}(x_F)=\lambda_{s_br}(x_F)=0\).  
In contrast, by \eqref{eq:two-terms}
\[  
\partial_{s_br}(x_F)=\partial_b\partial_r(x_F) =\pm c_F. 
\] 
This contradiction rules out the case of exactly one index lying in the set, so the sum is zero as desired.
As the sum is $0$, $\Lambda_{w}(\mathrm{SqfMon}(I))=0$.
As $\mathfrak{S}_{s_qw} \in \mathrm{SqfMon}(J_w)$, we have $\Lambda_{w}(\mathfrak{S}_{s_qw}) = 0 = 1$, and we have the contradiction we needed.
\end{proof}

We need one more lemma before proving the second implication 

\begin{proposition}
\label{prop:lift}
If \(w\in S_n(\mathcal{P})\) and \(x_F\) is squarefree, then
\[ 
\overline{x_F}\in J_w \quad\Longleftrightarrow\quad  x_F\in I_{\Delta_w}, 
\]
where $\overline{x_F}$ is the image of $x_F$ under the projection $\pi$.
\end{proposition}

\begin{proof}
If $w =\mathrm{id}$, then $I_{\Delta_w}$ and $J_w$ are each $(x_1,\dots,x_n)$ in their respective rings, so the result is immediate as they both contain all $x_F$ and do not contain $1$.
For the rest of the proof, let $w\neq \mathrm{id}$.

The reverse implication is also immediate since $\pi^{-1}(J_w)=I_{\Delta_w}+I_n$ by \Cref{lem:defined by inclusions}.  

For the forward implication, we work by contrapositive, so suppose that \(x_F\notin I_{\Delta_w}\).  
By definition of $I_{\Delta_w}$, we have
\[ F\subsetneq[i,j],\qquad [i,p]\nsubseteq F,\qquad [q,j]\nsubseteq F.\]
Without loss of generality, mod out by the variables $x_a$ for  \(a\in [i,j]^c\subset[n]\).  
By \Cref{lem:quotient-and-annihilators}, Part (1), the quotient is \(R_m\), where \(m=j-i+1\), and the two remaining generators\footnote{After relabeling $x_{k+i-1} = y_{k}$} of \(J_w\) are
\[  
y_G=y_1\cdots y_{p-i+1}\qquad\text{and}\qquad y_H=y_{q-i+1}\cdots y_m.
\]
Here \(p>i\) and \(q<j\), so both degrees are at least two.  
Further, as $p\neq q$, their overlap does not have size 1.
In terms of monomials, the three relations above become  
\[
y_F\mid y_{[1,m]},\qquad y_G\nmid y_F,\qquad y_H\nmid y_F,\qquad |G\cap H|\neq 1.
\]

By \Cref{lem:contraction}, $\overline{y_F}\not\in (\overline{y_G},\overline{y_H}) = J_w \subseteq R_m$.
Lifting back from the quotient, the contrapositive follows.

\end{proof}

Given \Cref{prop:1-implies-avoidance} and \Cref{prop:lift}, we can now prove the second implication of the main result. From here on, if $I$ is an ideal in $S$, we will use $\overline{I}$ to denote the corresponding ideal $\pi(I)$ in $R_n$.

\begin{proof}[Proof that \Cref{cond:B1} and $\mathcal{W}\subseteq S_n(\mathcal{P})$ implies \Cref{cond:B2}.]
Assume that $J_{\mathcal{W}}$ is generated by squarefree monomials, i.e., \(J_\mathcal{W} = \{\overline{m_1},\dots,\overline{m_r}\}\subseteq R_n\) where \(\{m_1,\dots,m_r\}\subset \mathbb{C}[x_1,\dots,x_n]\) are squarefree monomials, and $\mathcal{W}\subseteq S_n(\mathcal{P})$.

By Proposition \ref{prop:lift}, \(m_a\in I_{\Delta_w}\) for all \(w\in\mathcal{W}\) and \(a\in[r]\), so \(m_a \in  \bigcap_{w\in \mathcal{W}} I_{\Delta_w}\).
Thus, \(\overline{m_a}\in \overline{\left(\bigcap_{w\in \mathcal{W}}  I_{\Delta_w}\right)}\).
As all of the generators of $J_{\mathcal{W}}$ are contained in $\overline{\left(\bigcap_{w\in \mathcal{W}}  I_{\Delta_w}\right)}$, we have the inclusion $J_{\mathcal{W}}\subseteq \overline{\left(\bigcap_{w\in \mathcal{W}}  I_{\Delta_w}\right)}$.
Using that inclusion and that \(\overline{I_{\Delta_w}}=J_w\) by \Cref{lem:defined by inclusions}, we have the following string of inclusions:
\[ J_{\mathcal W} \subseteq \overline{\left(\bigcap_{w\in\mathcal W}I_{\Delta_w}\right)} \subseteq \bigcap_{w\in\mathcal W}\overline{I_{\Delta_w} }=\bigcap_{w\in\mathcal W}J_{w} = J_{\mathcal{W}}.\]
Thus, all of these are equal.

Finally, as $\pi^{-1}(J_w) = I_{\Delta_w}+I_n$ by \Cref{lem:defined by inclusions}, we have
\[
\begin{aligned}
\bigcap_{w\in\mathcal W}(I_{\Delta_w}+I_n)
&=\bigcap_{w\in\mathcal W}\pi^{-1}(J_w) =\pi^{-1}\!\left(\bigcap_{w\in\mathcal W}J_w\right)
=\pi^{-1}(J_{\mathcal W})\\
&=\pi^{-1}\left(\overline{\bigcap_{w\in\mathcal W}I_{\Delta_w}}\right)
=\left(\bigcap_{w\in\mathcal W}I_{\Delta_w}\right)+I_n.
\end{aligned}
\] 

Thus, \Cref{cond:B2} holds as desired.
\end{proof}

\subsection{Ideal equality implies prunable}
In this section, we prove the third implication of \Cref{thm:main-B} by induction on the length of $\mathcal{W}$.
As both \Cref{cond:B1,cond:B2} require the set to be contained in $S_n(\calP)$, and all of their other parts are trivial for singletons, that base case is trivial.
However, in the proof we need the additional base case of when the set has size 2.

\subsubsection{Ideal equality implies prunable for pairs.}
For pairs, \Cref{cond:B3}(b) is trivial, so it suffices to show \Cref{cond:B3}(a).

We first prove a more general lemma about intersections of ideals of simplicial complexes plus $I_n$ and the dimension of the intersection of the complexes.
When applied in our setting, this gives the dimension aspect of (3)(a). 
\begin{lemma}
\label{lem:dimension-complete}
If \(\Delta_{1},\Delta_{2}\) are simplicial complexes such that $(I_{\Delta_{1}}+I_n)\cap(I_{\Delta_{2}}+I_n)  =(I_{\Delta_{1}}\cap I_{\Delta_{2}}) +I_n$, then $ \dim(\Delta_{1}\cap\Delta_{2})  =\min\{\dim\Delta_{1},\dim\Delta_{2}\}$.
\end{lemma}

\begin{proof}
  Start with the exact sequence
    \[
    0\longrightarrow S/I_{\Delta_{1}\cup\Delta_{2}} \longrightarrow S/I_{\Delta_{1}}\oplus S/I_{\Delta_{2}} \longrightarrow S/I_{\Delta_{1}\cap\Delta_{2}} \longrightarrow0 
    \]
where $S$ is the polynomial ring.
Assume for contradiction that 
\[
d=\dim(\Delta_{1}\cap\Delta_{2}) <\min\{\dim\Delta_{1},\dim\Delta_{2}\}.
\]
Since \(\Delta_{1}\cap\Delta_{2}\) has no face of dimension \(d+1\), the element \(e_{d+2}\) vanishes in \(S/I_{\Delta_{1}\cap\Delta_{2}}\). 
Therefore, \((e_{d+2},0)\) is in the kernel of the surjective map.
By exactness, there exists an element \(\xi\in S/I_{\Delta_{1}\cup\Delta_{2}}\) mapping to \((e_{d+2},0)\).
Note, \((e_{d+2},0)\) is nonzero as $\Delta_{1}$ has a face of dimension at least $d+1$.

By the assumption on the ideal, we can tensor by $R_n$ to get the sequence
\[0\longrightarrow R_n/I_{\Delta_{1}\cup\Delta_{2}} \longrightarrow R_n/I_{\Delta_{1}}\oplus R_n/I_{\Delta_{2}} \longrightarrow R_n/I_{\Delta_{1}\cap\Delta_{2}} \longrightarrow0 \]

In the quotient sequence, we have \(\overline{(e_{d+2},0)} = \overline{(0,0)}\), so by injectivity \(\overline{\xi} =\overline{0}\).
Thus, \(\xi\in \overline{I_n}\subseteq S/I_{\Delta_{1}\cup\Delta_{2}}\), so we can write
\[ \xi=\sum_{j=1}^{d+1}e_jh_j+c e_{d+2} \]
where $c$ is a scalar.
Taking the two coordinate images in the original short exact sequence gives
\[ e_{d+2}=\sum_{j=1}^{d+1}e_jh_j+c e_{d+2} \quad\text{in }S/I_{\Delta_{1}}\qquad\text{and}\qquad 0=\sum_{j=1}^{d+1}e_jh_j+c e_{d+2} \quad\text{in }S/I_{\Delta_{2}}. \]

As $\Delta_{1}$ has a $d+1$-dimensional, $e_{d+2}$ is not in the ideal generated by the smaller elementary symmetric polynomials\footnote{Quotient out by the variables outside the the $d+1$ dimensional face and this becomes the standard statement on the ring $R_m$}.
Thus, by the 1st equality, \(c\neq 0\in S/I_{\Delta_1}\).
Applying the same logic to the 2nd equality, gives \(c=0\in S/I_{\Delta_2}\).
As $c$ was a scalar, those two statements are contradictory, and so the result follows.
\end{proof}

Next, we prove the purity aspect of (3)(a).
\begin{lemma}
\label{lem:crossing-class-complete}
Let \(\{u,v\}\subseteq S_n(P)\) be an antichain for which $(I_{\Delta_u}\cap I_{\Delta_v}) +I_n = (I_{\Delta_u}+I_n)\cap (I_{\Delta_v}+I_n)$.
Then \(\Delta_u\cap\Delta_v\) is pure.
\end{lemma}

\begin{proof}
    We work by contradiction so assume that \(\Delta_u\cap\Delta_v\) is not pure.
    As the intersection of doubly pruned complexes is a doubly pruned complex, \(\Delta_u\cap\Delta_v\) is doubly pruned, i.e. its minimal nonfaces are vertices outside of some interval $[i,j]$ and, possibly, two interval nonfaces of the form $[i,p]$ and $[q,j]$.
    Since it is not pure, $q\leq p$.
    As $\{u,v\}$ is an antichain, assume without loss of generality that $[i,p]$ is a nonface of $\Delta_u$ and a face of $\Delta_v$ and $[q,j]$ is a nonface of $\Delta_v$ and a face of $\Delta_u$.
    
    Define the permutation $w$ as follows: 
    \[
    w :=(1,\dots,i-1,i+1,\dots,q-1,j,q,\dots,p,i,p+1,\dots,j-1,j+1,\dots,n).
    \]
    We then have the rank inequality 
    \[
    r_w(p,i) = i-1 <i = r_u(p,i)
    \]
    so $w\not\leq u$.
    Similarly, we have 
    \[
    r_w(q-1,j-1)=q-2<q-1 = r_v(q-1,j-1)
    \]
    so $w\not\leq v$.
    Since $w\not\leq u,v$, it follows that $\mathfrak{S}_w\in J_u\cap J_v$ in $R_n$.
    
    By \Cref{lem:defined by inclusions}, $\pi^{-1}(J_u) = I_{\Delta_u}+I_n$ and $\pi^{-1}(J_v) = I_{\Delta_v}+I_n,$ so by our hypothesis,
    \[
    \pi^{-1}(J_u\cap J_v) \subseteq \pi^{-1}(J_u)\cap \pi^{-1}(J_v)
    = (I_{\Delta_u}+I_n)\cap(I_{\Delta_v}+I_n) 
    = (I_{\Delta_u}+I_{\Delta_v}) + I_n
    \]
    and in particular any lift of $\fraks_w$ to $S$ lies in $(I_{\Delta_u}\cap I_{\Delta_v})+I_n.$ 
    It now suffices to show that $\fraks_w\not\in\pi\left((I_{\Delta_u}\cap I_{\Delta_v})+I_n\right),$ i.e., that no lift of $\fraks_w$ lies in $(I_{\Delta_u}\cap I_{\Delta_v})+I_n.$
    Consider the permutation
    \[
    t := (1,\dots,i-1, 
    j-1,\dots,j-q+i+1,
    j, 
    j-q+i,\dots,j-p+i, 
    i, j-p+i-1,
    \dots,i+1,
    j+1,\ldots,n).
    \]
    Notice that $t\in S_n(\calP)$ is the unique strong Bruhat maximal element among any of the permutations $\sigma$ with $\sigma^{-1}(i)=p+1$, $\sigma^{-1}(j)=q-1$, such that $\sigma(k) = k$ for all $k\not\in[i,j]$. By \Cref{prop: Ehresmann's criterion}, placing the remaining values $i+1,\ldots,j-1$ in decreasing order in the remaining positions minimizes every entry of the rank matrix among such permutations, and $t\in\sn(\calP)$ by \Cref{lem:AvoidEquivDecrease}. 
    Since $w$ satisfies those conditions, $w\leq t$, and so $\fraks_w\not\in J_t,$ i.e., no lift of $\fraks_w$ lies in $\pi^{-1}(J_t) = I_{\Delta_t}+I_n.$ Thus, it suffices to show that $(I_{\Delta_u}+I_{\Delta_v})+I_n\subseteq I_{\Delta_t}+I_n,$ and for this it suffices to show that $I_{\Delta_u}+I_{\Delta_v}\subseteq I_{\Delta_t}.$
    
    As the vertex nonfaces of $u$ and $v$ are clearly both contained in the vertex nonfaces of $t$, we only need to show that the 4 least common multiples of nonsingleton generators in $I_{\Delta_u}$ and $I_{\Delta_v}$ are in $I_{\Delta_t}$ (recall that $I_{\Delta_u}\cap I_{\Delta_v}$ is generated by the least common multiple of pairs of generators, one from each ideal).

    The left end interval of $v$ either contains a vertex nonface of $t$ or contains the nonface $[i,p+1]$ so any least common multiple using that face is in $I_{\Delta_t}.$
    Similarly, the right end interval of $u$ either contains vertex nonface of $t$ or contains the nonface $[q-1,j]$ so any least common multiple using that face is in $I_{\Delta_t}$.
    Finally, the least common multiple of the generators coming from nonfaces $[i,p]$ and $[q,j]$ is the monomial $x_{[i,j]}$ which is clearly in $I_{\Delta_t}$.
    Thus, we have the needed inclusion, and the result follows.
\end{proof}

\begin{lemma}
\label{lem:2 implies 3 for pairs}
    If $\mathcal{W}\subseteq S_n$ satisfies \Cref{cond:B2} and $|\mathcal{W}|=2$, then it satisfies \Cref{cond:B3} as well.
\end{lemma}

\begin{proof}
    The previous two lemmas prove \Cref{cond:B3}(a). \Cref{cond:B3}(b) is trivial for pairs.
\end{proof}

\subsubsection{The proof for the general case}
With the case of pairs established, we move onto the general proof.

\begin{proof}[Proof that \Cref{cond:B2} implies \Cref{cond:B3}.]
    Let $\mathcal{W}= \{w_1,\dots,w_m\} \subseteq S_n$ satisfy \Cref{cond:B2}. 
    Reorder $\mathcal{W}$ lexicographically by the tuples $[i(w),j(w),-p(w),-q(w)]$ as defined in \Cref{def:Deltaw} and then relabel so that $w_i<w_{i+1}$ in that order.
    We work by induction on $m$.  

    The base case $m=1$ is trivial, and the base case $m=2$ is \Cref{lem:2 implies 3 for pairs}.
    
    Now assume that \Cref{cond:B2} implies \Cref{cond:B3} for all subsets of size at most $m-1$ where $m\geq 3$. 

    In order to utilize the induction hypothesis, we first want to show that $\{w_1,\dots,w_{m-1}\}$ satisfies \Cref{cond:B2}.
    A symmetric argument replacing each index $i$ with $m+1-i$ shows that $\{w_2,\dots,w_m\}$ also satisfies \Cref{cond:B2}.
    Putting those together with a simple induction argument, we see that $\{w_a,w_{a+1},\dots,w_{b-1},w_b\}$ satisfies \Cref{cond:B2} (and hence \Cref{cond:B3} by induction) for all $1\leq a<b\leq m$.

    Now to show that $\{w_1,\dots,w_{m-1}\}$ satisfies \Cref{cond:B2}, we must show that 
    \[
    \left(\bigcap_{i=1}^{m-1} I_{\Delta_{w_i}}\right)+I_n=\bigcap_{i=1}^{m-1} \left(I_{\Delta_{w_i}}+I_n\right)
    \]
    in the polynomial ring.
    As the kernel of the projection from the polynomial ring to $R_n$ is $I_n$ which is contained in both sides, showing that equality is equivalent to showing  
    \[
    \overline{\bigcap_{i=1}^{m-1} I_{\Delta_{w_i}}}=\bigcap_{i=1}^{m-1} \overline{I_{\Delta_{w_i}}}=\bigcap_{i=1}^{m-1}J_{w_i}.
    \]
    By hypothesis, 
    \[
    \overline{\bigcap_{i=1}^{m-1} I_{\Delta_{w_i}}}\supseteq \overline{\bigcap_{i=1}^{m} I_{\Delta_{w_i}}}=\bigcap_{i=1}^{m} J_{w_i}\subseteq \bigcap_{i=1}^{m-1} J_{w_i},
    \]
    so showing the desired equality is equivalent to showing
    \[\frac{\left(\displaystyle
\overline{\bigcap_{i=1}^{m-1}I_{\Delta_{w_i}}}\right)}
{\left(\displaystyle
\overline{\bigcap_{i=1}^{m}I_{\Delta_{w_i}}}\right)}
=
\frac{\displaystyle\bigcap_{i=1}^{m-1}J_{w_i}}
{\displaystyle\bigcap_{i=1}^{m}J_{w_i}}.\]
    As the identity is an inclusion $\bigcap_{i=1}^{m-1} I_{\Delta_{w_i}}\hookrightarrow I_{\Delta_{w_i}}$ for all $i$, it is also an injection $\overline{\bigcap_{i=1}^{m-1} I_{\Delta_{w_i}}}\hookrightarrow \bigcap_{i=1}^{m-1} \overline{I_{\Delta_{w_i}}}$, which in turn induces an injection 
    \[\frac{\left(~\displaystyle
\overline{\bigcap_{i=1}^{m-1}I_{\Delta_{w_i}}}~\right)}
{\left(~\displaystyle
\overline{\bigcap_{i=1}^{m}I_{\Delta_{w_i}}}~\right)}
\lhook\joinrel\xrightarrow{\hspace{.2in}}
\frac{\displaystyle\bigcap_{i=1}^{m-1}J_{w_i}}
{\displaystyle\bigcap_{i=1}^{m}J_{w_i}}.\]
    Thus, showing the equality is equivalent to showing that this map is surjective.
    Let $\overline{f}$ be a class in the codomain.
    As $I_{\Delta_{w_{m-1}}}$ surjects onto the numerator of the codomain, we may assume without loss of generality that $f \in I_{\Delta_{w_{m-1}}}$.
    As $I_{\Delta_{w_{m-1}}}$ is generated by squarefree monomials, if $x_H^{\vec{a}} = x_{h_1}^{a_1}\dots x_{h_K}^{a_k}$ is a monomial which appears in $f$ with nonzero coefficient, then $H\not\in \Delta_{w_{m-1}}$. 
    Let $g$ be the polynomial gotten from $f$ by setting all of the coefficients of monomials $x_H^{\vec{a}}$ to $0$ if $H \not\in \Delta_{w_m}$.
    For every $x_H^{\vec{a}}$ appearing in $g$ with nonzero coefficient, we then have $H \not\in \Delta_{w_{m-1}}$ and $H\in\Delta_{w_m}$.
    As $\Delta_{w_k}\cap \Delta_{w_m}\subseteq \Delta_{w_{m-1}}\cap \Delta_{w_m}$ by \Cref{lem:antichain ordering},  $H\not\in \Delta_{w_k}$ for all $k\in [1,m-1]$.
    Thus, $x_H\in \bigcap_{i=1}^{m-1}I_{\Delta_{w_i}}$ and therefore so is $g$.
    It follows that $\overline{g}\in \bigcap_{i=1}^{m-1}J_{w_i}$.
    As that is an ideal and (the image of) $\overline{f}$ is also in the ideal, $\overline{f-g} = \overline{f}-\overline{g}\in \bigcap_{i=1}^{m-1}J_{w_i}$.
    If $x_H^{\vec{a}}$ has nonzero coefficient in $\overline{f-g}$, by construction $H\not\in \Delta_{w_m}$.
    For such $H$, we have $x_H \in I_{\Delta_{w_m}}$.
    Thus $f-g\in I_{\Delta_{w_m}}$, so $\overline{f-g}\in\bigcap_{i=1}^m J_{w_i}$ so $\overline{f}$ and $\overline{g}$ represent the same class in the codomain, and by injectivity, in the domain.
    As $g \in \bigcap_{i=1}^{m-1}I_{\Delta_{w_i}}$ and maps to $\overline{f}$, the map is surjective. 
    Thus, we have the needed equality to show that $\{w_1,\dots,w_{m-1}\}$ satisfies \Cref{cond:B2} as desired.    

    Assume that $\dim(\Delta_{w_1})\geq \dim(\Delta_{w_m})$; if not, there is an entirely symmetric argument to what follows after replacing the indices $i$ with $m+1-i$. 
    As $\{w_1,\dots,w_{m-1}\}$ satisfies \Cref{cond:B2}, it satisfies \Cref{cond:B3} by induction. 
    Thus, it suffices to prove the parts of \Cref{cond:B3} involving $w_m$.
    (3)(b) of \Cref{cond:B3} is automatic for $w_m$ by \Cref{lem:antichain ordering}.
    Thus, it suffices to show (3)(a), i.e.,  that the intersection of $\Delta_{w_m}$ with the union of $\Delta_{w_1}$ through $\Delta_{w_{m-1}}$ is pure of dimension $\dim(\Delta_{w_m})$. 

    Let $\Gamma = \bigcup_{i=1}^{m-1}\Delta_{w_{i}}$. 
    By the argument above, $\{w_{m-1},w_m\}$ satisfies \Cref{cond:B2} and, hence, \Cref{cond:B3}.
    In particular, $\Delta_{w_m}\cap \Delta_{w_{m-1}}$ is pure (of dimension $\min(\dim(\Delta_{w_m}),\dim(\Delta_{w_{m-1}}))$).
    As $\Delta_{w_m}\cap \Gamma =\Delta_{w_m}\cap \Delta_{w_{m-1}}$, the same statement holds for it.
    $\Delta_{w_m}$ and $\Gamma$ satisfy the hypothesis of \Cref{lem:dimension-complete} by the following string of equalities which follow by standard facts about squarefree monomial ideals and the fact that $\{w_a,w_{a+1},\dots,w_{b-1},w_b\}$ satisfies \Cref{cond:B2} for all $1\leq a<b\leq m$, 
    \begin{align*}
    (I_{\Gamma}+I_n) \cap (I_{\Delta_{w_m}}+I_n) 
    &= \left(I_{\bigcup_{i=1}^{m-1}\Delta_{w_i}}+I_n\right) \cap (I_{\Delta_{w_m}}+I_n)\\
    &= \left(\left(\bigcap_{i=1}^{m-1}I_{\Delta_{w_i}}\right)+I_n\right) \cap (I_{\Delta_{w_m}}+I_n)\\
    &= \left(\bigcap_{i=1}^{m-1}\left(I_{\Delta_{w_i}}+I_n\right)\right) \cap (I_{\Delta_{w_m}}+I_n)\\
    &= \bigcap_{i=1}^{m}\left(I_{\Delta_{w_i}}+I_n\right)\\
    &= \left(\bigcap_{i=1}^{m}I_{\Delta_{w_i}}\right)+I_n\\
    &= \left(\left(\bigcap_{i=1}^{m-1}I_{\Delta_{w_i}}\right)\cap I_{\Delta_{w_m}}\right)+I_n\\
    &= \left(I_{\bigcup_{i=1}^{m-1}\Delta_{w_i}}\cap I_{\Delta_{w_m}}\right)+I_n\\
    &=(I_{\Gamma}\cap I_{\Delta_{w_m}})+I_n,
    \end{align*}
    and so  $\dim\left(\Delta_{w_m}\cap \Gamma\right) = \min(\dim(\Delta_{w_m}),  \dim(\Gamma))$ by \Cref{lem:dimension-complete}.

    As $\dim(\Gamma) \geq \dim(\Delta_{w_1})\geq \dim(\Delta_{w_m})$ by assumption, we have \[\dim\left(\Delta_{w_m}\cap \Gamma\right) = \min(\dim(\Delta_{w_m}),  \dim(\Gamma)) =\dim(\Delta_{w_m}).\]
    Thus, \Cref{cond:B3} holds for $\mathcal{W}$ as desired, and so the result follows by induction.
\end{proof}

\subsection{Prunable implies ideal equality}
\label{subsec:Prune_Equality}

We move onto the final implication of \Cref{thm:main-B}. 
We need four algebraic lemmas before our proof.
The first lemma relates the generators of $\Delta_w$ and the ideal of its intersection with $\Delta_v$.

\begin{lemma}\label{lem:mumbo-generator-jumbo} 
    Let $\Delta_v,\Delta_w$ be incomparable doubly-pruned complexes on $[n]$. Let $\Gamma = \Delta_v\cap \Delta_w$. If $\Gamma$ is a doubly-pruned complex that is pure of dimension $\dim(\Delta_w)$, then $I_{\Gamma} = (g_1,\dots, g_r)$ and $I_{\Delta_w} = (g_1,\dots, g_{r-1},x_B g_r)$ for some squarefree interval monomials $g_i$ such that $x_B$ is an interval squarefree monomial adjoining $g_r$ on one side.
\end{lemma}

Note, the proof of this lemma uses a few results from the next section, but they are completely independent of this entire section.

\begin{proof}   
    Suppose $\Gamma = \Delta_v\cap \Delta_w$ is doubly pruned and pure of dimension $\dim(\Delta_w)$. 
    By \Cref{lem: doubly pruned complexes are Schubert complexes}, \Cref{cor: purity of doubly pruned complexes}, and \Cref{rem: doubly pruned purity}, either $\Gamma = \Delta_u$ for some $u\in S_n(\calP)$ or $\Gamma = \mathscr{P}([i,j])$ for some $1\leq i \leq j \leq n$.
    We deal with each of these cases separately.

    First, assume that $\Gamma = \mathscr{P}([i,j])$ for some $1\leq i \leq j \leq n$, so in particular, it has dimension $j-i.$
    The complexes $\Delta_v$ and $\Delta_w$ have vertex nonfaces outside of intervals $[i_v,j_v]$ and $[i_w,j_w]$, and any vertex nonface of $\Gamma$ is a vertex nonface of $v$ or $w$.
    Thus, $[i,j] =[\max(i_v,i_w),\min(j_v,j_w)]$.
    If $[i,j] = [i_v,j_v]$, then as $\Delta_v$ has some nonvertex nonface, $\Gamma$ would also have some nonvertex nonface (or possibly smaller).
    Thus, $[i,j] \neq [i_v,j_v]$, and by symmetry, $[i,j] \neq [i_w,j_w]$.
    Up to symmetry, we may then assume that $i = i_w >i_v$ and $j = j_v <j_w$.
    Since $\Delta_w$ is pure and the same dimension as $\Gamma$, then either $[i_w,j_w] = [i,j+1]$ and its only nonvertex minimal nonface is $[q,j+1]$ or $[i_w,j_w] = [i,j+2]$ and its nonvertex minimal nonfaces are $[i,p]$ and $[q,j+2]$ with $p<q$.
    If $[i_w,j_w] = [i,j+2]$, then as $[i,p]$ is not a minimal nonface of $\Gamma$, $j<p<q<j+2$, which is a contradiction.
    If $[i_w,j_w] = [i,j+1]$, then the only nonface of $\Delta_w$ and $\Gamma$ that are different are the nonface $j+1 \in \Gamma$ and the nonface $[q,j+1]$ so the result follows in this case.
    
    Thus, we may assume that $\Gamma = \Delta_u$ for some $u\in S_n(\calP)$.
    By assumption, $\Delta_u$ is pure, and by \Cref{lem:pure-intersection}, $\Delta_w$ is pure. 
    Hence, for both $u$ and $w$, we either have $p<q$, $p =j$, or $q=i$.
    As $\Delta_u \subset \Delta_w$, $[i(u),j(u)]\subseteq [i(w),j(w)]$. 
    Similarly, as $\Delta_u \subset \Delta_v$, $[i(u),j(u)]\subseteq [i(v),j(v)]$.
    The rest of the proof is split into the nine cases define by the three possibilities for $u$ and $w$

    \medskip

    \noindent\textbf{Case 1:} $p(u)<q(u)$ and $p(w)<q(w)$.

    \medskip
    
    In this case,  $j(u)-i(u)-2 = \dim(\Delta_u) = \dim(\Delta_w) = j(w)-i(w)-2$, so $[i(u),j(u)] = [i(w),j(w)]$.
    We denote this interval as $J$ for the rest of this case.
    If $J \subseteq [i(v)+1,j(v)-1]$, then $\Delta_w\subset \Delta_v$, which is impossible by our assumption that they are incomparable. 
    Thus, either $J \subseteq [i(v),j(v)-1]$, $J \subseteq [i(v)+1,j(v)]$, or $J=[i(v),j(v)]$.
    
    If $J \subseteq [i(v),j(v)-1]$, then the vertex $j(v)$ is a minimal nonface of $\Delta_u$ and $\Delta_w$ but is a face of $\Delta_v$ so $j(u)=j(w)<j(v)$.
    If $i(v)\not\in J$, then $\Delta_w = \Delta_w$, which violates the incomparable hypothesis, so $i(v) = i(u) = i(w)$.
    Thus, $[q(v),j(v)]$ is a nonface of $\Delta_u$ that is not minimal.
    As a result, $q(u)=q(w)\geq i(w),i(v)$.
    As $i(u)=i(v)=i(w)$, $p(u) =\min(p(v),p(w))$.
    If $p(w)\leq p(v)$, then $\Delta_w\subset \Delta_v$, which is impossible again by assumption, so $p(w)>p(v)$.
    Thus, $p(u) = p(v)<p(w)\leq j(u)$.
    The minimal non-vertex nonfaces of $u$ are $[i(u),p(u)] = [i(v),p(v)]$ and $[q(u),j(u)]= [q(w),j(w)]$.
    Thus, a set of minimal generators of $I_{\Delta_w}$ is
    \[x_1,\dots,x_{i(u)-1},\qquad x_{j(u)+1},\dots,x_n,\qquad x_{[q(u),j(u)]},\qquad x_{[i(w),p(w)]}.\]
    As $x_{[i(w),p(w)]} = x_{[i(w),p(v)]}\cdot x_{[p(v)+1,p(w)]}$, the result follows in this case.
    
    The case where $[i(w),j(w)]= [i(u),j(u)] \subseteq [i(v)+1,j(v)]$ is similar.

    We are left with the case  $[i(w),j(w)]= [i(u),j(u)] = [i(v),j(v)]$, it follows that $p(u)=\min(p(v),p(w))< \max(p(v),p(w))$ and $\min(q(v),q(w))<q(u) = \max(q(v),q(w))$.
    If \newline $\min(p(v),p(w)) = p(v)$ and $\max(q(v),q(w)) = q(v)$, then $\Delta_u=\Delta_v$ which impossible.
    The symmetric statement holds for $w$.
    Thus, ($p(u) = p(v)<p(w)$ and $q(v)<q(w)=q(u)$) or ($p(u) = p(w)<p(v)$ and $q(w)<q(v)=q(u)$)
    If  $p(u) = p(v)<p(w)$ and $q(v)<q(w)=q(u)$, then the generators of $I_{\Delta_w}$ are
    \[x_1,\dots,x_{i(u)-1},\qquad x_{j(u)+1},\dots,x_n,\qquad x_{[q(u),j(u)]},\qquad x_{[i(w),p(w)]}.\]
    Again as $x_{[i(w),p(w)]} = x_{[i(w),p(v)]}\cdot x_{[p(v)+1,p(w)]}$, the result follows in this case.
    The case where $p(u) = p(w)<p(v)$ and $q(w)<q(v)=q(u)$ follows similarly.

    Thus, the result holds in the case $p(u)<q(u)$ and $p(w)<q(w)$.

    \medskip

    \noindent\textbf{Cases 2 \& 3: } $p(u)<q(u)$ and $p(w)=j(w)$ or ($p(u)<q(u)$ and $i(w)=q(w)$.

    \medskip
    
    Then  $j(u)-i(u)-2 = \dim(\Delta_u) = \dim(\Delta_w) = j(w)-i(w)-1$, so $j(u)-i(u)>j(w)-i(w)$ which is impossible as $[i(u),j(u)]\subseteq [i(w),j(w)]$.
    Thus, the result holds in these cases.

    \medskip
    
    \noindent\textbf{Case 4:} $p(u)=j(u)$ and $p(w)<q(w)$.

    \medskip
    
    In this case, $j(u)-i(u)-1 = \dim(\Delta_u) = \dim(\Delta_w) = j(w)-i(w)-2$.
    As \newline $[i(u),j(u)]\subseteq [i(w),j(w)]$, we get $[i(u),j(u)]= [i(w)+1,j(w)]$ or $[i(u),j(u)]= [i(w),j(w)-1]$.

    If $[i(u),j(u)]= [i(w)+1,j(w)]$, then the minimal generators of $I_{\Delta_w}$ are
    \[x_1,\dots,x_{i(u)-2},\qquad x_{j(u)+1},\dots,x_n,\qquad x_{[q(w),j(w)]},\qquad x_{[i(u)-1,p(w)]},\]
    while the minimal generators of $I_{\Delta_u}$ are 
    \[x_1,\dots,x_{i(u)-2},\qquad x_{j(u)+1},\dots,x_n,\qquad x_{[q(u),j(u)]},\qquad x_{i(u)-1}.\]
    It suffices to show that $x_{[q(w),j(w)]} = x_{[q(u),j(u)]}$.
    If those are not equal, then $[q(u),j(u)]$ and $i(u)-1$ must both be nonfaces of $\Delta_v$.
    Thus, $i(u)\leq i(v)$, and as $[i(u),j(u)]\subseteq [i(v),j(v)]$, $i(u)= i(v)$.
    Further, either $[q(u),j(u)] = [i(u),j(u)] = [i(v),p(v)]$ or $[q(u),j(u)] = [q(v),j(v)]$.
    If $[q(u),j(u)] = [q(v),j(v)]$, then all nonfaces of $\Delta_u$ are nonfaces of $\Delta_v$, i.e., $\Delta_v\subseteq \Delta_w$, which violates the antichain hypothesis.
    Thus, $[q(u),j(u)] = [i(u),j(u)] = [i(v),p(v)]$; if $j(u)=j(v)$, we again have a contradiction with the antichain hypothesis, so $j(v)>j(u)$.
    Thus, $i(w)+1=i(u)=q(u)=i(v)$ and $j(w) = j(u)<j(v)$.
    As $\Delta_u\subset \Delta_w$, $q(u)\geq q(w)$.
    If $q(u)>q(w)$, then $q(w)=i(w)< p(w)$, which is a contradiction.
    Thus, $q(u)=q(w)$ so $x_{[q(w),j(w)]}=x_{[q(u),j(u)]}$ as desired.

    If $[i(u),j(u)]= [i(w),j(w)-1]$, then as $p(w)<q(w)\leq j(w)-1 = j(u)=p(u)$ and $i(u)=i(w)$, $[i(w),p(w)]\subsetneq [i(u),j(u)]$, which is impossible unless $[q(u),j(u)]\subseteq [i(w),p(w)]$.
    Thus, $j(w)-1 = j(u)\leq p(w) <q(w) <j(w)$, which is impossible for integers.

    This concludes the case of $p(u)=j(u)$ and $p(w)<q(w)$.
    
    \medskip
    
    \noindent\textbf{Case 5:} $p(u)=j(u)$ and $p(w)=j(w)$. 

    \medskip
    
    In this case, $j(u)-i(u)-1 = \dim(\Delta_u) = \dim(\Delta_w) = j(w)-i(w)-1$.
    As $[i(u),j(u)] \subseteq [i(w),j(w)]$, we have $[i(u),j(u)] = [i(w),j(w)]$.
    Thus, the vertex nonfaces are the same for $\Delta_u$ and $\Delta_w$.
    They both only have one other minimal nonface, $x_{[q(u),j(u)]}$ and $x_{[q(w),j(w)]}$, respectively.
    If $q(w)\geq q(u)$, then $\Delta_w\subseteq \Delta_u$, which is impossible by the antichain assumption and the definition of $u$.
    Thus, as $q(w)<q(u)$, $x_{[q(w),j(w)]} = x_{[q(w),q(u)-1]}x_{[q(u),j(u)]}$ and so the result holds in this case.
    This concludes the case of $p(u)=j(u)$ and $p(w)=j(w)$.

    \medskip
    
    \noindent\textbf{Case 6:} $p(u)=j(u)$ and $i(w)=q(w)$. 
    
    \medskip
    
    By the same logic as the previous case, $[i(u),j(u)] = [i(w),j(w)]$.
    If $p(w)=j(w)$, then the previous case applies, so $p(w)<j(w)$.
    Then the vertex nonfaces of $\Delta_u$ and $\Delta_w$ agree, and their non-vertex minimal nonfaces are $[q(u),j(u)]$ and $[i(w),p(w)]$, respectively.
    But any nonface of $\Delta_w$ must be a nonface of $\Delta_u$ so $[i(u),p(w)] \subsetneq [i(u),j(u)]$ must be a nonface of $\Delta_u$ so we have a contradiction, so this case is impossible.

\medskip
 
 \noindent\textbf{Cases 7-9:}

\medskip

    The remaining cases of $i(u)=q(u)$ are symmetric to the previous three cases, so the result follows.
\end{proof}

Using the previous lemma, we can prove the following ideal containment.
\begin{lemma}\label{lem:colon-ideal-oscopy}
    With the same hypotheses and notation as \Cref{lem:mumbo-generator-jumbo}, letting \newline $G = (g_1,\dots, g_{r-1})$ be the corresponding ideal in $R_n$, we have
    \begin{equation}\label{eq:colon-dagger}
    (G:g_r)\subseteq G+(x_B).
    \end{equation}
\end{lemma}

\begin{proof}
 Let $\{g_1,\dots,g_r\}$ be the squarefree monomial generators of $I_{\Delta_v\cap\Delta_w}$ given by the minimal nonfaces as in \Cref{lem:mumbo-generator-jumbo} so that the generators of $I_{\Delta_w}$ are given by $\{g_1,\dots, g_{r-1},x_B g_r\}$. By the proof of \Cref{lem:mumbo-generator-jumbo}, we may order the $g_i$ so that $g_1,\dots, g_{r-2}$ are singletons.
 If $g_{r-1}$ is also a singleton, then modding out by $G$ reduces the desired statement to $\mathrm{Ann}_{R_{n-r+1}}(g_r) = (0:g_r) \subseteq (0)+x_B$, which is immediate by \Cref{lem:ann-squarefree-proofonly}.
 If $g_{r-1}$ is not a singleton, then by \Cref{lem:double pruned intersections} either:
\begin{itemize}
    \item $p(v)<p(w)$, $g_{r-1} = x_{[q(w),j(w)]}$, $g_r = x_{[i(w),p(v)]}$ and $B = [p(v)+1,p(w)]$, or
    \item $q(w)<q(v)$, $g_{r-1} = x_{[i(w),p(w)]}$, $g_r = x_{[q(v),j(w)]}$ and $B = [q(w),q(v)-1]$.
\end{itemize}
We shall prove the lemma in the first case. The proof of the second case is analogous and left to the reader. 

Quotient $R_n$ by the singleton generators $g_1,\dots, g_{r-2}$ to get a ring $R_m$, where $m= j(w)-i(w)+1$. For notational convenience, let us re-index the variables of $R_m$ to becomes $x_1,\dots, x_m$. After re-indexing, $g_r=x_{[1,p]}$, $B = [p+1,p']$  and $g_{r-1}=x_{[q, m]}$ (where $p'<q$). By 
\Cref{lem:local-endpoint-colon-proofonly}, we have
\[(g_{r-1}:g_r)\subseteq(g_{r-1},x_{[p+1,q-1]})\subseteq(g_{r-1},x_B), \]
where the second containment follows from $p'<q$. 
Lifting this containment back to $R_n$ proves \eqref{eq:colon-dagger}.
\end{proof}

We can use the previous lemma to prove the following lemma on the surjectivity of a map between Tor groups which will eventually let us show the ideal equality needed for \Cref{cond:B2}.

\begin{lemma}
\label{lem:one-sided-koszul-transfer}
Let $\Delta_v$ and $\Delta_w$ be (doubly-pruned) simplicial complexes coming from an antichain $\{v,w\}$ with $\Delta_w$ pure such that $\Delta_v\cap \Delta_w$ is a pure doubly-pruned complex of dimension $\dim(\Delta_w)$. 
Then the map
\[ \Tor_1^S(R_n,S/I_{\Delta_w})\longrightarrow\Tor_1^S(R_n,S/I_{\Delta_v\cap \Delta_w}) \]
is surjective.
\end{lemma}

\begin{proof}
The map in the statement of the lemma arises from tensoring the following short exact (Mayer-Vietoris) sequence by $R_n$:
    \[
    0\to S/I_{\Delta\cup\Delta_w}\to S/I_\Delta \oplus S/I_{\Delta_w}\to S/I_{\Delta\cap\Delta_w}\to 0.
    \]
    In particular, we get a long exact sequence on $\Tor$, part of which is
    \[
    \Tor_1^S(R_n,S/I_\Delta)\oplus \Tor_1^S(R_n,S/I_{\Delta_w})\to\Tor_1^S(R_n,S/I_{\Delta\cap\Delta_w})\to R_n/I_{\Delta\cup\Delta_w}\to R_n/I_\Delta \oplus R_n/I_{\Delta_w}.
    \]

As both $\Delta_w$ and $\Delta_w\cap \Delta_v$ are pure, their generators form a regular sequence, so the Tor groups can be computed using the Koszul complex of their respective regular sequences.
We will show surjectivity directly from this description.

Let $f_i$ be the standard coordinate basis of $(R_n)^r$. Using the same notation as in \Cref{lem:colon-ideal-oscopy}, let
\[z=\sum_{i\leq r}h_i f_i\]
be a degree-one cycle in $K_{R_n}(g_1,\ldots,g_{r-1},g_r)$, so that $\sum_{i\leq r}h_i g_i=0$ in $R_n$.
Modulo $G$, the cycle condition gives $\overline h_r \overline g_r=0$.
In the quotient, the containment in \Cref{lem:colon-ideal-oscopy} reduces to $(0:\overline{g_r}) \subseteq (\overline{x_B})$ so there exists $\overline q\in R_n/G$ such that $\overline{x_B} \,\overline{q}=\overline h_r$.
Lifting back to $R_n$, this gives
\[h_r-x_Bq=\sum_{i<r}b_i g_i \] for some $b_i$.

Define $y=\sum_{i<r}(h_i+b_i g_r)f_i+qf_r$.
Then $y$ is a cycle in $K_{R_n}(g_1,\dots, g_{r-1},x_Bg_r)$, because
\[ \sum_{i<r}(h_i+b_i g_r)g_i+qx_Bg_r = \sum_{i<r}h_ig_i+\left(\sum_{i<r}b_i g_i+qx_B\right)g_r=\sum_{i<r}h_i g_i+h_rg_r=0. \]
The image of $y$ in $K_{R_n}(g_1,\ldots,g_{r-1},g_r)$ is $\sum_{1\leq i< r} (h_i+b_ig_r)f_i + qx_Bf_r$, which differs from $z$ by 
\begin{align*}
\sum_{i<r} b_ig_rf_i + (qx_B-h_r)f_r &= \sum_{i<r} b_ig_rf_i - \sum_{i<r} b_ig_i
f_r\\
&= d\left(\sum_{i<r} b_i f_r\wedge f_i\right)
\end{align*}
a Koszul boundary.
Thus, $y$ maps to the equivalence class of $z$, so the map
\[ \Tor_1^S(R_n,S/I_{\Delta_w})\longrightarrow\Tor_1^S(R_n,S/I_{\Delta_v\cap \Delta_w}) \]
is surjective.
\end{proof}

\begin{lemma}\label{lem:deal-int-DDw}
    Let $\Delta,\Delta_w$ be simplicial complexes on $[n]$ such that $\Delta_w$ is doubly-pruned. If $\Delta\cap\Delta_w = \Delta_v\cap \Delta_w$ for some $v$ and if this intersection is pure of dimension $\dim(\Delta_w)$, then $(I_n+I_\Delta)\cap (I_n +  I_{\Delta_w}) = I_n + I_{\Delta \cup\Delta_w}$.
\end{lemma}

\begin{proof}
    We again start with the Mayer-Vietoris sequence and tensor by $R_n$ to get a long exact sequence in tor.
    Since $\Delta\cap\Delta_w = \Delta_v\cap \Delta_w$ satisfies the hypotheses of \Cref{lem:one-sided-koszul-transfer}, the map $\Tor_1^S(R_n,S/I_{\Delta_w})\to \Tor_1^S(R_n,S/I_{\Delta\cap\Delta_w})$ is surjective, hence the first map in the long exact sequence is surjective. Thus, $R_n/I_{\Delta\cup\Delta_w}\to R_n/I_\Delta \oplus R_n/I_{\Delta_w}$ is injective. Since the kernel of this map is $(I_n + I_\Delta)\cap (I_n+I_{\Delta_w})/(I_n + I_{\Delta\cup\Delta_w})$, then we have our desired equality of ideals.
\end{proof}

\begin{proof}[Proof of (3) implies (2)]
    Suppose $\mathcal{W} \subseteq S_n(\calP)$ and the elements of $\mathcal W$ admit an ordering $w_1,w_2,\dots, w_k$ satisfying the conditions in (3). We proceed by induction on $|\mathcal{W}|$. If $k=1$, then (2) is trivial. In the general case where $k>1$, then $w_1,\dots, w_{k-1}$ is an ordering that satisfies (3), so by induction
    \[
    \bigcap_{i=1}^{k-1} (I_n + I_{\Delta_{w_i}}) = I_n + \bigcap_{i=1}^{k-1}I_{\Delta_{w_i}}.
    \]
    Letting $\Delta = \bigcup_{i=1}^{k-1} \Delta_{w_i}$, then $\Delta\cap\Delta_{w_k} = \Delta_{w_t}\cap \Delta_{w_k}$ for some $t<k$ and is pure of dimension $\dim(\Delta_{w_k})$. Thus, $\Delta, \Delta_{w_k}$ satisfy the hypotheses of \Cref{lem:deal-int-DDw}, hence 
    \begin{align*}
    \bigcap_{w\in \mathcal{W}}(I_n + I_{\Delta_{w}}) &= \bigcap_{i=1}^{k-1} (I_n + I_{\Delta_{w_i}})\cap(I_n + I_{\Delta_{w_k}})\\
    &= (I_n + \bigcap_{i=1}^{k-1}I_{\Delta_{w_i}})\cap(I_n + I_{\Delta_{w_k}})\\
    &= (I_n + I_{\Delta})\cap(I_n + I_{\Delta_{w_k}})\\
    &= I_n + I_{\Delta\cup\Delta_{w_k}} = I_n + \bigcap_{w\in \mathcal{W}}I_{\Delta_w}\qedhere
    \end{align*}
\end{proof}

\subsection{The proof of \texorpdfstring{\Cref{thm:main-A}}{Theorem A}}
As we have proved \Cref{thm:main-B}, we can now complete the proof of \Cref{thm:main-A}.
\begin{proof}[Proof of \Cref{thm:main-A}]

    The equivalence of \Cref{cond:A1,cond:A2} is exactly the equivalence of \Cref{cond:B1,cond:B2} when $\mathcal{W} = \{w\}$ as then $J_\mathcal{W} = J_w$ and the equality in \Cref{cond:A2} is trivial.
    Then \Cref{lem: A3 implies A1} and \Cref{lem: A2 implies A3} give the remaining needed implications of \Cref{cond:A3} implies  \Cref{cond:A1} and \Cref{cond:A2} implies \Cref{cond:A3}, respectively. 
\end{proof}
\section{The combinatorics of doubly pruned simplicial complexes}
\label{sec: the combinatorics of doubly pruned simplicial complexes}
In this section we record the combinatorial and algebraic structure of the doubly pruned simplicial complexes themselves. We show that every such complex is vertex decomposable, hence shellable and sequentially Cohen-Macaulay. We also describe its facets explicitly, which gives a criterion for when $\Delta_w$ fails to be pure. We use this description in \Cref{subsec:hilb} to compute the Hilbert series of both $R_{\Delta_w}$ and $R_{\Delta_w,n}.$
 
For the rest of this section, write $\supp(\Delta) = \{v: \{v\}\in\Delta\}$ for the set of vertices of $\Delta.$ 
Recall the two subcomplexes that play a role in \Cref{def: vertex-decomposable}, that is $\lk_v(\Delta)$ and $\del_v(\Delta),$ the link and deletion of $v$ in $\Delta,$ respectively. 
Since $\lk_v(\Delta)\subseteq\del_v(\Delta),$ the shedding condition takes a simpler form.

\begin{remark}
    A vertex $v$ is a shedding vertex of $\Delta$ if and only if $G\cup\{v\}\not\in\Delta$ for every facet $G$ of $\del_v(\Delta).$ In particular, if $G$ is such a facet with $G\cup\{v\}\in\Delta,$ then $G\in\lk_v(\Delta),$ and any facet of $\lk_v(\Delta)$ containing $G$ also lies in $\del_v(\Delta)$ and hence equals $G$, so $G$ is a facet of both. The converse is immediate.  
\end{remark}

\begin{lemma}
\label{lem: two-min faces implies vertex-decomposable}
Let $\Delta$ be a simplicial complex on a finite set $V$ with at most two minimal nonfaces of cardinality at least two. Then $\Delta$ is vertex decomposable.
\end{lemma}
\begin{proof}
    We proceed by induction on the cardinality of the vertex set $V$ of $\Delta,$ with $\Delta = \varnothing$ or $\Delta$ being a simplex (including the singletons) being clear.
    If some minimal nonface is a singleton $\{u\},$ then no face of $\Delta$ contains $u,$ so $\Delta$ is also a simplicial complex on the smaller vertex set $V\setminus\{u\}$ with the same faces, and whether $\Delta$ is vertex decomposable does not depend on which of these vertex sets it is regarded. Thus, we may discard $u$ and so assume that no minimal nonface is a singleton.
    
    Now assume every minimal nonface has at least two elements. 
    Call them $A$ and $B$ (if a second minimal nonface exists). 
    If there are two minimal nonfaces with $A\cap B\neq\varnothing,$ choose $v\in A\cap B;$ otherwise choose $v\in A.$ 
    It follows that $v$ then lies in every minimal nonface unless $A$ and $B$ exist and are disjoint.
    If $v$ lives in every minimal nonface, then no subset of $V\setminus\{v\}$ contains a minimal nonface, so $\del_v(\Delta)$ is the simplex on $V\setminus\{v\},$ whose unique facet is $G = V\setminus\{v\}$ and $G\cup\{v\} = V\supseteq A$ is a nonface. 
    Otherwise, $A\cap B = \varnothing$ and $v\in A.$ 
    So, $\del_v(\Delta) = \{F\subseteq V\setminus\{v\}:B\not\subseteq F\},$ whose facets are the sets $G = V\setminus\{v,b\}$ with $b\in B.$  
    Since $b\not\in A,$ then $A\subseteq V\setminus\{b\} = G\cup\{v\}$ is a nonface. 
    In the first case, $\del_v(\Delta)$ is a simplex and has no minimal nonfaces at all, while in the second case $\del_v(\Delta) = \{F\subseteq V\setminus\{v\}:B\not\subseteq F\}$ has $B$ as its unique minimal nonface. So in both cases, $\del_v(\Delta)$ has at most one minimal nonface, and it is vertex decomposable by induction. Moreover, in both cases we have shown that $G\cup\{v\}\not\in\Delta$ for every facet $G$ of $\del_v(\Delta),$ so $v$ is a shedding vertex.
    Finally, $F\in\lk_v(\Delta)$ if and only if $F\cup\{v\}\in\Delta,$ i.e., if and only if $N\setminus\{v\}\not\subseteq F$ for every minimal nonface $N$ of $\Delta.$ So the minimal nonfaces of $\lk_v(\Delta)$ are the minimal elements of $\{A\setminus\{v\}, B\setminus\{v\}\}.$ Again, there are at most two, so $\lk_v(\Delta)$ is vertex-decomposable by induction. 
\end{proof}
\begin{theorem}
\label{thm: doubly pruned complexes imply vertex decomposable}
    Every doubly pruned simplicial complex is vertex decomposable. In particular, for every $w\in\sn(\calP),$ the complex $\Delta_w$ is vertex-decomposable, implying they are both shellable and sequentially Cohen-Macaulay.
\end{theorem}
\begin{proof}
    The minimal nonfaces of a doubly pruned complex of cardinality at least two are among $[i,p]$ and $[q,j],$ so \Cref{lem: two-min faces implies vertex-decomposable} applies. The remaining implications are immediate. 
\end{proof}
It is also possible to describe the facets of these complexes, which is what we do in \Cref{prop: facet description of doubly pruned complexes}. For a doubly pruned complex, normalize, so that $i<p$ and $q<j,$ which gives us the interval $[i,j].$ For $\Delta_w$ with $w\neq\id$ it is $[i(w),j(w)],$ since $p(w) = i(w)$ or $q(w) = j(w)$ would make $i(w)$ or $j(w)$ a fixed point of $w.$
\begin{proposition}
\label{prop: facet description of doubly pruned complexes}
Let $\Delta$ be a simplicial complex on vertex set $V$ with $W = \supp(\Delta)$ and with at most two minimal nonfaces $A,B$ of cardinality of at least two. 
Then $F$ is a facet of $\Delta$ if and only if $W\setminus F$ is a minimal subset of $W$ meeting each of $A$ and $B.$ Explicitly, when both $A$ and $B$ exist, the facets of $\Delta$ are 
\[
W\setminus\{v\}\quad (v\in A\cap B)\quad\text{and}\quad W\setminus\{a,b\}\quad (a\in A\setminus B,b\in B\setminus A).
\]
If only $A$ exists, the facets are $W\setminus\{v\}$ for $v\in A,$ and if neither exists, the unique facet is $W.$
Consequently, $\Delta$ is pure unless both $A$ and $B$ exist and $A\cap B\neq\varnothing,$ in which case $\Delta$ has facets of dimension $|W|-2$ and of dimension $|W|-3$ and $\dim(\Delta) = |W|-2.$
\end{proposition}
\begin{proof}
    Every minimal nonface is contained in $W$ or is a singleton outside of $W.$ 
    So for $F\subseteq W,$ by setting $T = W\setminus F,$ we have $F\in\Delta$ if and only if $T$ meets each of $A$ and $B.$ Complementation reverses the inclusion, which gives the first claim.
    A singleton $\{v\}$ meets both exactly when $v\in A\cap B,$ and singletons are automatically minimal. If $T$ is minimal with $|T|\geq 2,$ then $T\cap A\cap B = \varnothing,$ since one element of $T\cap A\cap B$ would already suffice. By choosing $a\in T\cap A$ and $b\in T\cap B,$ we obtain $a\in A\setminus B,b\in B\setminus A$ with $a\neq b,$ and $\{a,b\}\subseteq T$ meets both, so $T = \{a,b\}.$  Conversely, each such pair is minimal, as $\{a\}$ misses $B$ and \{b\} misses $A.$ 
    If only $A$ exists, then $T$ meets $A$ and $B$ exactly when $T\cap A\neq\varnothing,$ whose minimal solutions are the singletons $\{v\},$ $v\in A;$ if neither exists, $T = \varnothing$ is the unique minimal solution. Finally, $A$ and $B$ are incomparable, so whenever both exist, a facet of the second kind exists while a facet of the first kind exists exactly when $A\cap B\neq\varnothing.$
\end{proof}
\begin{corollary}
\label{cor: purity of doubly pruned complexes}
    Let $w\in\sn(\calP)$ with $w\neq\id$ and put $m = j(w) - i(w) + 1.$ Then $\Delta_w$ is not pure precisely when $i(w)<q(w)\leq p(w)<j(w).$ In this case, $\dim\Delta_w = m - 2$ and $\Delta_w$ has facets of dimension $m - 2$ and $m - 3.$
    If $q(w) = i(w)$ or $p(w) = j(w),$ then $\Delta_w$ is pure of dimension $m - 2.$ Otherwise, when $p(w) < q(w),$ it is pure of dimension $m - 3.$
\end{corollary}
\begin{proof}
    Here $W = [i(w),j(w)]$ and the minimal nonfaces of cardinality at least two are the minimal members of $\{[i(w),p(w)],[q(w),j(w)]\}.$ One of these intervals contains the other exactly when $q(w) = i(w)$ or $p(w) = j(w),$ leaving a single minimal nonface. Otherwise, both are minimal and meet in the interval $[q(w),p(w)]$ when $q(w)<p(w)$.
    To finish we apply \Cref{prop: facet description of doubly pruned complexes}.
\end{proof}

\begin{remark}
    \label{rem: doubly pruned purity}
    If a doubly pruned complex $\Delta$ is not $\Delta_w$ for some $w\in S_n(\calP)$, then it is pure if and only if $\Delta = \mathscr{P}([i,j])$, in which case it has dimension $j-i$. 
    If it is not pure, it has dimension $j-i-1$.
    In the notation of the corollary, these have dimensions $m-1$ and $m-2$, respectively.
\end{remark}

\begin{example}
\label{ex: 3412}
    Let $w = 3412 \in S_4(\calP)$, so that $w^{-1} = 3412$ as well. Its non-fixed
    points are $1,2,3,4$, so $i(w) = 1$ and $j(w) = 4$, giving $m = 4$, and
    \[
    p(w) = w^{-1}(i(w)) = w^{-1}(1) = 3,
    \qquad
    q(w) = w^{-1}(j(w)) = w^{-1}(4) = 2.
    \]
    By \Cref{def:Deltaw} the minimal nonfaces of $\Delta_w$ are the two intervals
    \[
    A = [i(w),p(w)] = \{1,2,3\}
    \qquad\text{and}\qquad
    B = [q(w),j(w)] = \{2,3,4\},
    \]
    with no singletons among them since $[i(w),j(w)] = [4]$. Here
    $A\cap B = \{2,3\}$, $A\setminus B = \{1\}$ and $B\setminus A = \{4\}$, so by
    \Cref{prop: facet description of doubly pruned complexes} the minimal subsets of
    $[4]$ meeting both $A$ and $B$ are $\{2\}$, $\{3\}$ and $\{1,4\}$, and the facets
    of $\Delta_w$ are their complements
    \[
    \{1,3,4\},\qquad \{1,2,4\},\qquad \{2,3\},
    \]
    of dimensions $2,2$ and $1$. Thus $\Delta_w$ is not pure and
    $\dim\Delta_w = 2$, in agreement with
    \Cref{cor: purity of doubly pruned complexes}, since
    $i(w) = 1 < 2 = q(w) \leq p(w) = 3 < 4 = j(w)$ and $m-2 = 2$. Pictorially,
     \begin{center}
    \begin{tikzpicture}[scale=0.4, vtx/.style={circle,fill=black,inner sep=1.7pt}]
      \coordinate (v1) at (0,2.6);
      \coordinate (v3) at (-2.2,-1);
      \coordinate (v2) at (2.2,-1);
      \coordinate (v4) at (0,0.45);
      \fill[faceblue] (v1) -- (v3) -- (v4) -- cycle;
      \fill[faceblue] (v1) -- (v2) -- (v4) -- cycle;
      \draw[thick] (v1) -- (v3) -- (v2) -- cycle;
      \draw[thick] (v1) -- (v4);
      \draw[thick] (v3) -- (v4);
      \draw[thick] (v2) -- (v4);
      \node[vtx] at (v1) {}; \node[vtx] at (v2) {};
      \node[vtx] at (v3) {}; \node[vtx] at (v4) {};
      \node[above]                     at (v1) {\small $1$};
      \node[below right]               at (v2) {\small $2$};
      \node[below left]                at (v3) {\small $3$};
      \node[above right,inner sep=2pt] at (v4) {\small $4$};
      \node[anchor=east] at (-2.75,0.8) {\large $\Delta_w \;=$};
    \end{tikzpicture}
    \end{center}
    Following the proof of
    \Cref{lem: two-min faces implies vertex-decomposable}, the shedding vertex is any $v \in A\cap B$, say $v = 2$. Then $\del_2(\Delta_w)$ is the simplex on $\{1,3,4\}$, while $\lk_2(\Delta_w)$ has minimal nonfaces $A\setminus\{2\} = \{1,3\}$ and $B\setminus\{2\} = \{3,4\}$ and hence facets $\{1,4\}$ and $\{3\}$. The unique facet $\{1,3,4\}$ of $\del_2(\Delta_w)$ satisfies $\{1,3,4\}\cup\{2\} = [4]\supseteq A$, a nonface, confirming that $2$ is a shedding vertex.
\end{example}

\begin{lemma}
\label{lem:pure-intersection}
    Let $u,v\in\sn(\calP).$ If $\Omega = \Delta_u\cap\Delta_v$ is pure of dimension $\dim\Delta_u,$ then $\Delta_u$ is pure.
\end{lemma}
\begin{proof}
    By \Cref{prop:intersection is doubly pruned}, $\Omega$ is doubly pruned. Set $W = \supp(\Delta_u) = [i,j]$, $W' = \supp(\Omega) = [i',j']$, $m = |W|$, and $m' = |W'|.$ Since $\Omega\subseteq\Delta_u,$ then $W'\subseteq W$ and $m'\leq m.$

    Suppose $\Delta_u$ is not pure, and write $A = [i,p]$ and $B = [q,j]$ for its two minimal nonfaces of cardinality at least two. Then $u\neq\id$, so \Cref{cor: purity of doubly pruned complexes} gives $i<q\leq p<j$ and $\dim\Delta_u = m-2.$ By \Cref{prop: facet description of doubly pruned complexes}, the pure complex $\Omega$ has dimension $m'-1,m'-2,$ or $m'-3,$ since it has zero, one, or two minimal nonfaces of cardinality at least two. Since $m' - 3 = m - 2$ would give $m' = m + 1 > m,$ this case is impossible, so only two cases remain.

    If $\dim(\Omega) = m'-2$, then it has one minimal nonface $C$ of cardinality at least 2 and $m' = m,$ so $C\subseteq A$ and $C\subseteq B,$ which implies that $C\subseteq [q,p].$ However, $C$ has left endpoint $i$ or right endpoint $j,$ so $i\in[q,p]$ or $j\in[q,p].$ This forces $q\leq i$ or $j\leq p$ and contradicts $i<q$ and $p<j.$ If $\Omega$ is a simplex on $W'$ and $m' = m - 1,$ then $W' = [i+1,j]$ or $W' = [i,j-1],$ and the minimal nonfaces of $\Omega$ are the singletons outside $W'.$ But both $A$ and $B$ must contain one of them, which forces $q = i$ in the first case and $p = j$ in the second, once again forcing a contradiction.
\end{proof}

\subsection{The Hilbert series for the Stanley-Reisner rings of doubly pruned simplicial complexes}
\label{subsec:hilb}

Let $S = \bK[x_1,\ldots,x_n]$, let $w\in\sn(\calP)$, and let $I_{\Delta_w}$ be the Stanley-Reisner ideal of $\Delta_w.$ 
In this section, we give formulations of the Hilbert series for both $R_{\Delta_w} = S/I_{\Delta_w}$ and of $R_{\Delta_w,n} = S/(I_{\Delta_w}+I_n),$ where $I_n$ is generated by the elementary symmetric functions.
Recall the definition of $\Delta_w$ for $w\in S_n(\calP)$ from Definition~\ref{def:Deltaw}. As in earlier sections, we write $i = i(w), j = j(w), p = p(w),$ and $q = q(w).$ We also set $A := [i,p]$ and $B := [q,j],$ two interval nonfaces of $\Delta_w.$
When $A$ and $B$ meet, we write $C := A\cap B$ and $\gamma = := \#A + \#B - m,$ so that $\#C = \#A\cap B.$
The Stanley-Reisner ideal is then
\[
I_{\Delta_w} = (x_1,\ldots,x_{i-1},x_{j+1},\ldots,x_n,x_{A},x_{B}),
\]
where $x_{[c,d]} = x_cx_{c+1}\cdots x_d.$
We now define three parameters used throughout this section:
\begin{itemize}
    \item $m = j - i + 1$, the size of the support of $\Delta$,
    \item $\#A = p - i + 1 = \#[i,p],$
    \item $\#B = j - q + 1 = \#[q,j].$
\end{itemize}
The two interval nonfaces are disjoint precisely when $p < q$ (equivalently $\#A + \#B \leq m$) they overlap when $q \leq p$ (equivalently $\#A + \#B > m$), in which case their overlap has size $\#A + \#B - m.$ 
If neither interval contains the other, then in the overlapping case we have $A\cup B = [i,j],$ i.e. the union already used the realized vertices. One interval contains the other exactly when $q = i$ or $p = j,$ equivalently when $\#B = m$ or $\#A = m;$ by \Cref{cor: purity of doubly pruned complexes} this is exactly the case in which $\Delta_w$ has a single minimal nonface of cardinality at least two, and we call it the \emph{nested case.} In the nested case, we write $[c,d]$ for that nonface, so $[c,d] = A$ if $\#B = m$ and $[c,d] = B$ if $\#A = m,$ and we set $k = \#[c,d] = \min(\#A,\#B).$ The nested case is not exotic as it contains every simple transposition $s_r,$ for which $[i,j] = [r,r+1]$ and $p = j,$ and it is the majority of $\sn(\calP)$ for small $n.$ It is, however, the easy case, since then $R_{\Delta_w} = S'/(x_{[c,d]})$ is a hypersurface, with $S' = \bK[x_i,\ldots,x_j]$ as in $\Cref{lem: minimal free resolution of R_Delta}$ below, and we treat it separately in each statement of this section. Both Hilbert series computed in this section rest on a single resolution, which we record in the next lemma.
\begin{lemma}
\label{lem: minimal free resolution of R_Delta}
    Suppose that neither of $A,B$ contain the other. 
    Let 
    \[
    S' = \bK[x_i,\ldots,x_j]\cong S/(x_1,\ldots,x_{i-1},x_{j+1},\ldots,x_n) \]
    and set $f = x_{A},g = x_{B}\in S',$ with degrees $\#A,\#B,$ respectively. Writing $f' = \lcm(f,g)/f$ and $g'= \lcm(f,g)/g,$ 
    \begin{align}
        0\longrightarrow S'(-min(\#A + \#B, m))
        \xlongrightarrow{
            \begin{psmallmatrix}
            f' \\ -g' 
            \end{psmallmatrix}
        } 
        S'(-\#A)\oplus S'(-\#B)
        \xlongrightarrow{
            \begin{psmallmatrix}
            f & g
            \end{psmallmatrix}
        } 
        S'
        \longrightarrow S'/(f,g)
        \longrightarrow 0
    \end{align}
    is a minimal free resolution of $S'/(f,g) = R_{\Delta_w}.$
\end{lemma}

\begin{proof}
    We begin by noting that since neither $A$ nor $B$ contains the other, then $\lcm(f,g) = x_{A\cup B}$ has degree $\min(\#A+\#B,m).$
    Now write $h = \gcd(f,g),$ so $f = hg'$ and $g = hf'$ with $\gcd(f',g') = 1.$
    Since $S'$ is a UFD, then if $\alpha f + \beta g = 0,$ then $\alpha g' = -\beta f'.$ 
    Furthermore, since $f',g'$ are coprime, then $\alpha = cf'$ and $\beta = -cg'$ for some $c\in S'.$
    Thus, the syzygy module of $(f,g)$ is free of rank one and generated by $(f',-g')$ and the sequence is exact at the middle term. Furthermore, since $S'$ is a domain, then we have injectivity at the left-hand map.
    The sequence is minimal: since $A\not\subseteq B$ and $B\not\subseteq A,$ then both $f'$ and $g'$ have positive degree, so every entry of both matrices lies in the graded maximal ideal. 
    Finally, if the two intervals overlap then neither is nested in the other, so together they reach endpoints $i,j$ and $A\cup B = [i,j],$ an interval of size $m.$ 
    However, if they are disjoint then $\#(A\cup B) = \#A + \#B.$ In both cases, $\#(A\cup B) = \min(\#A + \#B, m).$
\end{proof}

\begin{theorem}
\label{thm: hilbert series - single Schubert variety}
Let $w \in \sn(\calP)$, and let $m, \#A,$ and $\#B$ be as previously defined. Then the Stanley-Reisner ring $R_{\Delta_w} = S/I_{\Delta_w}$ has Hilbert series
\[
  \hilb\big(R_{\Delta_w},\, t\big)
  = \frac{1 - t^{\#A} - t^{\#B} + t^{\min(\#A + \#B,\, m)}}{(1-t)^{m}}
\]
Equivalently:
\begin{itemize}
\item if the interval nonfaces are disjoint \textup{(}$\#A + \#B \leq m$, i.e.\ $p < q$\textup{)},
\[
  \hilb\big(R_{\Delta_w},t\big) = \frac{(1-t^{\#A})(1-t^{\#B})}{(1-t)^{m}};
\]
\item if they overlap \textup{(}$\#A + \#B > m$, i.e.\ $p \leq q$\textup{)},
\[
  \hilb\big(R_{\Delta_w},t\big) = \frac{1 - t^{\#A} - t^{\#B} + t^{m}}{(1-t)^{m}};
\]
\item in particular, in the nested case ($\#A = m$ or $\#B = m$), with $k = \min(\#A,\#B),$
\[
\hilb\big(R_{\Delta_w},t) = \frac{1 - t^k}{(1-t)^m}.
\]
\end{itemize}
For $w = \mathrm{id}$ we have $\Delta_w = \{\emptyset\}$, so that
$I_{\Delta_w} = (x_1,\dots,x_n)$, $R_{\Delta_w} = \bK$, and
$\hilb(R_{\Delta_w},t) = 1$.
\end{theorem}

\begin{proof}
    In the nested case $R_{\Delta_w} = S'/(x_{[c,d]})$ is a hypersurface of degree $k$ in $m$ variables, so $\hilb\big(R_{\Delta_w},t) = (1 - t^k)/(1-t)^m,$ since either $\#A = m $ or $\#B = m.$
    Now suppose neither interval contains the other.
    Taking the alternating sum of the Hilbert series along the minimal free resolution described in \Cref{lem: minimal free resolution of R_Delta} and using $\hilb(S',t) = (1-t)^{-m}$ gives the displayed formula. In the disjoint case $\min(\#A + \#B,m) = \#A + \#B$ and $1 - t^{\#A} - t^{\#B} + t^{\#A + \#B} = (1-t^{\#A})(1-t^{\#B}),$ which gives us (a). 
    We note that this would also be clear without appealing to the minimal free resolution. 
    Since $f,g$ share no variables, then $g$ is a nonzerodivizor on $S'/(f),$ so $(f,g)$ forms a regular sequence, meaning $S'/(f,g)$ is a complete intersection, giving us the result.
    Now, in the overlapping case, $\min(\#A + \#B, m) = m,$ giving (b). The case for $w = \id$ is immediate.
\end{proof}

The resolution in \Cref{lem: minimal free resolution of R_Delta} also allows us to compute the depth and dimension of $R_{\Delta_w}.$
\begin{corollary}
\label{cor: len-depth-dim R_Delta}
    Let $w\neq\id\in\sn(\calP).$ The minimal free resolution of $R_{\Delta_w}$ over $S$ has 
    \begin{itemize}
        \item length $n - m +2$ (length $n - m + 1$ in the nested case)
        \item $\mathrm{depth}\,R_{\Delta_w} = m - 2$  ($\mathrm{depth}\,R_{\Delta_w} = m - 1$ in the nested case) 
        \item $\dim R_{\Delta_w} = 
        \begin{cases}
        m - 2 & p < q \\
        m - 1 & q\leq p.
        \end{cases}
        $
    \end{itemize}
    Consequently, when $p < q$ or in the nested case, $R_{\Delta_w}$ is Gorenstein. When $q\leq p$ and neither interval contains the other, then it is one short of being Cohen-Macaulay. 
\end{corollary}
\begin{proof}
    In the nested case, $R_{\Delta_w} = S'/(x_{[c,d]})$ and the minimal free resolution over $S$ is the Koszul complex on the regular sequence $x_1,\ldots,x_{i-1},x_{j+1},\ldots,x_n,x_{c,d]},$ of length $n-m+1;$ so $R_{\Delta_w}$ is a complete intersection, hence Cohen-Macaulay and Gorenstein, with $\depth R_{\Delta_w} = \dim R_{\Delta_w} = n - (n - m + 1) = m - 1.$ 
    Now suppose neither interval contains the other.
    The length is $n - m +2,$ since the minimal free resolution over $S,$ formed by tensoring the minimal free resolution from \Cref{lem: minimal free resolution of R_Delta} with a Koszul complex has $n-m$ generators $x_1,\ldots,x_{i-1},x_{j+1},\ldots,x_n.$ 
    The Auslander-Buchsbaum formula gives us that the depth is $n - (n - m + 2) = m - 2,$ and $\dim R_{\Delta_w} = \dim\Delta_w+1$ is given by \Cref{cor: purity of doubly pruned complexes}. 
    Finally, it is Gorenstein when $p < q$ since it is both Cohen-Macaulay ($\mathrm{depth}R_{\Delta_w} = \dim R_{\Delta_w}$) and a complete intersection ring. 
    When $q\leq p,$ $\mathrm{depth}R_{\Delta_w} = \dim R_{\Delta_w}-1.$
\end{proof}

As of this point, we have computed everything for the Stanley-Reisner ring of $\Delta_{w},$ $w\in\sn(\calP).$ We mostly do this for completeness. However, the ring that is more relevant to this paper is $R_{\Delta_w,n} = S/(I_\Delta+I_n).$ These are not simply a choice forced by the Borel presentation of the cohomology ring $H^*(G/B),$ but instead are a natural choice to reduce any Stanley-Reisner ring. Consider the image of the elementary symmetric functions generating $I_n$ under $S/I_{\Delta_w}$ -- these are exactly the \emph{universal system of parameters}. This system of parameters have made appearances in throughout the literature, but were initially named by Herzog and Moradi, and simultaneously defined by Adams and Reiner.
\begin{definition}[\cite{herzog-moradi,adams-reiner}]
\label{def: universal system of parameters}
    For a simplicial complex $\Delta$ of dimension $d-1,$ the \emph{universal system of parameters} $\theta = (\theta_1,\ldots,\theta_d)$ of $R_\Delta$ are 
    \[
    \theta_k = \sum_{\substack{F\in\Delta \\ \#F = k}}x_F\in R_\Delta,\quad k = 1,\ldots,d.
    \]
    They form a homogeneous system of parameters.
\end{definition}
Since $R_{\Delta_w,n} = S/(I_\Delta+I_n) = R_{\Delta_w}/(\theta),$ we are able to use this presentation to compute its Hilbert series. Since $R_{\Delta_w,n}$ is isomorphic to cohomology ring of the Schubert variety $X_w$ and \Cref{lem:defined by inclusions}, we translate the Hilbert series of this ring into a sum over the Bruhat interval below $w:$ since by \Cref{prop: schuberts form a basis for Rn} and since $J_w$ is spanned by Schubert polynomials \Cref{eq: Jw} the ring $R_n/J_w$ has a $\bK$-basis the images of the Schubert polynomials $\fraks_v$ with $v\leq w,$ and $\deg\fraks_v = \ell(v),$ we have
\[
\hilb(R_{\Delta_w,n},t) = \sum_{v\leq w}t^{\ell(v)}.
\]
Note that $\hilb(\bK[\theta_1,\ldots,\theta_d],t) = \prod_{r=1}^d(1-t^r)^{-1}.$
\begin{theorem}
\label{thm: hilbert series for R_Delta_n}
    Let $w\in\sn(\calP),$ $w\neq\id,$ with $m, \#A, \#B$ as above, with overlap $\gamma = \#A + \#B - m.$
    \begin{enumerate}
        \item If $p < q$ (disjoint, $\gamma\leq 0$), then $\theta_1,\ldots,\theta_{m-2}$ is a regular sequence on $R_{\Delta_{w}},$ and 
        \[
        \hilb(R_{\Delta_{w,n}},t) = [\#A]_t[\#B]_t[m-2]_t!.
        \]
        \item If $q\leq p$ and neither interval contains the other (overlapping, $i < q\leq p < j,$ $\gamma\geq 1$), then $\theta_1,\ldots,\theta_{m-1}$ is not a regular sequence on $R_{\Delta_w}$ and 
        \[
        \hilb(R_{\Delta_{w,n}},t) 
        = [\gamma]_t[m-1]_t! + t^\gamma~[m-\#B]_t[m-\#A]_t[m-2]_t!.
        \]
        \item If $q = i$ or $p = j$ (nested, $\#A = m$ or $\#B = m$), then $\theta_1,\ldots,\theta_{m-1}$ is a regular sequence on $R_{\Delta_w}$ and, with $ k = \min(\#A,\#B),$
        \[
        \hilb(R_{\Delta_{w,n}},t) = [k]_t[m - 1]_t!. 
        \]
    \end{enumerate}
\end{theorem}
\begin{proof}
    \begin{enumerate}
        \item By \Cref{cor: len-depth-dim R_Delta}, $R_{\Delta_w}$ is Cohen-Macaulay of Krull dimension $m - 2.$ 
        Thus it's a homogeneous system of parameters $\theta_1,\ldots,\theta_{m-2}$ is a regular sequence and 
        \begin{align*}
            \hilb(R_{\Delta_{w,n}},t) 
                & = \hilb(R_{\Delta_{w}},t)\prod_{r=1}^{m-2}(1-t^r) \\
                & = \frac{(1-t^{\#A})(1-t^{\#B})}{(1-t)^m}\cdot(1-t)^{m-2}[m-2]_t! \\
                & = [\#A]_t[\#B]_t[m-2]_t!.
        \end{align*}
        In the string of equalities, we use the standard equalities
        \[
        1-t^p = (1-t)[\#A]_t,
        \quad
        1-t^q = (1-t)[\#B]_t,
        \quad
        \prod_{r=1}^{m-2}(1-t^r) = (1-t)^{m-2}[m-2]_t!.
        \]
        
        \item In this case, $\dim R_{\Delta_w} = m - 1$ while $\depth\, R_{\Delta_w} = m - 2$ by \Cref{cor: len-depth-dim R_Delta}, so there is no homogeneous system of parameters of length $m - 1$ that also forms a regular sequence. 
        For the Hilbert series, recall that in the overlapping case $i<q\leq p<j,$ since neither nonface contains the other. Our strategy is as follows. We find an ideal $Z\subseteq R_{\Delta_w}$ such that both $Z$ and the quotient $R_{\Delta_w}/Z$ are, up to a shift, face rings of complexes whose Hilbert series modulo $\theta$ we already know. We also require of $Z$ that $R_{\Delta_w}/Z$ is Cohen-Macaulay, hence a free $\bK[\theta]$-module. The short exact sequence $0\to Z\to R_{\Delta_w}\to R_{\Delta_w}/Z\to 0$ then remains exact after passing modulo $\theta,$ and $\hilb(R_{\Delta_w,n},t)$ is the sum of the Hilbert series $Z/\theta Z$ and of $(R_{\Delta_w}/Z)/\theta(R_{\Delta_w}/Z).$
        Write 
        \[
        C := A\cap B=[q,p],
        \quad\#C = \gamma,
        \]
        and cover $\Delta_w$ by two subcomplexes, the faces that do not contain $C$ and the star of $C:$
        \[
        \Gamma_C := \{G\in\Delta_w:C\not\subseteq G\},
        \quad
        \star(C) := \{G:G\cup C\in\Delta_w\}
        \]
        so that $\Delta_w = \Gamma_C\cup\star(C)$ (the two subcomplexes are not disjoint; they meet in the faces $G\in\star(C)$ with $C\not\subseteq G).$ 
        Let $Z := (x_C)\subseteq R_{\Delta_w}$ be the ideal generated by the span of the monomials whose support contains $C,$ so that
        \[
        0 \longrightarrow 
        Z \longrightarrow 
        R_{\Delta_w}\longrightarrow 
        R_{\Delta_w}/Z \longrightarrow 
        0.
        \]
        First, we observe that both the quotient $R_{\Delta_w}/Z$ and the ideal $Z$ correspond to complexes that we have already seen in this paper.
        Beginning with the quotient, we notice that a monomial dies when it is divisible by $x_C,$ that is, when its support contains $C,$ thus $R_{\Delta_w}/Z = R_{\Gamma_C}.$     
        Since $C$ is contained in both $A$ and $B$, a subset of $[i,j]$ not containing $C$ contains neither $A$ nor $B$ and hence is a face $\Delta_w,$ so $\Gamma_C = \{G\subseteq[i,j]:C\not\subseteq G\}.
        $ and the minimal nonfaces of $\Gamma_C$ are $C$ together with the singletons outside $[i,j],$ with one minimal nonface of of cardinality of at least two.
        By \Cref{lem: two-min faces implies vertex-decomposable}, $\Gamma_C$ is vertex-decomposable, hence shellable. Moreover, by \Cref{prop: facet description of doubly pruned complexes} its facets are $[i,j]\setminus\{v\}$ with $v\in C,$ so it is pure of dimension $m - 2.$
        A pure shellable complex is Cohen-Macaulay, so $\theta_1,\ldots,\theta_{m-1}$ is a regular sequence on $R_{\Gamma_C}.$ 
        Thus, $R_{\Gamma_C} = S'/(x_C)$ is a hypersurface, where $S' = \bK[x_i,\ldots,x_j]$ as in \Cref{lem: minimal free resolution of R_Delta}, so $\hilb(R_{\Gamma_C},t) = (1-t^\gamma)/(1-t)^m.$ Since quotienting by a regular sequence of elements in degrees $1,\ldots, m - 1$ multiplies the Hilbert series by $\prod_{r = 1}^{m-1}(1-t^r),$ we obtain 
        \[
        \hilb(R_{\Gamma_C}/\theta R_{\Gamma_C},t) = 
        \frac{1-t^\gamma}{(1-t)^m}
        \prod_{r = 1}^{m - 1}(1-t^r)
        = [\gamma]_t[m - 1]_t!.
        \]

        We now move onto identifying the ideal $Z = (x_C).$
        The monomials of $Z$ are the $x_C x^{\vec{a}}$ with $C\cup\supp(\vec{a})\in\Delta_w,$ that is, with $\supp(\vec{a})\in\star(C).$ 
        Thus, division by $x_C$ is an isomorphism $Z\cong R_{\star(C)}(-\gamma)$ of graded vector spaces. 
        Moreover, it is also an isomorphism of modules. 
        For a face $F$ of $\Delta_w$ with $\#F = k$ and $g = x_C x^{\vec{a}}\in Z,$ the product $x_Fg$ survives exactly when $F\cup C\cup\supp(\vec{a})\in\Delta_w,$ so the terms of $\theta_k$ with $F\not\in\star(C)$ never survive. 
        Furthermore, for $F\in\star(C),$ the condition is simply the multiplication rule of $R_{\star(C)}.$ 
        Thus, for $G\subseteq[i,j]$ we have $G\cup C\in\Delta_w$ if and only if $[i,q-1]\not\subseteq G$ and $[p+1,j]\not\subseteq G.$ 
        That is, $\star(C)$ is the doubly pruned complex with support $[i,j]$ whose interval nonfaces are $[i,q-1]$ and $[p+1,j]$ of sizes $m - \#B$ and $m - \#A.$
        These are disjoint, since $q - 1 < p + 1,$ so the first case outlined in this theorem and already proved above can be applied to $\star(C).$ 
        This gives that
        \[\hilb(R_{\star(C),n},t) = [m-\#B]_t[m-\#A]_t[m-2]_t!.\]
        In particular, $\star(C)$ has Krull dimension $m - 2,$ so it has no faces of size $m - 1,$ and therefore $\theta_{m-1}$ acts on $Z$ as the empty sum $\theta_{m-1}(\star(C)) = 0.$ 
        Hence, $(\theta_1,\ldots,\theta_{m-1})Z = (\theta_1,\ldots,\theta_{m-2})Z,$ where the second sequence is the universal system of parameters of $R_{\star(C)}.$ 
        Therefore, $Z/\theta Z\cong (R_{\star(C)}/(\theta_1,\ldots,\theta_{m-2}))(-\gamma) = R_{\star(C),n}(-\gamma),$ and so
        \[
        \hilb(Z/\theta Z,t) = t^\gamma~[m-\#B]_t[m-\#A]_t[m-2]_t!.
        \]

        Finally, remembering that $R_{\Gamma_C}$ is Cohen-Macaulay and so a free $\bK[\theta]$-module, we 
        combine these Hilbert series calculations and identifications in order to reduce the earlier short exact sequence to another short exact sequence by considering it modulo $\theta:$
        \[
        0 \longrightarrow
        Z/\theta Z \longrightarrow
        R_{\Delta_w,n} \longrightarrow
        R_{\Gamma_C}/\theta R_{\Gamma_C} \longrightarrow
        0.
        \]
        Since Hilbert series are additive along short exact sequences, we obtain the required formula.

        \item By \Cref{cor: len-depth-dim R_Delta}, $R_{\Delta_w} = S'/(x_{[c,d]})$ is Cohen-Macaulay of Krull dimension $m - 1,$ and $\theta_r = 0$ for $r\geq m$ since $\Delta_w$ has no face of dimension $m - 1.$ Thus, its homogeneous system of parameters $\theta_1,\ldots,\theta_{m-1}$ is a regular sequence, and by \Cref{thm: hilbert series - single Schubert variety} 
        \[
        \hilb(R_{\Delta_{w,n}},t) = \hilb(R_{\Delta_w},t)\prod_{r=1}^{m-1}(1-t^r) = \frac{1-t^k}{(1-t)^m}\cdot(1-t)^{m-1}[m-1]_t! = [k]_t[m-1]_t!.
        \]
    \end{enumerate}
\end{proof}

Note that the formula of Part (2) of \Cref{thm: hilbert series for R_Delta_n} is valid in the nested case as well since there is one of $[m-\#A]_t,[m-\#B]_t$ is $[0]_t = 0$ and $\gamma = k.$ However, the argument of Part (2), is not, which is why we treat the proof for the nested case separately. For instance, when $w = s_r$ we have $m = k = 2$ and $\hilb(R_{\Delta_{w,n}},t) = [2]_t[1]_t! = 1 + t,$ the Poincar{\'e} polynomial of $X_{s_r} \cong \mathbb{P}^1,$ while for $w = w_0,$ the complex $\Delta_w$ is the boundary of the simplex on $[n], k = m = n,$ and $[n]_t[n-1]_t! = [n]_t!.$

\begin{example}
    Let $w = 45312\in\mathrm{S}_5.$ Then $w(3) = 3$ is the only fixed point. Then $i = 1, j = 5, p = w^{-1}(1) = 4,$ and $q = w^{-1}(5) = 2,$ so $A = [1,4]$ and $B = [2,5].$ 
    Since we are in the overlapping case: $q\leq p,$ then $ m = 5,\#A = \#B = 4,$ and $o = \#A + \#B - m = 3.$ \Cref{thm: hilbert series for R_Delta_n} gives
    \begin{align*}
        \hilb(R_{\Delta_{w,n}},t) 
        & = [3]_t[4]_t! + t^3[1]_t[1]_t[3]_t! \\
        & = 1 + 4t + 9t^2 + 15t^3 + 18t^4 + 16t^5 + 10t^6 + 4t^7 + t^8.
    \end{align*}
\end{example}

\subsection{The Adams-Reiner conjecture for doubly pruned simplical complexes}

Utilizing the universal system of parameters (\Cref{def: universal system of parameters}) Adams and Reiner posed a conjecture \cite[Conj.~6.1]{adams-reiner} that has only been proven for small families of simplicial complexes, but is open in general. It says that the shape of the minimal free resolution of $R_\Delta$ over the parameter ring $\bK[\theta]$ is governed by the barycentric subdivision $\Sd\Delta,$ colored so that the barycenter of a face $F$ of a simplicial complex $\Delta$ of dimension $d$ has color $\# F$ in $\Sd\Delta.$ Through their ``colorful" Hochster formula \cite[Thm.~3.3]{adams-reiner} and under their grading specialization $\epsilon_k\mapsto k,$ where $\epsilon_k$ is the $k$-th standard basis vector in $\bN^d$ it reads
\begin{equation}
\label{eq: Adams-Reiner conj}
    \dim_\bK
    \Tor_r^{\bK[\theta]}(R_\Delta,\bK)_{N} =
    \sum_{\substack{
    S\subseteq [d] \\
    \wt S = N
    }}
    \dim_\bK \Tilde{H}^{\# S - r - 1}(\Sd_S(\Delta);\bK).
\end{equation}
Here $d = \dim\Delta + 1,$ the weight $\wt S: = \sum_{k\in S}k$ records the internal degree, and $\Sd_S\Delta$ is the \emph{rank-selected} subcomplex of the barycentric subdivision $\Sd\Delta,$ (i.e., the order complex of $\{F\in\Delta:F\neq\varnothing,\# F\in S\}$).
We prove it for the first class where $R_\Delta$ fails to be free over $\bK[\theta]$ of simplicial complexes, namely those with at most two minimal nonfaces of cardinality of at least two, which by definition includes every doubly pruned complex, in particular, every $\Delta_w$. 
All (co)homology is reduced and taken with coefficients in $\bK.$ Since the complexes involved are finite and $\bK$ is a field, $\Tilde{H}^i(\delta;\bK)$ is dual to $\Tilde{H}_i(\Delta;\bK)$ by the universal coefficient theorem, so the two have the same dimension. Since \Cref{eq: Adams-Reiner conj} is only concerned with dimensions, we compute its RHS in homology, where the Mayer-Vietoris sequence is available.

We use only use three known facts within this section. First, \Cref{eq: Adams-Reiner conj} holds when $R_\Delta$ is Cohen-Macaulay \cite[Cor.~6.7]{adams-reiner} since both sides vanish for $r\geq 1$ giving for every Cohen-Macaulay complex $\Delta$ with $d = \dim\Delta + 1,$
\begin{equation}
\label{eq: conjecture in the CM case}
    \sum_{T\subseteq[d]}t^{\wt T}\dim\Tilde{H}_{\#T-1}(\Sd_T(\Delta); \mathbb{K}) = \hilb(R_{\Delta,n},t).
\end{equation}

Second, for such a Cohen-Macaulay complex $\Delta$ and $T\subseteq[d]$ the complex $\Sd_T(\Delta)$ is again Cohen-Macaulay of dimension $\# T - 1$ \cite{munkres,stanley-balanced}. In particular, $\Tilde{H}_i(\Sd_T(\Delta)) = 0$ for $i\neq\# T - 1.$ Third, the two sides of \Cref{eq: Adams-Reiner conj} have the same alternating sum over $r,$ both equal to $\hilb(R_\Delta,t)\prod_{r=1}^d(1-t^r)$ \cite[Prop~2.5,~Eq.~30]{adams-reiner}.

As stated before, a complex is a simplex, then its face ring is simply a polynomial ring. If the complex has exactly one minimal nonface, then its face ring is a hypersurface. If the complex has two disjoint minimal nonfaces, then its face ring is a complete intersection ring, and therefore Cohen-Macaulay. Thus, for the rest of this subsection, we can assume that the complex has two minimal nonfaces, each of cardinality of at least two, that are not disjoint (so they are overlapping). Call these sets $A$ and $B$ and call the overlap $C := A\cap B\neq\varnothing.$ Write $W = \supp(\Delta), m = \#W, \gamma = \#C, \alpha = \#(A\setminus B),$ and $\beta = \#(B\setminus A),$ all positive. Since $W\setminus\{v\}$ is a face exactly when $v\in C,$ then we have $\dim\Delta = m - 2$ and $d = m - 1.$ \Cref{rmk: conj into Delta w notation} translates this back into $\Delta_w$ notation.

The proof of \Cref{thm: hilbert series for R_Delta_n}(2) applies verbatim in this generality since it uses only that $A$ and $B$ are incomparable minimal nonfaces meeting in $C.$ 
The complexes $\Gamma_C,$ $\star(C),$ and the ideal $Z = (x_C)$ are defined there, and $\Delta = \Gamma_C\cup\star(C).$

Within $W = \supp(\Delta),$ the only minimal nonface of $\Gamma_C$ is $C.$
The complex $\star(C)$ has the two minimal nonfaces $A\setminus B$ and $B\setminus A;$ since $I_{\Gamma_{C}\cap\star(C)} = I_{\Gamma_C}+I_{\star(C)}.$ Finally, the complex $\Gamma_C\cap\star(C)$ has all three. 
Within each of these three lists, the sets are pairwise disjoint, so their monomials form a regular sequence and $R_{\Gamma_C}, R_{\star(C)},$ and $R_{\Gamma_C\cap\star(C)}$ are complete intersections, Cohen-Macaulay of Krull dimensions $d,d-1,$ and $d - 2,$ with
\begin{equation}
\label{eq: hilb of star and C intersect star}
    \hilb(R_{\star(C),n}, t) = [\alpha]_t[\beta]_t[m-2]_t!,
    \quad
    \hilb(R_{\Gamma_C\cap\star(C),n}, t) = [\gamma]_t[\alpha]_t[\beta]_t[m-3]_t!
\end{equation}

\begin{proposition}
\label{prop: dimensions formulas for reduced homologies}
    Let $S\subseteq [d]$ be nonempty, set $s = \#S$ and $S' = S\setminus\{d\},$ and $T = S\cap[d-2].$ Then $\Tilde{H}_i(\Sd_S(\Delta)) = 0$ for $i < s-2,$ and also for $i = s - 2$ when $d\not\in S.$ 
    \begin{itemize}
        \item If $d,d-1\in S,$ then
        \begin{align*}
        \dim\Tilde{H}_{s-2}\left(\Sd_S(\Delta)\right) 
        & = \dim\Tilde{H}_{s-2}\left(\Sd_{S'}(\star(C))\right) \\
        & + \dim\Tilde{H}_{s-3}\left(\Sd_{S'\setminus\{d-1\}}(\Gamma_C\cap\star(C))\right), 
    \end{align*}
    \item If $d\in S$ and $d - 1\not\in S,$ then
    \begin{align*}
        \dim\Tilde{H}_{s-2}\left(\Sd_S(\Delta)\right) 
        & = \dim\Tilde{H}_{s-2}\left(\Sd_{S'}(\star(C))\right) \\
        & - \dim\Tilde{H}_{s-2}\left(\Sd_{S'}(\Gamma_C\cap\star(C))\right).
    \end{align*}
    \end{itemize}
\end{proposition}
\begin{proof}
    A chain of faces of $\Delta$ with largest face $M$ lies in $\Gamma_C$ when $C\not\subseteq M$ and in $\star(C)$ when $C\subseteq M.$ 
    Thus, $\Sd_S(\Gamma_C)\cup\Sd_S(\star(C))$ with intersection $\Sd_S(\Gamma_C\cap\star(C)),$ and so we can apply Mayer-Vietoris.
    Moreover, $\star(C)$ has no face of size $d$ and $\Gamma_C\cap\star(C)$ has no face of size $d - 1,$ so $\Sd_S(\star(C)) = \Sd_{S\cap[d-1]}(\star(C))$ and $\Sd_S(\Gamma_C\cap\star(C)) = \Sd_T(\Gamma_C\cap\star(C)).$

    By the second fact listed a priori, the homology of the three complexes is thus concentrated in the single degrees $s-1,$ $\#(S\cap[d-1])-1,$ and $\#T-1,$  each between $s - 3$ and $s - 1.$
    For $i < s - 2,$ both groups on either side of $\Tilde{H}_i(\Sd_S(\Delta))$ vanish, which proves our first claim.
    If $d\not\in S,$ the degrees are $s - 1,$ $s - 1,$ and $s - 2$ or $s - 1.$ So, the groups on either side of $\Tilde{H}_{s - 2}(\Sd_S(\Delta))$ vanish as well.

    Let $d\in S$ so that $\# S' - 1 = s - 2.$
    \begin{itemize}
        \item If $d - 1\in S,$ then the three degrees are $s - 1, s - 2,$ and $s - 3$ so the sequence in degree $s - 2$ is
        \[
        0 
        \longrightarrow
        \Tilde{H}_{s-2}\left(\Sd_{S'}(\star(C))\right)
        \longrightarrow
        \Tilde{H}_{s-2}\left(\Sd_S(\Delta)\right)
        \longrightarrow
        \Tilde{H}_{s-3}\left(\Sd_T(\Gamma_C\cap\star(C))\right)
        \longrightarrow
        0.
        \]
        \item If $d-1\not\in S,$ then $S\cap[d-1] = T = S',$ and the degrees are $s-1,s-2,$ and $s-2,$ so the sequence in degree $s-2$ is
        \[
        \Tilde{H}_{s-2}\left(\Sd_T(\Gamma_C\cap\star(C))\right)
        \xlongrightarrow{\phi}
        \Tilde{H}_{s-2}\left(\Sd_{S'}(\star(C))\right)
        \longrightarrow
        \Tilde{H}_{s-2}\left(\Sd_S(\Delta)\right)
        \longrightarrow
        0,
        \]
        with $\phi$ induced by inclusion. 
        Every simplex of $\Sd_{S'}(\star(C))$ is a chain with cardinalities in $S'$ and so has dimension of at most $s - 2.$
        Thus, every $(s-1)$-simplex of $\Sd_S(\Delta)$ already lies in $\Sd_S(\Gamma_C).$
        Since neither complex has a simplex of dimension $s,$ then $\Tilde{H}_{s-1}(\Sd_S(\Gamma_C))\to\Tilde{H}_{s-1}\left(\Sd_S(\Delta)\right)$ is surjective. The connecting map out of $\Tilde{H}_{s-1}(\Sd_S(\Delta))$ is therefore zero, which means $\phi$ is injective, and $\Tilde{H}_{s-2}(\Sd_S(\Delta)) = \mathrm{coker}\phi.$
    \end{itemize}
\end{proof}

\begin{theorem}
\label{thm: adams-reiner conj for doubly pruned complexes}
    The Adams-Reiner conjecture (restated in \Cref{eq: Adams-Reiner conj}) holds over any field for every simplicial complex with at most two minimal nonfaces of cardinality at least two. In particular, for every doubly pruned complex, and so for every $\Delta_w$ with $w\in\sn(\calP).$ Moreover, $\Tor_r^{\bK[\theta]}(R_\Delta,\bK) = 0$ for $r\geq 2,$ and $\Tor_1$ vanishes unless the two nonfaces meet. In the latter case,
    \[
    \sum_N
    \dim_\bK \Tor_1^{\bK[\theta]}(R_\Delta,\bK)_N t^N 
    = t^{m-1+\gamma}~[\alpha]_t[\beta]_t[m-2]_t!.
    \]
\end{theorem}
\begin{proof}
    Call $P_r$ and $Q_r$ the generating functions in $t^N$ of the LHS and the RHS of \Cref{eq: Adams-Reiner conj}, respectively. We have completed the Cohen-Macaulay cases above, so we progress by assuming the remaining hypothesis. In this proof, we compute the LHS and the RHS, we show all terms on both sides are $0$ for $r\geq 2,$ and then compare the remaining generating functions. 
    
    \textbf{LHS:} As $\Gamma_C$ is Cohen-Macaulay of Krull dimension $d,$ it is $\bK[\theta]$-free, so the long exact sequence induced on $0\rightarrow Z\rightarrow R_\Delta\rightarrow R_{\Gamma_C}\rightarrow 0$ gives that $\Tor_r^{\bK[\theta]}(R_\Delta,\bK)\cong\Tor_r^{\bK[\theta]}(Z,\bK)$ for $r\geq 1,$ where $Z = (x_C)\cong R_{\star(C)}(-\gamma)$ with $\theta'$ acting as the universal system of parameters of $\star(C)$ and $\theta_d$ as zero. 
    Since $R_{\star(C)}$ is $\bK[\theta']$-free, the Koszul complex on $(\theta',\theta_d)$ has acyclic rows and zero differential in the $\theta$ direction. Thus, $\Tor_{\geq 2}(Z,\bK) = 0$ and $\Tor_1(Z,\bK)\cong (R_{\star(C)}/(\theta'))(-\gamma-d).$
    Hence, $P_r = 0$ for $r\geq 2,$ and by \Cref{eq: hilb of star and C intersect star}
    \begin{equation}
        P_1(t) = t^{\gamma+d}~\hilb(R_{\star(C),n},t) = t^{m-1+\gamma}~[\alpha]_t[\beta]_t[m-2]_t!\neq 0.
    \end{equation}

    \textbf{RHS:} By \Cref{prop: dimensions formulas for reduced homologies}, $Q_r = 0$ for $r\geq 2.$ Furthermore, we know that the only those $S$ containing $d$ contribute to $Q_1.$ Write such an $S$ as $S'\sqcup\{d\}$ with $S'\subseteq[d-1]$ so that $t^{\wt S} = t^dt^{\wt S'}.$
    Summing the first term of $\Cref{prop: dimensions formulas for reduced homologies}$ over these $S$ gives a $t^d$ term times a sum over all $S'\subseteq [d-1].$ In the second term, those $S'$ that contain $d-1$ become those subsets $T=S'\setminus\{d-1\}$ of $[d-2],$ which carry a factor of $t^{d-1},$ with a positive sign, while those avoiding $d-1$ are already such a $T$ and carry a negative sign. Together, they contribute the factor $t^{d-1}-1.$ Hence,
    \begin{align*}
        Q_1(t) = 
        & t^d~\sum_{S'\subseteq[d-1]}t^{\wt S'}
        \dim\Tilde{H}_{\#S'-1}\left(\Sd_{S'}(\star(C))\right) \\
        & + t^d(t^{d-1}-1)\sum_{T\subseteq [d-2]}t^{\wt T}
        \dim\Tilde{H}_{\#T-1}\left(\Sd_T(\Gamma_C\cap\star(C))\right).
    \end{align*}
    Both sums are instances of \Cref{eq: conjecture in the CM case}, so by \Cref{eq: hilb of star and C intersect star} they can be computed explicitly. 
    Using $1 - t^{d-1} = (1-t)[m-2]_t,$ then $[m-2]_t[m-3]_t! = [m-2]_t!$ and $(1-t)[\gamma]_t = 1 - t^\gamma.$ So,
    \[
    Q_1(t) = t^d~[\alpha]_t[\beta]_t[m-2]_t!
    \left(1-(1-t)[\gamma]_t\right)
    = t^{d+\gamma}~[\alpha]_t[\beta]_t[m-2]_t!,
    \]
    which is precisely $P_1(t).$ 
    Thus, $P_1 = Q_1,$ and since $P_r = Q_r = 0$ for $r\geq 2$ the third fact from above forces $P_0 = Q_0,$ as well.
\end{proof}

\begin{remark}
\label{rmk: conj into Delta w notation}
Let $w\in\sn(\calP),$ $w\neq\id,$ give $\Delta_w$ that sits in the overlapping case where $i < q \leq p < j,$ i.e., $q\leq p$ and neither interval contains the other, so that $A$ and $B$ are the two minimal nonfaces of $\Delta_w$ of cardinality of at least two. 
Discard any singleton nonfaces, which do not affect either side of \Cref{eq: Adams-Reiner conj}. 
Then $\Delta_w$ satisfies the hypothesis with $A = [i,p], B = [q,j], C = [q,p], W = [i,j],$ and $m = j - i + 1.$ 
Since incomparability forces $i < q\leq p < j,$ then
\[
\alpha = \#[i,q-1] = m - \#B,
\quad
\gamma = \#[q,p],
\quad
\beta = \#[p+1,j] = m - \#A.
\]
Thus, \Cref{thm: adams-reiner conj for doubly pruned complexes} proves \Cref{eq: Adams-Reiner conj} for $\Delta_w.$ Its row $r = 0$ is Part (2) of \Cref{thm: hilbert series for R_Delta_n} can now be recovered from the rank-selected homology of $\Sd\Delta_w$ rather than from the ideal $(x_C).$ Its row $r=1$ is new and is
\[
\sum_N
\dim_\bK\Tor_1^{\bK[\theta]}(R_{\Delta_w},\bK)_N~t^N
=
t^{m-1+\gamma}~[m-\#B]_t[m-\#A]_t[m-2]_t!
\]
and has total dimension $(m-\#B)(m-\#A)(m-2)!.$ 
By our calculation of $P_1(t)$ in the previous proof, it measures the failure of $R_{\Delta_w}$ to be free over $\bK[\theta].$
\end{remark}

\section{An example: The poset of simplicial complexes for \texorpdfstring{$w\in\mathrm{S}_4(\calP)$}{w avoiding P}}
\label{sec: example}
\input{Example-S4}

\newpage
\bibliographystyle{plain}
\bibliography{bibliography}
\end{document}